\documentstyle [12pt]{article} 
\def\picture#1by#2(#3){ 
\vbox to #2 { 
  \hrule width #1 height 0pt depth 0pt \vfill \special{picture #3}} 
} 
 
\def\scaledpicture#1by#2(#3scaled#4){{ 
\dimen0=#1  \dimen1=#2 
\divide\dimen0 by 1000 \multiply\dimen0 by #4 
\divide\dimen1 by 1000 \multiply\dimen1 by #4 
\picture \dimen0 by \dimen1 (#3 scaled #4)}} 
\def\dfigure#1by#2(#3scaled#4offset#5:#6) 
  {\medskip 
   \vglue 2mm minus 2mm 
   $$ 
     \hbox{ 
       \hglue#5 
       {\scaledpicture #1 by #2 (#3 scaled #4)} 
     } 
   $$ 
   \par\goodbreak 
   \vglue 2mm minus 2mm 
   \medskip}

\def\qmod#1#2{{\hbox{}^{\displaystyle{#1}}}\!\big/\!\hbox{}_{
\displaystyle{#2}}}

\newfam\msbfam

\font\twelmsb=msbm10 at 12pt
\font\tenmsb=msbm10
\font\sevenmsb=msbm10 at 7pt
\font\fivemsb=msbm10 at 5pt

\textfont\msbfam=\twelmsb
\scriptfont\msbfam=\sevenmsb
\scriptscriptfont\msbfam=\fivemsb

\def\Bbb{\fam\msbfam\twelmsb}

\def\C{{\Bbb C}}

\def\G{{\Bbb G}}
\def\H{{\Bbb H}}

\def\N{{\Bbb N}}

\def\R{{\Bbb R}}
\def\Z{{\Bbb Z}}

\def\union{\mathop{\bigcup}}
\def\qed {\hfill\vrule height6pt width6pt depth0pt \bigskip}

\def\map{\longrightarrow}
\def\textmap#1{\mathop{\vbox{\ialign{
                                ##\crcr
    ${\scriptstyle\hfil\;\;#1\;\;\hfil}$\crcr
    \noalign{\kern-1pt\nointerlineskip}
    \rightarrowfill\crcr}}\;}}

\def\textlmap#1{\mathop{\vbox{\ialign{
                                ##\crcr
    ${\scriptstyle\hfil\;\;#1\;\;\hfil}$\crcr
    \noalign{\kern-1pt\nointerlineskip}
    \leftarrowfill\crcr}}\;}}

\newfam\meuffam
\font\twelvmeuf=eufm10 at 12 pt
\font\tenmeuf=eufm10
\font\sevenmeuf=eufm7

\textfont\meuffam=\twelvmeuf
\scriptfont\meuffam=\tenmeuf
\scriptscriptfont\meuffam=\sevenmeuf

\def\germ{\fam\meuffam\tenmeuf}

\def\eg{{\germ e}}

\def\g{{\germ g}}

\def\mg{{\germ m}}

\def\pg{{\germ p}}
\def\qg{{\germ q}}

\begin{document}
\def\Pr{{\rm Pr}}
\def\tr{{\rm Tr}}
\def\End{{\rm End}}
\def\Spin{{\rm Spin}}
\def\U{{\rm U}}
\def\SU{{\rm SU}}
\def\SO{{\rm SO}}
\def\PU{{\rm PU}}
\def\GL{{\rm GL}}
\def\spin{{\rm spin}}
\def\u{{\rm u}}
\def\su{{\rm su}}
\def\so{{\rm so}}
\def\ub{\underbar}
\def\pu{{\rm pu}}
\def\Pic{{\rm Pic}}
\def\Iso{{\rm Iso}}
\def\NS{{\rm NS}}
\def\deg{{\rm deg}}
\def\Hom{{\rm Hom}}
\def\h{{\germ h}}
\def\Herm{{\rm Herm}}
\def\Vol{{\rm Vol}}
\def\pf{{\bf Proof: }}
\def\id{{\rm id}}
\def\i{{\germ i}}
\def\im{{\rm im}}
\def\rk{{\rm rk}}
\def\ad{{\rm ad}}
\def\h{{\bf H}}
\def\coker{{\rm coker}}
\def\dv{\bar\partial}
\def\Ad{{\rm Ad}}
\def\RSU{\R SU}
\def\ad{{\rm ad}}
\def\dva{\bar\partial_A}
\def\da{\partial_A}
\def\p{\partial\bar\partial}
\def\sp{\Sigma^{+}}
\def\sm{\Sigma^{-}}
\def\spm{\Sigma^{\pm}}
\def\smp{\Sigma^{\mp}}
\def\oo{{\scriptstyle{\cal O}}}
\def\ooo{{\scriptscriptstyle{\cal O}}}
\def\sw{Seiberg-Witten }
\def\pa{\partial_A\bar\partial_A} 
\def\Dr{{\raisebox{0.15ex}{$\not$}}{\hskip -1pt {D}}}
\def\gr{{\scriptscriptstyle|}\hskip -4pt{\g}}
\def\subsetint{{\  {\subset}\hskip -2.45mm{\raisebox{.28ex}
{$\scriptscriptstyle\subset$}}\ }}
\def\nr{\parallel}
\newtheorem{sz}{Satz}[section]
\newtheorem{thry}[sz]{Theorem} 
\newtheorem{pr}[sz]{Proposition}  
\newtheorem{re}[sz]{Remark} 
\newtheorem{co}[sz]{Corollary} 
\newtheorem{dt}[sz]{Definition}
\newtheorem{lm}[sz]{Lemma} 
\newtheorem{cl}[sz]{Claim}

\title{Moduli spaces of $PU(2)$-monopoles}
\date{ }
\author{Andrei Teleman}
\maketitle

\section{Introduction}

The
most natural way to prove the equivalence between Donaldson theory and Seiberg-Witten
theory is to consider a suitable moduli space of "non-abelian monopoles". In [OT5]  it was
shown that an
$S^1$-quotient of a moduli space of quaternionic monopoles should give an homological
equivalence between a fibration over  a union of Seiberg-Witten moduli spaces  and  a fibration
over   certain $Spin^c$-moduli spaces  [PT1]. 

By the same method, but using moduli spaces of $PU(2)$-monopoles instead of quaternionic
monopoles,  one should be able to express any Donaldson   invariant
  in terms of Seiberg-Witten invariants associated with the \underbar{twisted} abelian 
monopole equations of [OT6]. In  [T1], [T2], we have shown that this idea can be
further generalized
  to express   Donaldson-type  invariants associated with higher symmetry
groups in terms of new Seiberg-Witten-type invariants. 

The strategy has a very general algebraic-geometric analogon, which we call the
"Master Space" strategy. This procedure, developed by Ch.
Okonek and the author     [OT7], [OST] reduces the problem of the
computation of certain numerical invariants of a GIT moduli space to similar computations
on simpler    moduli spaces. One "couples" the given GIT problem to a  simpler one (having the
same symmetry group), and then studies the "Master Space"  associated with the coupling as a
$\C^*$-space. The fixed  point locus of the $\C^*$-action consists of the original moduli
space and a union of  simpler  ones. Then one can use the $S^1$-quotient of the master space
to define a  homological equivalence between a projective fibration  over the initial moduli
space and a projective fibration over the other components of the fixed point locus. In the
GIT-framework, as in the gauge theoretical one,  the
technical difficulty is the same: the master space can be singular. The present paper deals
with this difficulty in the gauge theoretical situation.\vspace{2mm}\\
\ \ \ A  program for proving the equivalence between Donaldson theory and
Seiberg-Witten theory, which also uses moduli
spaces of non-abelian monopoles, is due to Pidstrigach and Tyurin [PT2], and was already
announced by Pidstrigach in a Conference  at the Newton  Institute in Cambridge, in December
1994.

There are, however, several important differences between 
Pidstrigach-Tyurin's original
approach, and the strategy developed  by  Ch.
Okonek in collaboration with  the author,   which is the strategy we follow in the present
paper.

First, our equations have a gauge group of the form $SU(E)$ and hence the moduli spaces which
we construct are $S^1$-spaces; in contrast, the Pidstrigach-Tyurin equations  [PT2] have a
gauge group of the form $U(E)$. Whereas we fix the connection in
the determinant line bundle, they only fix the curvature of this connection. 
If $H_1(X,\Z)=0$,
their moduli space is the $S^1$-quotient of ours.   On the other hand, the 
$S^1$-operation
plays a very important role in our strategy: The description of the ends around the
 abelian locus
at infinity uses in an essential way the $S^1$-equivariance of the local models.

Second, we do not follow   Pidstrigach-Tyurin's program to  prove generic
 regularity results.
We show (see section 3.1) that the  proofs of the transversality theorems 
which they use
[PT2] to get generic regularity  are incomplete, by indicating counterexamples
 to one of the
statements on which these proofs are based\footnote{ This gap    as
well as the difficulty of the problem was pointed out by the author   during the
Workshop "4-dimensional manifolds", Oberwolfach, March 1996.   }.

It is interesting to notice
that, in fact, {\it any}  non-abelian solution of  the equations in the
K\"ahler case gives a counterexample to their statement.  This same
statement was also used by the authors in the their definition of the
$Spin^c$-polynomial invariants [PT1], on which was based their
approach to prove the Van de Ven conjecture. 

The transversality problem is very complicated, for the $PU(2)$-monopole
equations as well as for the  non-abelian $Spin^c$-equations. The difficulty is
 the same in
both cases: in the non-abelian points with degenerate spinor component 
transversality cannot
be proved using only perturbations with 0-order operators.  

In [T1] the author tried to
use perturbations with first order operators, and proved that perturbations
 of this type lead
to transversality at least away from the solutions which are abelian on a 
non-empty open
set. However, in order to have a  
complete transversality result away from the abelian locus, one would need a
unique continuation theorem which seems to be  difficult to get because of the 
perturbed
symbol.

Another way to achieve transversality is to use an
infinite family of   "holonomy perturbations" [FL].

The present paper solves two fundamental problems concerning the moduli
 spaces of
$PU(2)$-monopoles: {\it generic regularity} and {\it compactification}.

First we prove an $S^1$-equivariant generic-smoothness theorem:  we
define perturbations of the equations which lead to $S^1$-spaces which,
for generic choices of the perturbing parameters    are smooth, at least
outside the "Donaldson locus"  (the vanishing locus of the projection on the
spinor component) and of the abelian locus (Theorem 3.19).   The proof of 
the generic-smoothness theorem is \ub{not} a pure transversality
argument; it combines a standard transversality argument with a new
method to control  the exceptions to transversality.

{\it Our result shows that one does get  regularity for a generic 
choice of a system
$(g,\sigma,\beta,K)$, consisting of a metric, a compatible 
$Spin^{U(2)}(4)$-structure
$\sigma$ and an order 0- perturbation $(\beta,K)$ of the type 
considered in [PT2].}

We also obtain   generic regularity results for the normal bundles of
the Donaldson locus and the abelian locus within the moduli space
(Theorem 3.21, Proposition 3.22).   Similar results, but obtained using
quite different lines of reasoning, were obtained by Feehan [F] in a
preprint   distributed around the same time as the first version of the
present paper.

Therefore one can go forward towards a proof of the Witten conjecture 
(see for instance [OT5] for a detailed description of the strategy) using 
relatively simple equations.

Note however that the generic
regularity results which we prove for the ASD-$Spin^c$-equations, do
\ub{not} automatically solve the transversality problem needed in order
to give sense  to the $Spin^c$-polynomial invariants, and to
use them effectively.  For this purpose one would need  a pure
transversality argument for the ASD-$Spin^c$-equations.\footnote
{In order to have {\it
well defined} invariants, one needs a {\it smooth} parameterized
moduli space ( [DK], p 143, 149).  
Moreover,   the K\"ahlerian parameters  are {\it all}
\ub{non}-generic in our sense; on the other hand, \ub{all} 
computations needed in order to get a proof of the Van de Ven conjecture
using $Spin^c$-invariants, must be done in the K\"ahler case. } 
This seems to be difficult. 
The 
theory of $Spin^c$-polynomial invariants, and the  attempt
to prove the Van de Ven conjecture
using these invariants, should be
 therefore
revised.

We get our result in two steps. In a first step we prove that, using only  
the perturbations
$(\beta,K)$, one
can prove the following partial transversality result: If the 
Seiberg-Witten map extended to the parameterized moduli space is not a
submersion in a point $(A,\Psi,\beta,K)$, then the spinor component
$\Psi$ must be degenerate. This is very easy to see.

In the second step we prove that, if we also let the $Spin^{U(2)}(4)$-structure
 (together with
the metric) vary , then the moduli space 
$\widetilde{{\cal D}{\cal M}}^*_X$ of
solutions with non-trivial but degenerate spinor component in the
enlarged parameterized moduli space $\widetilde{\cal M}^*_X$ has
infinite codimension in every non-abelian point. Using this, we can
show (by "weakening" locally the degeneracy equation) that every
non-abelian   point
$[p]$ in $\widetilde{{\cal D}{\cal M}}^*_X$ has a  neighbourhood
$U_{[p]}$ which is a closed analytic subspace of a manifold $V_{[p]}$
which is Fredholm of negative index over the enlarged parameter space. 
Taking a countable subcover $(V_{[p]_i})_{i\in\N}$, and using the fact
that Fredholm maps are locally proper ([Sm]),  we prove that the set of
parameters for which there exists a non-abelian solution with
non-trivial degenerate spinor component is of the first category.  
The
desired set of "generic parameters" is then obtained by intersecting   the
complement of this set  with the set of regular values of the projection of
$\widetilde{{\cal M}}^*_X\setminus\widetilde{{\cal D}{\cal M}}^*_X$
on the parameter space.

We believe that this method is in fact a very general one; it can be
summarized as follows: Prove first a partial transversality result
 using  
perturbations with 0-order differential operators, and show then 
that the space of
solutions which are exceptions to transversality has infinite 
codimension  if one
introduces new variable parameters. Such a result is to be 
 expected
 provided the  "exceptional
solutions" , the ones which are exceptions to transversality, 
solve an overdetermined elliptic
system.

In particular, the method can be applied to obtain generic 
regularity  along the Donaldson
and the abelian locus. More precisely,  the moduli space of 
solutions (with non-vanishing
spinor component) of the Dirac-ASD system of [PT1] becomes 
smooth of expected
dimension  for generic perturbations. The same property has the
 complement of the
zero-section in the fibration of "normal infinitesimal deformations" 
over the subspace of
abelian solutions associated with an abelian reduction of the 
$Spin^{U(2)}(4)$-bundle.

In this way we obtain   perturbed moduli spaces which are   smooth
 except in the 
abelian points and in the Donaldson-points. These points remain   exceptions to 
transversality, and in general, regularity (smoothness and expected dimension) {\it
cannot} be achieved in   these points by using $S^1$-equivariant perturbations.

The second purpose of the paper, the existence of an  "Uhlenbeck
compactification" for the perturbed moduli spaces, is achieved in section
4 (see Theorem 4.24).  A different proof of the "Uhlenbeck
compactification" can be found in [FL].  

Our arguments follow the same strategy as in the instanton case [DK],
which can be summarized as follows:\\
\\
Local estimates --  Regularity --  Removable Singularities --  Compactification.
\\
\\
 Some care must be taken, since the monopole equations are  only "scale
invariant", not conformal invariant as in the instanton case. On the other
hand, many of the results in [DK] were obtained by cutting off the
solutions and transferring the problem from the 4-ball to the 4-sphere,
and then using the conformal invariance of the equations.  

 Our proof uses
the same method, but endows the sphere with a metric with non-negative
sectional curvature which is flat in a neighbourhood of the north pole.
With this choice, the
 corresponding first order elliptic operators ($\Dr$, $d^*+d^+$,  \dots )
 are still injective.   
For the local computations  we work with pairs whose connection
component is in Coulomb gauge in the sense of [DK], so that  all the results
in [DK] about connections in Coulomb gauge apply automatically. 
Therefore, we do not use the Coulomb gauge condition for pairs  which
follows from the elliptic complex of the $PU(2)$-monopole equations
(compare with [FL]).

A short version of our proof of the Uhlenbeck compactification appeared in
[OT5], and a very detailed version of it can be found   in [T1]. The
existence of an Uhlenbeck compactification for   moduli
spaces of non-abelian monopoles was predicted by Pidstrigach and
Tyurin in [PT2].

Note that in order to prove the equivalence between the  Donaldson and the
Seiberg-Witten theories, it now remains only to give explicit descriptions
of the ends of the moduli space along the abelian locus, and to calculate the
corresponding contributions.   

My own strategy to study the ends of the moduli spaces of
$PU(2)$-monopoles is based on the analytical results in [T3].  The
$PU(2)$-monopole equations are not conformally invariant, so it
is difficult   to use the method developed in the case of  instantons [DK]
(which consists of
 identifying the  solutions concentrated in a  point with
the solutions on the connected sum of $X$ with $S^4$). We use a new
strategy  [T4] which is still based on the gluing method. We
obtain concentrated solutions by gluing (non-concentrated) solutions on
$X$ corresponding to lower topological data, with  concentrated instantons
on the  {\it tangent spaces}, and then we deform the  obtained
almost-solutions
 into solutions.  This last step makes use of the classical   Fredholm
$L^p$ theory on
$X$, as well as  of  the  Fredholm $L^p$-theory on the  tangent spaces
(instead of
$S^4$) which is developed in the quoted paper.

Progress on this problem, using different methods, was also announced by 
Feehan and Leness.

I would like to thank professor Ch. Okonek for encouraging me to write this
paper, for the careful reading,   and for his suggestions. I would also
like to thank  professor S. T. Yau for suggesting me to submit the paper
to AJM. Finally I thank  the referee for the very careful checking of the
technical arguments and for his valuable observations.

\section{$PU(2)$-monopoles}
\subsection{The $Spin^{U(2)}$ group and $Spin^{U(2)}$-structures}

For a more detailed  presentation of the theory of 
$Spin^{U(2)}$-structures we refer the interested reader to 
[T1], [T2].
In these papers we also introduce the concept of $Spin^G$-structures and
$G$-monopole equations for quite general compact Lie groups $G$.
 
The group $Spin^{U(2)}$ is defined by 
$$Spin^{U(2)}:=Spin\times_{\Z_2} U(2)\ .
$$
Using the natural isomorphism 
$\qmod{U(2)}{\Z_2}\simeq PU(2)\times S^1$, 
we get the exact
sequences
$$1\map Spin\map Spin^{U(2)}\textmap{(\bar\delta,\det)}
PU(2)\times S^1\map 1 
$$
$$1\map U(2)\map Spin^{U(2)}\stackrel{\pi}{\map}  SO\map 1
\eqno{(1)}$$
$$1\map \Z_2\map Spin^{U(2)}\textmap{(\pi,\bar\delta,\det)}  
SO\times PU(2)\times S^1
 \map 1 \ .
$$

Let $X$ be a compact manifold and $P^u$  a   $Spin^{U(2)}$-bundle  over
  $X$. 
We consider the
following associated bundles
$$\pi(P^u):=P^u\times_\pi SO, \ \bar\delta(P^u):=
P^u\times_{\bar\delta} PU(2),\
\det(P^u):=P^u\times_{\det} S^1,
$$
$$\G_0:=P^u\times_{Ad\circ\bar\delta} SU(2);\
\gr_0:=P^u\times_{ad\circ\bar\delta} SU(2)
\ ,
$$
where
$Ad:PU(2)\map Aut(SU(2))$,
$ad:PU(2)\map so(su(2))$ are induced by the adjoint morphism $SU(2)\map Aut(SU(2))$,
$SU(2)\map so(su(2))$.

The group of sections ${\cal G}_0:=\Gamma(X,\G_0)$ can be identified 
with the group of
automorphisms of $P^u$ over $\pi(P^u)\times_X\det(P^u)$. After 
suitable Sobolev completions it
becomes  a Lie group, whose Lie algebra is the corresponding completion of 
$A^0(\gr_0)$.\\ 

Let $P$ be a $SO$ bundle over $X$. A $Spin^{U(2)}$-structure in $P$ is a
 morphism
 $P^u \map P$
of type $\pi$, where $P^u$ is a $Spin^{U(2)}$-bundle. Two 
$Spin^{U(2)}$-structures
 $P^u\map P$,
$P'^u\map P$ in $P$ are called equivalent if the bundles  $P^u$, $P'^u$ are 
isomorphic 
over $P$.
A $Spin^{U(2)}(n)$-structure in an oriented Riemannian 4-manifold $(X,g)$
 is a
$Spin^{U(2)}(n)$-structure in the bundle $P_g$ of oriented coframes. 
\\

We refer to [T1], [T2] for the following classification result:
\begin{pr} Let $P$ be a principal $SO$-bundle, $\bar P$ a $PU(2)$-bundle, and
$L$ a Hermitian line bundle  over $X$.\\
 i) $P$ admits a $Spin^{U(2)}$-structure $P^u
\rightarrow  P$  with 
$$P^u\times_{\bar \delta}PU(2)\simeq\bar P\ ,\ \ P^u\times_{\det}\C\simeq L$$
if and only if $w_2(P)=w_2(\bar P)+\overline c_1(L)$, where 
$\overline c_1(L)$ is
the mod 2 reduction of $c_1(L)$ .\\
 ii) If the base $X$ is a compact oriented 4-manifold, then the 
map 
$$P^u\longmapsto \left([P^u\times_{\bar\delta} PU(2)],
[P^u\times_{\det}\C]\right)$$
 defines a 1-1 correspondence between the set of  isomorphism classes of
 $Spin^{U(2)}$-struc\-tures in  $P$  and   the set of pairs of isomorphism
 classes
$([\bar P],[L])$, where
$\bar P$ is a $PU(2)$-bundle and $L$ an $S^1$-bundle with $w_2( 
P)=w_2(\bar P)+\overline c_1(L)$. The latter set can be identified with
$$\{(p,c)\in H^4(X,\Z)\times H^2(X,\Z) |\ p\equiv (w_2(P)+ 
\bar c)^2\ {\rm mod}\ 4\} 
$$
\end{pr}
\qed

The group $Spin^{U(2)}(4)$ can be written as 
$$Spin^{U(2)}(4)=\qmod{SU(2)_+\times
SU(2)_-\times U(2)}{\Z_2}\ ,$$
hence it comes with natural orthogonal representations
$$\ad_{\pm}: Spin^{U(2)}(4)\map so(su(2)), 
$$
defined by the adjoint representations of $SU(2)_{\pm}$, and with 
natural unitary
representations
$$\sigma_{\pm}:Spin^{U(2)}(4)\map
U(\H_{\pm}\otimes_\C \C^2)
$$
obtained by coupling  the canonical representations of $SU(2)_{\pm}$ 
 with the
canonical representation of $U(2)$.

We denote by $\ad_{\pm}(P^u)$, $\Sigma^{\pm}(P^u)$ the 
corresponding associated
vector bundles. The Hermitian 4-bundles $\Sigma^{\pm}(P^u)$ 
are called the spinor
 bundles of
$P^u$, and the sections in these bundles are called spinors.\\

We refer to [T2] for the following simple  result
\begin{pr} Let $P$ be an $SO(4)$-bundle whose second Stiefel-Whitney
 class admits
integral lifts. 

There is a 1-1 correspondence between isomorphism classes of
  $Spin^{U(2)}$-structures   in $P$ and   equivalence classes of pairs
consisting of a $Spin^c(4)$-structure $P^c\map  P$ in
$P$ and a  $U(2)$-bundle $E$. Two pairs are considered equivalent if, 
after tensoring
the first one with a line bundle, they become isomorphic over $P$. 
\end{pr}

Suppose that $P^u$ is associated with the pair
$(P^c,E)$, and let $\Sigma^{\pm}_c$ be the spinor
bundles corresponding to $P^c$. Then  the  associated bundles 
$P^u\times_\pi\R^4$, $\Sigma^{\pm}(P^u)$,
$\bar\delta(P^u)$,
$\det(P^u)$, ${\G(P^u)}$, $\G_0(P^u)$
  can be expressed in terms of the pair $(P^c,E)$ as follows:
$$P^u\times_\pi\R^4=\RSU(\Sigma^+_c,\Sigma^-_c)\ ,
$$
$$\Sigma^{\pm}(P^u)=[\Sigma^{\pm}_c]^{\vee}\otimes
E=\Sigma^{\pm}_c\otimes E^{\vee}\otimes[\det(P^u)]    ,\  
\bar\delta(P^u)\simeq
\qmod{P_E}{S^1},\ \ad_{\pm}(P^u)=su(\Sigma^{\pm}_c)$$   $$  \det(P^u)\simeq \det
(P^c)^{-1}\otimes (\det E),\
 \G_0(P^u)=SU(E),\ \gr_0(P^u)=su(E)\ .
 $$

Here we denoted by  $\RSU(\Sigma^+_c,\Sigma^-_c)$ the   bundle   of real
multiples of
$\C$-linear isometries  of determinant 1  from $\Sigma^+_c$ to $\Sigma^-_c$. 
The
euclidean structure and the orientation in this bundle are fibrewise defined by 
the Pauli
matrices associated with a pair  of frames $(e_1^{\pm},e_2^{\pm})$ in 
$\Sigma^{\pm}_c$,
satisfying $e_1^+\wedge e_2^+=e_1^-\wedge e_2^-$.

The data of a $Spin^{U(2)}(4)$-structure $P^u\map P$ in an $SO(4)$-bundle 
$P$ is
equivalent to the data of an orientation preserving linear isometry
$$\gamma:P\times_{SO(4)}{\R^4}\map
P^u\times_\pi\R^4=\RSU(\Sigma^+_c,\Sigma^-_c)\subset \Hom_{\G_0}
(\Sigma^+(P^u),
\Sigma^-(P^u))
$$
which will be called  the \underbar{Clifford} \underbar{map} of 
the structure. 

The Clifford map $\gamma$ induces  isomorphisms
$$\Gamma_{\pm}:\Lambda^2_{\pm}(P\times_{SO(4)}{\R^4})\map
su(\Sigma^{\pm}_c)=\ad_{\pm}(P^u)\ ,
$$
which multiply the norms by 2 ([DK] p. 77, [OT1]).

The following simple remark will play a fundamental role in this 
paper:\\

{\it Suppose that
$\Lambda$ is a real oriented 4-bundle, and $\gamma:\Lambda\map
P^u\times_\pi\R^4$ an orientation preserving linear isomorphism.
Then
$\gamma$ defines an Euclidean structure $g_\gamma$ on
$\Lambda$ such that
$\gamma$ becomes the Clifford map of a
$Spin^{U(2)}(4)$-structure in $(\Lambda,g_\gamma)$.}
\subsection{The $PU(2)$-monopole equations}

Let $\sigma:P^u\map P_g$ be a $Spin^{U(2)}(4)$-structure in the 
oriented compact
 Riemannian
4-manifold $(X,g)$. Fix a connection $a\in{\cal A}(\det(P^u))$. Using the 
third exact
 sequence in
(1), we see that the data of a connection $A\in {\cal A}(\bar\delta(P^u))$ is 
equivalent 
to the
data of a connection $B_{A,a}$ in $P^u$ which lifts the Levi-Civita connection
 in $P_g$ 
and the
fixed connection $a$ in $\det(P^u)$  (via the maps $P^u\map P_g$ and  
$P^u\map \det(P^u)$
respectively).  The
Dirac operator $\Dr_{A,a}$ associated with the pair $(A,a)$ is the first 
order elliptic 
operator
$$\Dr_{A,a}:A^0(\Sigma^{\pm}(P^u))\textmap{\nabla_{B_{A,a}}} 
A^1(\Sigma^{\pm}(P^u))\stackrel{\gamma}{\map}A^0(\Sigma^{\mp}(P^u))
$$
Regarded as operator
$\Sigma^+(P^u)\oplus\Sigma^-(P^u)\map\Sigma^+(P^u)\oplus\Sigma^-(P^u)$,
 the Dirac
operator
$\Dr_{A,a}$ is also selfadjoint.

We define the quadratic map $\Sigma^{\pm}(P^u)\map \ad_+(P^u)\otimes\gr_0$,
$\Psi\longmapsto (\Psi\bar\Psi)_0$ by
$$(\Psi\bar\Psi)_0:=pr_{\ad_+(P^u)\otimes\gr_0}(\Psi\otimes\bar\Psi) \ ,
$$
where $pr_{\ad_+(P^u)\otimes\gr_0}$ denotes the orthogonal  projection
$$\Herm(\Sigma^+(P^u))\map{\ad_+(P^u)\otimes\gr_0}\ .$$

We introduce now the $PU(2)$-Seiberg-Witten equations $SW^\sigma_a$ 
associated to the
pair $(\sigma,a)$, which are equations for a pair $(A,\Psi)$ formed by  a
$PU(2)$-connection $A\in{\cal A}(\bar\delta(P^u))$ and a positive spinor
 $\Psi\in
A^0(\Sigma^+(P^u))$:
$$\left\{\begin{array}{ccc}
\Dr_{A,a}\Psi&=&0\\
\Gamma(F_A^+)&=&(\Psi\bar\Psi)_0
\end{array}\right. \eqno{(SW^\sigma_a)}
$$

The natural symmetry group of the equations is the gauge group ${\cal
G}_0:=\Gamma(X,\G_0)$. We denote by ${\cal M}^\sigma_a$ the
 moduli space 
$${\cal M}^\sigma_a:=\qmod{\left[{\cal A}(\bar\delta(P^u))\times
A^0(\Sigma^+(P^u))\right]^{SW^\sigma_a}}{{\cal G}_0} \ ,
$$
where $\left[{\cal A}(\bar\delta(P^u))\times
A^0(\Sigma^+(P^u))\right]^{SW^\sigma_a}$ denote the space of solutions
 of the equations
$(SW^\sigma_a)$. Using the well-known Kuranishi method one can endow ${\cal
M}^\sigma_a$ with the structure of a ringed space, which has locally the form
$\qmod{Z(\theta)}{G}$, where $G$ is a closed subgroup of $SU(2)$ acting on finite
dimensional vector spaces $H^1$, $H^2$, and 
$Z(\theta)$ is the real analytic space cut-out by a $G$-equivariant real analytic map
$H^1\supset U \stackrel{\theta}{\map}  H^2\  $
(see [OT5],  [T1], [T2] for details).

\section{Smooth moduli spaces}

\subsection{The difficulty}

Equations for pairs $(A,\Psi)$, where $A$ is a unitary connection 
with fixed
determinant connection and
$\Psi$ a non-abelian  Dirac spinor   have been already considered 
  [PT1], [PT2]. For
instance, the definition of $Spin^c$-polynomial invariants  starts 
with the
construction of the moduli space of solutions of the   
$(ASD-Spin^c)$-equations
$$\left\{\begin{array}{lll}
\Dr_{A} \Psi&=&0 \ ,\ \ \ \Psi\ne 0\\
F_A^+&=&0\ .
\end{array}\right.$$

The proofs of the corresponding transversality results are incomplete. 
They are based 
on the following false statement ([PT2], [PT1]):\\ \\
{\bf (A)} Let $P^c\map P_g$ be a $Spin^c(4)$-structure with spinor
 bundles
$\Sigma^{\pm}(P^c)$ on a Riemannian 4-manifold $(X,g)$,   $E$ 
a Hermitian 2-bundle
on $X$, and  $A$   a unitary connection in  $E$. If the
$\Dr_{A}$-harmonic non-vanishing positive spinor
$\Psi\in A^0(\Sigma^+(P^c)\otimes E)$ is fibrewise degenerate 
  considered as morphism
$E^{\vee}\map \Sigma^+$, then
$A$ is reducible. 
 \\

In the proof of this assertion  ([PT1] p. 277) it was used   that, in the presence of  
a $Spin^c(4)$-structure, the  Clifford pairing $(\alpha,\sigma)\longmapsto
\gamma(\alpha)\sigma$  between 1-forms and positive spinors   has fibrewise no
  divisors of zero. This is  true for real 1-forms, but not for complex ones. 

Counterexamples are easy to find:

Every holomorphic section in a holomorphic Hermitian 2-vector bundle ${\cal E}$
on a K\"ahler surface can be regarded as  a degenerate  harmonic  positive spinor in
$\Sigma^+_{\rm can}\otimes {\cal E}$, where $\Sigma^+_{\rm
can}=\Lambda^{00}\oplus\Lambda^{02}$ is the positive spinor bundle of the
canonical $Spin^c(4)$-structure in $X$, if we endow ${\cal E}$ with the Chern
connection given by the holomorphic structure. Therefore any indecomposable
holomorphic 2-bundle ${\cal E}$ with $H^0({\cal E})\ne 0$ gives a counterexample to
the  assertion {\bf (A)}.\\

Note that these counterexamples occur precisely in the K\"ahler framework, where all  
explicit computations of moduli spaces and invariants were carried out.

\subsection{Partial transversality results}

Let $\sigma:P^u\map P_g$ be a $Spin^{U(2)}(4)$-structure  
 on $(X,g)$, denote by
$$\gamma:\Lambda^1\map\Hom(\Sigma^+,\Sigma^-)$$
 be the associated Clifford map, and 
let $C_0$ be a fixed  $SO(4)$-connection in 
$P^u\times_\pi SO(4)\simeq P_g$ (not
necessarily the Levi-Civita connection). We fix again a connection 
$a\in{\cal A}(\det(P^u))$. For
any connection $A\in {\cal A}(\bar\delta(P^u))$ we have an associated
 Dirac operator
$$\Dr_{a,A}^0=\gamma\cdot \nabla_{C_0,a,A}\ ,$$ 
where $\nabla_{C_0,a,A}:A^0(\Sigma^+)\map A^1(\Sigma^+)$ 
is the covariant 
derivative
associated with the connection in $P^u$ which lifts the triple $(C_0,a,A)$.

The role and the  properties of these slightly more general   Dirac 
operators will 
be cleared
up in the next section, where   $C_0$ will be a \ub{fixed} 
${\cal C}^{\infty}$-connection in  the
fixed bundle $P^u\times_\pi SO(4)$, but the metric $g$ and the
 Clifford map $\gamma$ will be
variable ${\cal C}^k$-parameters. 

Recall that one has  a canonical embedding
$P^u\times_\pi\C^4\subset\Hom(\Sigma^+,\Sigma^-)$, 
and that $\sigma$ defines an
isomorphism $\Lambda^1_\C\textmap{\simeq}P^u\times_\pi\C^4$.
 We consider the following 
equations
$$\left\{
\begin{array}{lcc}
 \Dr_{a,A}^0(\Psi)+\beta(\Psi)&=&0\\
\Gamma(F_A^+)&=&K(\Psi\bar\Psi)_0
\end{array}\right.\ ,
$$
which are  equations for a system 
$$(A,\Psi,\beta,K)\in {\cal A}(\bar\delta(P^u))\times
A^0(\Sigma^+)\times A^0(P^u\times_\pi\C^4)
\times \Gamma(X,GL(\ad_+)).$$
Complete the configuration space ${\cal A}:={\cal A}(\bar\delta(P^u))\times
A^0(\Sigma^+)$ with respect to a large Sobolev norm $L^2_l$, and the
 parameter space 
$${\cal
Q}:=A^0(P^u\times_\pi\C^4)\times \Gamma(X,GL(\ad_+))$$
 with respect to the Banach norm
${\cal C}^k$, $k\gg l$. 

The perturbations $(\beta,K)$ were also considered by Pidstrigach and Tyurin in
[PT2] in their attempt to get transversality for their version of non-abelian 
monopole equations.

  An $SU(2)\times SU(2)\times SU(2)$-reduction of $P^u$ over an open set 
$U\subset X$ induces
isomorphisms $\Sigma^{\pm}(P^u)|_U\simeq S^{\pm}\otimes E$ where 
$S^{\pm}$, $E$ are
 $SU(2)$-bundles. A spinor $\Psi\in \Sigma^+(P^u)$ will be called
\ub{degenerate} \ub{in} $x\in X$ if, with respect to an  $SU(2)\times 
SU(2)\times
SU(2)$-reduction around $x$, 
$\Psi_x\in S_x\otimes E_x=S_x^{\vee}\otimes E_x$ 
has
rank $\leq 1$. $\Psi$ will be called degenerate on $V\subset X$ if it
 is degenerate in
every point of $V$. 

A pair $(A,\Psi)\in {\cal A}(\bar\delta(P^u))\times A^0(\Sigma^+)$ will be 
called \ub{abelian} if
the connection $A$ is reducible, and the spinor $\Psi$ is contained in one of the
$A$-invariant summands of $\Sigma^+$.

If $(A,\Psi)\in {\cal A}(\bar\delta(P^u))\times A^0(\Sigma^+)$ is an 
abelian pair, then
$\Psi$ is clearly degenerate on $X$. However, the counterexamples in the
 previous section 
show that there exist  non-abelian pairs with non-trivial Dirac-harmonic 
spinor-component
which is  degenerate on $X$.

Let $sw=sw_{g,\sigma,C_0}:{\cal A}_l\times{\cal Q}^k\map
 A^0(\Sigma^-)_{l-1}\times
A^0(\ad_+\otimes\gr_0)_{l-1}$ be the map defined by the left hand 
side of the equations
above, and let
$${\cal N}^*:=\qmod{[{\cal A}_l^*\times{\cal Q}^k]\cap
 sw^{-1}(0)}{{\cal G}_{l+1}}$$
 be the moduli
space of solutions with non-trivial spinor-component. 
${\cal N}^*$ is the vanishing locus of the induced section $\bar {sw}$  in the Banach
bundle 
$$[{\cal A}_l^*\times{\cal Q}^k]\times_{{\cal G}_{l+1}}[A^0(\Sigma^-)_{l-1}\times
A^0(\ad_+\otimes\gr_0)_{l-1}]$$
over the Banach manifold ${\cal B}^*:=
\qmod{{\cal A}_l\times{\cal Q}^k}{{\cal G}_{l+1}}$
which is defined by $sw$.

The purpose of this section is to prove the following partial transversality theorem 
\begin{thry} If $sw$ is not a submersion in a solution $p=(A,\Psi,\beta,K)\in 
{\cal
A}_l^*\times{\cal Q}^k$, then $\Psi$ must be degenerate on $X$. In particular, 
${\cal N}^*$ is
smooth away from the closed subspace of solutions with globally degenerate spinor
 component.
\end{thry}
\pf Let   $(\Phi,S)\in
A^0(\Sigma^-)_{l-1}\times A^0(\ad_+\otimes\gr_0)_{l-1}$ a pair
 which is 
$L^2$-orthogonal to
$\im (d_p)$. Using the perturbation $\beta$ we get immediately that
$Re(\dot\beta,\Phi\otimes\bar\Psi)$ vanishes for every variation
 $\dot\beta$ of
 $\beta$. 
With respect to any  local $SU(2)\times SU(2)\times SU(2)$-reduction 
$(S^{\pm},E)$ of
$P^u|_U$ ($U$ an open set) the contraction of $\Phi\otimes\bar\Psi$ with 
the Hermitian
metric in $E$ must vanish, which shows that pointwise 
$\Psi(v^+)\bot \Phi(v^-)$
 for every 
$v^{\pm}\in S^{\pm}_u$,
 $u\in U$. If $\Psi$ has rank 2 in a point $x\in X$, then $\Phi$ must vanish
 identically on a
neighbourhood of $x$.

Also, if $\Psi$ has rank 2 in $x$, then $(\Psi\bar\Psi)_0$ has rank 3 in $x$ as map
$\ad_{+,x}^{\vee}\map \gr_{0,x}$, hence the same argument as above shows that $S$ 
must
vanish on a neighbourhood of $x$. Therefore the pair $(\Phi,S)$ must be zero on a
neighbourhood of $x$.

We can assume that $A$ is the Coulomb gauge with respect  to a smooth connection 
$A_0$. 
Therefore, by Agmon-Douglis-Nirenberg's non-linear-elliptic regularity theorems 
(see for
instance [B],  p. 467, Theorem 41), it follows that
$(A,\Psi)$ is a smooth pair (if the Clifford map 
$\Lambda^1\map P^u\times\R^4$ 
had only
class ${\cal C}^k$, we would have got a ${\cal C}^{k+1-\varepsilon}$-pair, which 
is enough to
complete the argument).  Using now variations
$(\dot A,\dot\Psi)$, we see that  $(\Phi,S)$ must satisfy an elliptic system of the 
form
$$\tilde D^1_{A,\Psi}[\tilde D^1_{A,\Psi}]^*(\Phi,S)=0 \ .
$$
Here $\tilde D^1_{A,\Psi}$ is the first derivative in $(A,\Psi)$ of the map  
$\tilde sw$
obtained by dividing by 2 the second component of $sw$ such that the symbol of $\tilde
D^1_{A,\Psi}[\tilde D^1_{A,\Psi}]^*$ becomes a scalar, and Aronszajin's
 theorem applies. It
follows that $(\Phi,S)=0$, because it vanishes on a non-empty open set.

\begin{re} The same result holds if ${\cal Q}^k$ is replaced by  any product 
${\cal Q}^k\times
{\cal R}$ of ${\cal Q}^k$ with a Banach manifold $R$,   and $sw$ by a smooth
 map $sw':{\cal
A}_l
\times {\cal Q}^k\times {\cal R}$ whose restriction  to any fibre ${\cal A}_l \times
{\cal Q}^k\times \{r\}$ has the form $sw_{g,\sigma,C_0}$ for a metric $g$, a
$Spin^{U(2)}$-structure $\sigma$ in  $(X,g)$, and an
$SO(4)$-connection $C_0$.
\end{re}

An easy way to parameterize the space of pairs consisting of a metric and a
$Spin^{U(2)}(4)$-structure will be given in the next section.

\subsection{$PU(2)$-monopoles with degenerate spinor component. Generic regularity}

Let   $P^u$ be a $Spin^{U(2)}$-bundle.   Suppose 
that  the spinor $\Psi\in A^0(\Sigma^+)$ is degenerate on a whole neighbourhood of 
 a point
$x\in X$ but $\Psi_x\ne 0$, and let $A\in\bar\delta(P^u))$ be
 a $PU(2)$-connection. 
The pair
$(A,\Psi)$ will be called \ub{non-abelian} \ub{in}
\ub{$x$} if the second fundamental form of the line subbundle $L\subset E$ generated
 by $\Psi$
  around $x$ is non-zero in $x$.

We recall that if $P^u$ is associated with a pair $(P^c,E)$, where $P^c$ is a 
$Spin^c(4)$ bundle
$P^c$ of spinor bundles $\Sigma^{\pm}_c$ and $E$ is a $U(2)$-bundle, then
$\Sigma^{\pm}=[\Sigma^{\pm}_c]^{\vee}\otimes E=\Sigma^{\pm}_c\otimes
E^{\vee}\otimes\det(P^u)$ and
$P^u\times_\pi\R^4=\RSU(\Sigma^+_c,\Sigma^-_c)\subset
\Hom(\Sigma^+_c,\Sigma^-_c)\subset\Hom(\Sigma^+,\Sigma^-)$. The 
euclidean structure and
the orientation in the real 4-bundle $\RSU(\Sigma^+_c,\Sigma^-_c)$ are fibrewise 
defined by
the Pauli matrices associated with frames $(e_1^{\pm},e_2^{\pm})$ 
of $\Sigma^{\pm}_x$
satisfying
$e_1^+\wedge e_2^+=e_1^-\wedge e_2^-$.
\begin{dt} Let $P^u$ be a $Spin^{U(2)}$-bundle with $P^u\times_\pi\R^4
\simeq \Lambda^1$. A
Clifford map is an orientation preserving linear isomorphism
$$\gamma:\Lambda^1\map
P^u\times_\pi\R^4=\RSU(\Sigma^+_c,\Sigma^-_c)\subset
\Hom(\Sigma^+,\Sigma^-)\ .$$ 
\end{dt}

 Every ${\cal C}^k$ Clifford map $\gamma:\Lambda^1 \map P^u\times_\pi \R^4$   
defines a
${\cal C}^k$ metric $g_\gamma$ on $X$ which makes $\gamma$ an 
isometry, so that
$\gamma:\Lambda^1 \map P^u\times_\pi
\R^4\subset
\Hom(\Sigma^+,\Sigma^-)$ becomes the Clifford map of a 
$Spin^{U(2)}$-structure
$\sigma_\gamma$ in $(X,g_\gamma)$.

This formalism will play an important role in this paper. The space 
$$Clif:=\Gamma(X,{\rm
Iso}_+(\Lambda^1,P^u\times_\pi\R^4))$$
of Clifford maps parameterizes the set of pairs consisting of a metric and a
$Spin^{U(2)}(4)$-structure for that metric. Note that the metric 
   determines
a $Spin^{U(2)}$-structure with a given bundle $P^u$  only up to an 
$SO(4)$-gauge
transformation of the cotangent bundle.

As in the previous section fix  a ${\cal C}^{\infty}$ $SO(4)$-connection 
$C_0$ in
$P^u\times_\pi SO(4)$. To any pair of connections $(a,A)\in 
{\cal A}(\det(P^u))\times{\cal
A}(\bar\delta(P^u))$ we associate    a Dirac operator
$\Dr_{\gamma,a,A}^0$   using the Clifford map $\gamma$ and the lift
$\nabla_{C_0,a,A}:A^0(\Sigma^+)\map A^1(\Sigma^+)$  of
$(C_0,a,A)$: 
$$\Dr_{\gamma,a,A}^0=\gamma\cdot \nabla_{C_0,a,A}\ . 
$$
This Dirac operator does not coincide with the   standard Dirac operator 
$\Dr_{\gamma,a,A}$
associated with $(A,a)$ and the $Spin^{U(2)}$-structure on $(X,g_\gamma)$
 defined by
$\gamma$, because
$\gamma^{-1}(C_0)$ may be different from the Levi-Civita connection in
$(\Lambda^1,g_\gamma)$; however, it has the same symbol as the standard one.  
 The
advantage of using these Dirac operators, is that they depend in a very simple way on
$\gamma$ and that they are operators with ${\cal C}^k$-coefficients if $\gamma$ is
 of class
${\cal C}^k$. The coefficients of the Levi-Civita connection in 
$(\Lambda^1,g_\gamma)$ are 
in general only of class ${\cal C}^{k-1}$, and the coefficients of the 
induced Levi-Civita
connection in
$P^u\times_\pi \R^4$ are also of class ${\cal C}^{k-1}$, so that the 
coefficients of the
standard Dirac operator $\Dr_{\gamma,a,A}$  have   a regularity-class  
  smaller by 1 than the regularity class of $\gamma$.

The use of these Dirac operators, whose coefficients do not contain 
 the derivatives
 of the
Clifford map, is essential in our proofs.
\begin{re} There exists a section $\beta=\beta(\gamma,C_0)\in {\cal
C}^{k-1}(P^u\times_\pi\C^4)$ such that
$\Dr_{\gamma,a,A}^0=\Dr_{\gamma,a,A}+ \beta $.
\end{re}

To see this, let $C_\gamma$ be the $SO(4)$-connection in 
$P^u\times_\pi\R^4$
 induced via
$\gamma$ by the Levi Civita connection in $(\Lambda^1,g_\gamma)$. 
 The difference
$\alpha:=\nabla_{C_\gamma,a,A}-\nabla_{C_0,a,A}$ is an 
$\ad_+$-valued 1-form of class
${\cal C}^{k-1}$, hence an element in  
$${\cal C}^{k-1}(\Lambda^1(\ad_+))=
{\cal C}^{k-1}(\Lambda^1(su(\Sigma^+_c)))\subset {\cal
C}^{k-1}(\Lambda^1(\End(\Sigma^+)))$$
 which  does not depend on $(A,a)$. In local
coordinates, $\alpha$ has the form
$\alpha=\sum u^i\otimes \alpha_i $, with local sections 
$\alpha^i$  in $su(\Sigma^+_c)$. Its
contraction with
$\gamma$ has locally the form
$\sum\gamma(u^i)\circ 
\alpha_i$, and defines a ${\cal C}^{k-1}$-section  $\beta$ in  
$\Hom(\Sigma^+_c,\Sigma^-_c)=P^u\times_\pi\C^4 $.
\qed

Consider the following $PU(2)$-monopole equations
$$
\left\{\begin{array}{lcl}
\Dr_{\gamma,a,A}^0\Psi&=&0\\
\Gamma_\gamma(F_A)&=&(\Psi\bar\Psi)_0 \ .
\
\end{array}\right.\eqno{(SW_a)}
$$
for a triple $(A,\Psi,\gamma)\in{\cal A}(\bar\delta(P^u))\times A^0(\Sigma^+)\times Clif$.
The map   
$$\Gamma_\gamma:\Lambda^2\map
\End(\Sigma^+_c)\subset\End(\Sigma^+)$$
 is determined by
$\gamma$ via the formula
$$\Gamma_\gamma(u\wedge
v)=\frac{1}{2}\left(-\gamma(u)^*\gamma(v)+\gamma(v)^*\gamma(u)\right)
$$
  and vanishes identically on $\Lambda^2_{-,g_\gamma}$, so that we 
could have  written 
 $F_A^{+_{g_{\gamma}}}$ instead of $F_A$ in the second equation. 
In the form above it will be
easier to compute the derivative  with respect to $\gamma$.

Complete the configuration space 
${\cal A}:={\cal A}(\bar\delta(P^u))\times A^0(\Sigma^+)$ 
with respect to a large Sobolev norm $L^2_l$ and the space of
 Clifford maps $Clif$ with
respect  to the Banach norm ${\cal C}^k$, $k\gg l$.

Before stating  the main result of this section, we begin with two simple remarks
\begin{re} Let $A$, $F$ be subspaces of a normed space $H$ with $F$ finite
dimensional.   Then
$$\overline{A+F}=\bar A +F\ .
$$
\end{re}
Indeed, $\overline{A+F}\supset 
 \bar A$, and $\overline{A+F}\supset  F$, hence
$\overline{A+F}\supset \bar A +F$. To prove the opposite inclusion, it is enough to
notice that  $\bar A+ F\supset A+F$ and to prove that $\bar A+ F$ is closed. Let
$q: H\map  \qmod{H}{\bar A}$ be the canonical projection. The right hand
space is also normed, hence $q(F)\subset
\qmod{H}{\bar A}$ is   closed (being finite dimensional), and therefore 
$q^{-1}(q(F))=\bar A+ F$
is closed in
$H$, since $q$ is continuous. This proves the remark
\qed

\begin{re} Let $f:H_1\map H_2$ be a continuous operator with closed
image and finite dimensional kernel between Banach spaces , and let
$A\subset H_1$ be a closed subspace. Then $f(A)$ is closed.
\end{re}
\pf $f$ factorizes as $H_1\stackrel{p}{\map} \qmod{H_1}{\ker f}
 \stackrel{\simeq}{\map}
f(H_1)\hookrightarrow H_2$, where the middle arrow is an isomorphism
  by the Banach
Theorem. Therefore it is enough  to show that  $p(A)$ is closed, or equivalently that
$p^{-1}(p(A))=A+\ker f$ is closed. But this follows by  the remark above.
\qed

Let $[{\cal A}_l\times Clif^k]^{SW}$  be the   space of solutions
$(A,\Psi,\gamma)$ of the equations above, and let  
$[{\cal A}_l\times Clif^k]^{SW}_U$ be the  
subspace of solutions whose spinor component is 
degenerate on the  open set  $U$.  

 The space $[{\cal A}_l\times Clif^k]^{SW}_U$ is a closed real
analytic subspace of the space $[{\cal A}_l\times Clif^k]^{SW}$, since it is
  the vanishing locus
of the  (real analytic)   map 
$${\cal A}_l\map A^0(\Sigma^+)_l\stackrel{\det}{\map}
A^0(\det(P^u))_l\stackrel{res_U}{\map}A^0(\det(P^u)|_U)_l \ .
$$ 
We can now state the main result of this section. 
\begin{thry}
  Let $\theta=(A,\Psi,\gamma)\in [{\cal A}_l\times
Clif^k]^{SW}_U$, and suppose that for a point $u\in U$,  one has $\Psi_u\ne 0$, 
and the pair
$(A,\Psi)$ is non-abelian in
$u$. Then the image of the Zariski tangent space $T_{\theta}[{\cal A}_l\times
Clif^k]^{SW}_U$ under the projection 
$$T_{\theta}[{\cal A}_l\times
Clif^k]^{SW} \map T_\gamma(Clif^k)={\cal C}^k(\Hom(\Lambda^1,
P^u\times_\pi\R^4))$$
   has infinite
codimension.
\end{thry}

For the proof of the theorem, we need some preparations:  

Note first (using [DK], p. 135) that we may assume that the
Sobolev connection $A$ is in Coulomb gauge with respect to a smooth connection 
$A_0$ and a
fixed  smooth metric
$g_0$, i.e.
$$d_{A_0}^{*_{g_0}}(A-A_0)=0 \ .
$$
Put $\alpha:=A-A_0$, hence $F_A=d_{A_0}\alpha+\alpha\wedge\alpha+F_{A_0}$. 
The 
differential operator $\Gamma_\gamma\circ d_{A_0}+d_{A_0}^{*_{g_0}}$ is 
{\it elliptic}
 although the metrics $g_0$ and $g_\gamma$ may be   different, and it has  
 coefficients of class ${\cal C}^k$.  Note also that $\Gamma_\gamma\circ
d_{A_0}+d_{A_0}^{*_{g_0}}$ is an operator between ${\cal C}^{\infty}$-bundles.

The Dirac operator $\Dr_{\gamma,a,A_0}^0=\Dr_{\gamma,a,A}^0-\gamma(\alpha)$ has  
coefficients of class ${\cal C}^k$. Therefore, the  pair $(\alpha,\Psi)$ is a solution of the
non-linear elliptic system
$$\left\{\begin{array}{lcl}
\Dr_{\gamma,a,A_0}^0\Psi+\gamma(\alpha)\Psi&=&0\\
\Gamma_\gamma(d_{A_0}\alpha+\alpha\wedge\alpha+F_{A_0})&=&(\Psi\bar\Psi)_0\\
d^{*_{g_0}}_{A_0}\alpha&=&0 \ \ \ \ .
\end{array}\right.
$$
Writing the left hand side as a function of $x^j$, $\alpha^k$, 
$\Psi^l$, $\partial_j \alpha^k$,
$\partial_j \Psi^l$ (with respect to a   smooth chart and  bundle trivializations), 
we see that
this function has class ${\cal C}^k$ in this system of variables  (in fact it is polynomial
of degree 2 in the last four group of variables).  It follows, by
Agmon-Douglis-Nirenberg's non-linear-elliptic regularity theorems ([B], p. 467, 
Theorem 41)
that
$\alpha$,
$\Psi$, hence also the pair $(A,\Psi)$,  have class
${\cal C}^{k+1-\varepsilon}$.  It would have class ${\cal C}^{k+1}$ if we had 
chosen a
non-integer index $k=[k]+\varepsilon$, i.e. if we had worked with the 
H\"older space ${\cal
C}^{[k],\varepsilon}$.\\

Let $sw:{\cal A}_l\times Clif^k\map 
A^0(\Sigma^-)_{l-1}\times A^0(\ad_+\otimes\gr_0)_{l-1}$ be the 
map given by the left hand
side  of the equations $(SW_a)$, and  put $\det_U:=res_U\circ \det$. 

The tangent space
$T_\theta[{\cal A}_l\times Clif^k]^{SW}_U$ is the space of solutions $(\dot
A,\dot\Psi,\dot\gamma)$ of the linear system
$$\left\{
\begin{array}{lcl}
\frac{\partial sw}{\partial(A,\Psi)}|_\theta(\dot A,\dot\Psi)+\frac{\partial
sw}{\partial\gamma}|_\theta(\dot\gamma)&=&0\\
d_\Psi(\det_U)(\dot \Psi)&=&0 \ .
\end{array}
\right.
$$

 Denote  by 
$$D_l^U=:\ker[d_\Psi({\det}_U)]\subset
A^0(\Sigma^+)_l$$
  the Zariski tangent space at
$\Psi$ to the space
${\cal D}_l^U:=\det_U^{-1}(0)$ of $L^2_l$ positive spinors which are 
degenerate  on $U$.\\

Theorem 3.7 can now be reformulated   as follows
\begin{pr}  The subspace
$$\left(\frac{\partial sw}{\partial\gamma}|_\theta\right)^{-1} 
\left(\frac{\partial sw}{\partial(A,\Psi)}|_\theta \ (A^1(\gr_0)_l\times
D_l^U)\right)
$$
has infinite codimension in ${\cal C}^k(\Hom(\Lambda^1,P^u\times_\pi \R^4))$.
\end{pr}

In order to prove Proposition 3.8 we start by giving explicit formulas for the partial
derivatives above.  

The derivative  with respect to $\gamma$,
$$\left(\frac{\partial sw}{\partial\gamma}|_\theta\right):{\cal
C}^k(\Hom(\Lambda^1,P^u\times_\pi \R^4))\map A^0(\Sigma^-)_{l-1}\times
A^0(\ad_+\otimes\gr_0)_{l-1}\ ,
$$
is given by 
$$\left(\frac{\partial
sw}{\partial\gamma}|_\theta\right)(\dot\gamma)=
\left(\matrix{\dot\gamma(\nabla_{C_0,a,A}\Psi)\cr \cr
\frac{d}{d\gamma}(\Gamma_\gamma(F_A))(\dot\gamma)}\right)\ .
\eqno{(1)}$$
The derivative with respect to the pair  $(A,\Psi)$, 
$$\left(\frac{\partial sw}{\partial(A,\Psi)}|_\theta\right):  A^1(\gr_0)_l\times
A^0(\Sigma^+)_l
\map A^0(\Sigma^-)_{l-1}\times
A^0(\ad_+\otimes\gr_0)_{l-1}\ ,
$$
 is  
$$\left(\frac{\partial sw}{\partial(A,\Psi)}|_\theta\right)(\dot A,\dot \Psi)=
\left(\matrix{\Dr_{\gamma,a,A}^0\dot \Psi+\gamma(\dot
A)\Psi\cr \cr
 \Gamma_\gamma(d_A\dot A)-[(\dot\Psi\bar
\Psi)_0+(\Psi\bar{\dot\Psi})_0]}\right) \ .
\eqno{(2)}$$

The next two lemmata will translate the problem into a similar one  which 
involves only
Sobolev completions.

Let $j^k_{l-1}$ be the compact embedding 
$$j^k_{l-1}:{\cal
C}^k(\Hom(\Lambda^1,P^u\times_\pi \R^4)) \map
A^0(\Hom(\Lambda^1,P^u\times_\pi \R^4))_{l-1}\ .
$$
\begin{lm}\hfill{\break}
1.  The linear operator $\left(\frac{\partial
sw}{\partial\gamma}|_\theta\right)$ has a continuous extension to the Sobolev completion
$A^0(\Hom(\Lambda^1,P^u\times_\pi \R^4))_{l-1}$. More
precisely, formula (1) defines a linear continuous map 
$$a_{l-1}: A^0(\Hom(\Lambda^1,P^u\times_\pi \R^4))_{l-1}\map
A^0(\Sigma^-)_{l-1}\times A^0(\ad_+\otimes\gr_0)_{l-1}
$$
such that 
$$\left(\frac{\partial
sw}{\partial\gamma}|_\theta\right)=a_{l-1}\circ j^k_{l-1}\ .
$$

2. The space $ \frac{\partial
sw}{\partial(A,\Psi)}|_\theta \ (A^1(\gr_0)_l\times D_l^U) $ is closed in
$$A^0(\Sigma^-)_{l-1}\times A^0(\ad_+\otimes\gr_0)_{l-1} .$$
\end{lm}
\pf 
1. The first assertion follows easily, since $\nabla_{C_0,a,A}\Psi$ and $F_A$ have
regularity class  ${\cal C}^{k-\varepsilon}$, and $\gamma$ has regularity
 class ${\cal C}^k$.
Therefore, working  in local
${\cal C}^\infty$-coordinates, the expression
$$(\frac{d}{d\gamma}\Gamma) (\dot\gamma)(F_A)= 
\frac{d}{d\gamma}
\left( \frac{1}{2}[-{\gamma(u^i)}^*\gamma(u^j)+
(\gamma(u^j)^*\gamma(u^i)]\otimes
F_{A,ij}\right)(\dot \gamma)
$$
is a linear operator  of order 0 with ${\cal C}^{k-\varepsilon}$
 coefficients in the variable
$\dot\gamma$. 
\\

2.    Decompose $A^1(\gr_0)_l\times A^0(\Sigma^+)_{l}$ as
$$A^1(\gr_0)_l\times A^0(\Sigma^+)_{l}=
D^0_{(A,\Psi)} [A^0(\gr_0)_{l+1}]\oplus\ker [D^0_{(A,\Psi)}]^*=
\im D^0_{(A,\Psi)}\oplus\ker [D^0_{(A,\Psi)}]^*  
$$
where $D^i_{(A,\Psi)}$ are the differential operators in the  fundamental elliptic
complex associated with  the pair $(A,\Psi)$ and the metric $g_\gamma$.  
The decomposition 
is
$L^2_{g_\gamma}$-orthogonal.

The subspace $A^1(\gr_0)_l\times D_l^U\subset A^1(\gr_0)_l\times
A^0(\Sigma^+)_{l} $ is closed, and contains the first summand   $\im D^0_{(A,\Psi)}$  by  the
gauge-invariance property of the degeneracy-condition. Using the fact that
$D^1_{(A,\Psi)}\circ D^0_{(A,\Psi)}=0$, we get
$$ \frac{\partial
sw}{\partial(A,\Psi)}|_\theta \ (A^1(\gr_0)_l\times
D_l^U)= D^1_{(A,\Psi)}
\left[(A^1(\gr_0)_l\times D_l^U)\cap\ker (D^0_{(A,\Psi)})^*\right]= $$
$$=D^1_{(A,\Psi)}|_{_{\ker (D^0_{(A,\Psi)})^*}} 
\left[(A^1(\gr_0)_l\times
D_l^U)\cap\ker (D^0_{(A,\Psi)})^*\right]\ .
$$

But 
$D^1_{(A,\Psi)}|_{\ker (D^0_{(A,\Psi)})^*}:
{\ker (D^0_{(A,\Psi)})^*}\map
A^0(\Sigma^-)_{l-1}\times A^0(\ad_+\otimes\gr_0)_{l-1}$ 
is Fredholm and the
subspace $\left[(A^1(\gr_0)_l\times
D_l^U)\cap\ker (D^0_{(A,\Psi)})^*\right]$ of  
$\ker (D^0_{(A,\Psi)})^*$ is closed, so that the
assertion follows from Remark 3.6.
\begin{lm} If 
$$V:=\left(\frac{\partial sw}{\partial\gamma}|_\theta\right)^{-1} 
\left(\frac{\partial sw}{\partial(A,\Psi)}|_\theta \ (A^1(\gr_0)_l\times
D_l^U)\right)
$$
had finite codimension in 
${\cal C}^k(\Hom(\Lambda^1,P^u\times_\pi \R^4))$, then
$$V_{l-1}:=a_{l-1}^{-1} 
\left(\frac{\partial sw}{\partial(A,\Psi)}|_\theta \ (A^1(\gr_0)_l\times
D_l^U)\right)
$$
would have finite codimension in $A^0(\Hom(\Lambda^1,P^u\times_\pi
\R^4))_{l-1}$.
\end{lm}
\pf   
Suppose   there exists a finite dimensional subspace $F$ of   the space  
${\cal
C}^k(\Hom(\Lambda^1,P^u\times_\pi \R^4))$, such that
$$V+F={\cal
C}^k(\Hom(\Lambda^1,P^u\times_\pi \R^4))\ .
$$

Then  we have 
$$j^k_{l-1}(V)+j^k_{l-1}(F)=j^k_{l-1}({\cal
C}^k(\Hom(\Lambda^1,P^u\times_\pi \R^4))\subset 
A^0(\Hom(\Lambda^1,P^u\times_\pi
\R^4))_{l-1}\ ,
$$
hence
$$\overline{j^k_{l-1}(V)+j^k_{l-1}(F)}=A^0(\Hom(\Lambda^1,P^u\times_\pi
\R^4))_{l-1}\eqno{(4)}
$$
by the density property of   smooth sections in any Sobolev completion.

  Therefore, under the hypothesis of the
lemma,  and using (4) and Remark 3.5, one gets 
$$\overline{j^k_{l-1}(V)}+j^k_{l-1}(F)=
A^0(\Hom(\Lambda^1,P^u\times_\pi
\R^4))_{l-1}\ .\eqno{(5)}
$$

On the other hand, we know that  $\frac{\partial
sw}{\partial\gamma}|_\theta=a_{l-1}\circ j^k_{l-1}$. Therefore 
$$V=[j^k_{l-1}]^{-1}(V_{l-1})\ , 
$$
which shows that $j^k_{l-1}(V)\subset V_{l-1}$. But $V_{l-1}$ is closed  by
Lemma 3.9., hence $\overline{j^k_{l-1}(V)}\subset V_{l-1}$. From (5) 
it follows  that 
$$V_{l-1}+j^k_{l-1}(F)=A^0(\Hom(\Lambda^1,P^u\times_\pi
\R^4))_{l-1}
$$
which proves   Lemma 3.10.
\qed

The proof of Proposition 3.8 is now reduced to showing that $V_{l-1}$
 cannot have finite
codimension in $A^0(\Hom(\Lambda^1,P^u\times_\pi
\R^4))_{l-1}$. To prove this, we show that the sections in $V_{l-1}$ 
must fulfill a
very restrictive condition, which is not of finite codimension.

Let $v\in V_{l-1}$. Then, by definition
$$a_{l-1}(v)\in \frac{\partial sw}{\partial
(A,\Psi)}|_{\theta}(A^1(\gr_0)_{l}\times D_l^U)\ ,
$$
hence there exists a pair $(\dot A,\dot \Psi)\in A^1(\gr_0)_{l}\times D_l^U$  such
that

$$
\left\{\begin{array}{lcl}
\Dr_{\gamma,a,A}^0\dot \Psi+\gamma(\dot
A)\Psi&=& v(\nabla_{C_0,a,A}\Psi)\\ \\
\Gamma_\gamma(d_A\dot A)-[(\dot\Psi\bar
\Psi)_0+(\Psi\bar{\dot\Psi})_0]&=& 
\frac{d}{d\gamma}(\Gamma_\gamma(F_A))(v)
\ .
\
\end{array}\right.
$$
\\

Consider now small balls $U_1$, $U_2$ centered in $u$ such that 
$U_1\subsetint U_2\subset U$,
and such that the following two conditions hold:\\ 
 
1.  $\Psi$ is nowhere vanishing on $U_2$.\\
 
Let $S^{\pm}$, $E$ be the trivial $SU(2)$-bundles
associated with a $SU(2)\times SU(2)\times SU(2)$-reduction 
of $P^u|_{U_2}$. The
connection $C_0$ induces
${\cal C}^{\infty}$-connections in $S^{\pm}$, and  the pair $(A,a)$
 induces a
connection $B_A$ (with ${\cal C}^{k+1-\varepsilon}$-coefficients) in
 $E$ which
lifts the connection $A|_{U_2}$ in $ \bar\delta(P^u)|_{U_2}=
\qmod{P_E}{S^1}$ and the
connection $a|_{U_2}$ in $\det(P^u)|_{U_2}=\det(E)$.   Since $\Psi$ has
rank 1 in every point of
$U_2\subset U$,  it defines   a 
${\cal C}^{k+1-\varepsilon}$-splitting $E=L\oplus M$ with
$\Psi|_{U_2}\in A^0(S^+\otimes L)$. \\

2. The second fundamental form 
$b\in {\cal C}^{k+1-\varepsilon}(\Lambda^1_\C)$  of $L$ with
respect to the unitary connection $B_A$  (or, equivalently, 
with respect to $A$) is nowhere
vanishing on $U_2$.\\

Let $l,\ m$ be ${\cal C}^{k+1-\varepsilon}$ sections of
$E$ giving unitary frames in $L$ and $M$.  Then we can write   
$\Psi|_{U_2}=s_0^+\otimes l$, where
$s_0^+$ is a nowhere vanishing 
${\cal C}^{k+1-\varepsilon}$-section of $S^+$.  
Once we have
fixed this trivialization of $E$, we can identify the
 connections with the associated
connection matrices, and   write $B_A=A+\frac{1}{2}a\id$ 
\\

Recall that   $b$ is defined by   $b:=(\nabla_{B_A}l,m)$, and  
for any section $\varphi l$ of $L$
one has $\nabla_{B_A}(\varphi l)=
\nabla_{B_L}(\varphi l)+ \varphi b\otimes m$, where   $B_L$
(resp. $B_M$) are  the connections induced by
 $B_A$ in $L$ (resp. $M$). 

By the Dirac harmonicity condition, one has, 
taking the component of
$\Dr_{\gamma,a,A}^0\Psi$ in $S^-\otimes M$,
$$\gamma(b)(s_0^+)=0\ .
$$

Denote by $S_0$ the rank 1 subbundle of $S^+$ generated 
by the section $s_0^+$, and by
$S^{\bot}_0$ its orthogonal complement. Let $\Psi_t$ be a
 path of spinors with $\Psi_0=\Psi$
and
$\det(\Psi_t)=0$. Derivating it in 0, we get that the component 
of $\dot\Psi_0$ in
$S^\bot_0\otimes M$ must vanish. Therefore, the restriction
$\dot\Psi|_{U_1}$ of an element
$\dot\Psi\in D_l^U=T_\Psi({\cal D}^U_l)$ must have the form 
$$  \dot\Psi|_{U_1}=\dot
\sigma^+\otimes l+
\dot\zeta s_0^+\otimes m\ ,\ \  \dot \sigma^+\in
L^2_l(S^+|_{U_1}) \ ,\ \ \dot\zeta\in L^2_l(U_1,\C) \ .$$

Take now the component in $(S^-\otimes M)|_{U_1}$ of the
 restriction of the first
equation  to $U_1$.  Put $\nabla_{B_M}(m)=\lambda\otimes m$, 
where $\lambda$ is a ${\cal
C}^{k-\varepsilon}$ pure imaginary 1-form. 

One gets   the following equation on   $U_1$:
$$\Dr_{\gamma}^0(\dot\zeta s_0^+) +\dot\zeta\gamma(\lambda)(s_0^+) 
+\gamma(b)(\dot
\sigma^+) +\gamma(\dot A^2_1)(s_0^+)=v(b)(s_0^+) \ .
\eqno{(6)}$$

Here $\Dr_{\gamma}^0:A^0(S^+)_{s.}\map A^0(S^-)_{s-1}$ , $s\leq k$,
 stands for the   Dirac
operator associated with the Spin(4) structure on $(U_2, g_{\gamma})$ 
defined by $\gamma$
and  the $SO(4)$-connection $C_0|_{U_2}$ in $\RSU(S^+,S^-)$. 
$\Dr_{\gamma}^0$ is a  first
order elliptic operator with ${\cal C}^k$-coefficients. 
The complex 1-form
$\dot A^2_1$ is the component of
$\dot A$ written in the matricial form with respect to the 
decomposition $E=L\oplus M$.

The idea  to prove Proposition 3.8 is the following:

By the  properties 1., 2. above it follows that, varying $v$ in 
the equation (6), one can
get \ub{all} the
$L^2_{l-1}$-sections of the \ub{rank-2} bundle
$(S^-\otimes M)|_{U_1}$. But on the left of the same
 equation one has a differential operator
of order 1 with ${\cal C}^{k-\varepsilon}$ coefficients in $(\dot
\zeta,
\dot\sigma^+,\dot A^2_1)$ which has a  
\ub{non}-\ub{surjective} symbol: only the
complex valued function $\dot \zeta$, which
is a section in a \ub{rank-1} bundle on
$U_1$, is derivated on the left.  

The problem comes down to showing that the map 
$L^2_l\map L^2_{l-1}$ associated with such
an operator,
 cannot have a range of finite codimension.\\

We define the following operators:

$$res_{U_1}:A^0(\Hom(\Lambda^1,P^u\times_\pi \R^4))_{l-1}\map
A^0(\Hom(\Lambda^1,P^u\times_\pi \R^4)_{U_1})_{l-1}\ , $$
$$ev_{b,s_0^+}:A^0(\Hom(\Lambda^1,
P^u\times_\pi \R^4)_{U_1})_{l-1} \map
A^0(S^-|_{U_1})_{l-1}\ ,\ v'\longmapsto v'(b)(s_0^+)\ ,$$
$$[\Dr_{\gamma}^0]^-:A^0(S^-|_{U_1})_{l-1}
\map  A^0(S^+|_{U_1})_{l-2} \ ,
$$
$$pr^\bot:A^0(S^+|_{U_1})_{l-2}\map
A^0(S_0^\bot|_{U_1})_{l-2}\ .
$$
Here $[\Dr_{\gamma}^0]^-$ is the Dirac operator
 associated with   the connection $C_0$ and
the Clifford map
$\gamma^-:\Lambda^1\map \RSU(S^-,S^+)$ given by
$$\gamma^-(u)=-\gamma(u)^*\ .
$$

In general, the operator $[\Dr_{\gamma}^0]^-$ is not the
 formal adjoint of $\Dr_\gamma^0$,
because $\gamma^{-1}(C_0)$ can have non-vanishing torsion, 
but it has the same symbol as
$[\Dr_\gamma^0]^*$ and it is an operator with ${\cal C}^k$-coefficients. 
The associated
Laplacian 
$[\Dr_{\gamma}^0]^-\circ \Dr_\gamma^0$ has scalar symbol 
 given by $\xi\mapsto
-g_\gamma(\xi,\xi)\id_{S^+}$.
\begin{lm} \hfill{\break}
1. The operators $res_{U_1}$,  $pr^\bot$ are surjective.\\
2. The image of the operator $[\Dr_{\gamma}^0]^-:A^0(S^-|_{U_1})_{l-1}\map 
A^0(S^+|_{U_1})_{l-2}$ has finite codimension.\\
3. The operator $ev_{b,s_0^+}$ is  surjective. \\
\end{lm}

\pf 1. The surjectivity of $res_{U_1}$ follows from the 
extension theorems for Sobolev
spaces ([Ad], p. 83); the surjectivity of $pr^{\bot}$ is 
obvious. 

2. The fact that the image  of
$[\Dr_{\gamma}^0]^-:A^0(S^-|_{U_1})_{l-1}\map  A^0(S^+|_{U_1})_{l-2}$ has finite
codimension follows from the general theory of  elliptic
 operators (see for instance [BB]); 
It can also be directly verified as follows: We may suppose that
$X$ is the 4-sphere
$S^4$ and that
$S^\pm|_{U_1}$ are the restrictions to $U_1$ of the spinor
 bundles $S'^{\pm}$ associated
with a
$Spin(4)$-structure on $S^4$ whose Clifford map
$\gamma'$ extends $\gamma|_{U_1}$. We can also find a
 connection $C_0'$ in the
associated $SO(4)$-bundle extending
$C_0|_{\bar U_1}$.  

The image of $[\Dr_{\gamma}^0]^-$ contains the image of
 the composition
$res_{U_1}\circ[\Dr_{\gamma'}^0]^-$, where
$[\Dr_{\gamma'}^0]^-:A^0(S'^-)_{l-1}\map  A^0(S'^+)_{l-2}$ 
is the Dirac operator on the sphere
associated with $(\gamma')^-$ and $C_0'$. But $res_{U_1}$ is 
surjective and
$[\Dr_{\gamma'}^0]^-$ is Fredholm.

Note that  $[\Dr_{\gamma}^0]^-$   is in fact
\ub{surjective}, if $U_1$ is sufficiently small. 

3. The surjectivity of $ev_{b,s_0^+}$ is the crucial point  in which 
the fact that $s^+_0$ and
$b$ are nowhere vanishing on $U_2$ is used in an essential way.

We begin by  choosing a \ub{smooth} Clifford map
$$\gamma_0:\Lambda^1_{U_2}
\map  P^u|_{U_2}\times_\pi\R^4$$
such that $\gamma_0(b): S^+\map  S^-$ is an isomorphism in every
 point $u\in U_2$. 

This can be  achieved as follows: We know that $\gamma(b)(s_0^+)=0$,
 so the determinant
$\det(\gamma(b))$ of the induced morphism $\gamma(b):S^+\map  S^-$ 
must vanish. Therefore 
$g_\gamma^\C(b)=\det(\gamma(b))=0$, hence the real 
forms ${\rm Re}(b)$, ${\rm Im}(b)$ have
pointwise in
$U_2$ the same (non-zero !) $g_\gamma$-norm and are
 pointwise $g_\gamma$-orthogonal. It
suffices to choose $\gamma_0$ such that ${\rm Re}(b)$, ${\rm Im}(b)$
 are nowhere
$g_{\gamma_0}$-orthogonal on $U_2$.
With this choice $\gamma_0(b)(s_0^+)$ will be a nowhere
 vanishing section of $S^-$ on
$U_2$. 

Let now $s'\in A^0(S^-|_{U_1})_{l-1}$ be an arbitrary 
$L^2_{l-1}$-negative spinor. 

One can find a
unique $L^2_{l-1}$ section 
$\delta\in A^0(\RSU(S^-,S^-)|_{U_1}))_{l-1}$, such that
$\delta(\gamma_0(b)(s_0^+))=s'$: To see this, one 
uses the bilinear bundle map
$$ \RSU(S^-,S^-)\times S^-\map S^-\ .$$
The section $\delta$ is obtained by  fibrewise dividing 
(in the quaternionic sense) $s'$ by  the
smooth nowhere vanishing spinor $\gamma_0(b)(s_0^+)$ which is a ${\cal
C}^{k-\varepsilon}$-section on $U_2\supset \bar U_1$.

One also has a bilinear bundle map 
$$\RSU(S^+,S^-)\times \RSU(S^-,S^-)\map \RSU(S^+,S^-)
$$
which in local coordinates looks like quaternionic multiplication.

Now define  the $L^2_{l-1}$-morphism $v':\Lambda^1_{U_1}\map
\RSU(S^+|_{U_1},S^-|_{U_1})$   by
$$v'(\alpha):=
\delta\cdot [\gamma_0(\alpha)]\ ,\ \forall\ \alpha\in \Lambda^1_{U_1}\ .$$
This morphism defines a section in
$$A^0(\Hom(\Lambda^1_{U_1},P^u|_{U_1}\times_\pi\R^4)_{l-1}
=A^0(\Hom(\Lambda^1_{U_1},\RSU(S^+,S^-)|_{U_1})_{l-1}$$
 which acts   on complex 1-forms $\alpha$  by
$$v'(\alpha)(\cdot)=\delta[\gamma_0(\alpha)(\cdot)] \ .
$$
In particular, $v'(b)(s_0^+)=\delta[\gamma_0(b)(s_0^+)]=s'$.  
\qed

 After these preparations we can finally prove Proposition 3.8.

\pf  
We have to show that $V_{l-1}$ has infinite codimension in
$$A^0(\Hom(\Lambda^1,P^u\times_\pi \R^4))_{l-1}\ .$$

Take $v\in V_{l-1}$ and apply  
$[pr^\bot\circ [\Dr_{\gamma}^0]^-]$  to both sides
of (6). 

On the left, the only term containing second order derivatives of the sections 
$(\dot\zeta,\dot \sigma^+,\dot A^2_1)$ is 
$$[pr^\bot\circ
[\Dr_{\gamma}^0]^-](\Dr_{\gamma}^0(\dot\zeta s_0^+))\ .$$
 But, denoting   by $i_0$ the
bundle inclusion     $U_1\times\C \map S^+|_{U_1}$, 
$z\longmapsto z s_0^+$, one sees
that the 2-symbol of the composition
$$pr^\bot\circ \left[[\Dr_{\gamma}^0]^-
\circ \Dr_{\gamma}^0\right]\circ i_{0}$$
vanishes, since the symbol of the  Laplacian
$ [\Dr_{\gamma}^0]^-\circ \Dr_{\gamma}^0$ is scalar.

Therefore, applying $[pr^\bot\circ [\Dr_{\gamma}^0]^-]$
 on the left, 
one gets an expression
containing only first order derivatives of the Sobolev $L^2_l$ 
sections  $(\dot \zeta,
\dot\sigma^+,\dot A^2_1)$, hence an $L^2_{l-1}$-section of $S_0^{\bot}$.

 On the other hand applying $[pr^\bot\circ [\Dr_{\gamma}^0]^-]$
 on the right of (6), one
gets precisely 
$$\left[pr^\bot\circ [\Dr_{\gamma}^0]^-\circ ev_{b,s_0^+}\circ
res_{U_1}\right](v)\ .$$

Now consider the operator  
$$P:=\left[pr^\bot\circ
 [\Dr_{\gamma}^0]^-\circ ev_{b,s_0^+}\circ
res_{U_1}\right]:A^0(\Hom(\Lambda^1,
P^u\times_\pi \R^4))_{l-1}\map
A^0(S^\bot_0 |_{U_1})_{l-2}$$
and the following exact sequence
$$0\rightarrow\qmod{\im(P)}{P(V_{l-1})}\rightarrow
\qmod{A^0(S^\bot_0)_{l-2}}{P(V_{l-1})}\rightarrow   
\coker(P)\rightarrow 0\ .
$$

We have seen that $P(V_{l-1})$ is contained in 
 $A^0(S^\bot_0)_{l-1}$, which has
infinite codimension in $A^0(S^\bot_0)_{l-2}$.
\footnote{We used here
the following simple remark: The space of  $L^2_{l-1}$-sections
in the space of  $L^2_{l-2}$ sections in a bundle has infinite
codimension. Note that  $L^2_{l-1}$ is 
nonetheless dense in
$L^2_{l-2}$.} 

Therefore
$\qmod{A^0(S^\bot_0)_{l-2}}{P(V_{l-1})}$ has 
infinite dimension.  By Lemma 3.11
$\coker(P)$ has finite dimension, so that 
 $\qmod{\im(P)}{P(V_{l-1})}$ must have
infinite dimension.   But  $\qmod{\im(P)}{P(V_{l-1})}$ is 
a quotient of 
$$\qmod{A^0(\Hom(\Lambda^1,P^u\times_\pi 
\R^4))_{l-1}}{V_{l-1}}\ ,
$$
so that the latter must  also have infinite  dimension.
\qed
\\

 Let ${\cal M}^*$, ${\cal D}{\cal M}^*_U$ be the moduli spaces 
$${\cal M}^*:=\qmod{[{\cal A}_l^*\times 
Clif^k]^{SW_a}}{{\cal G}_{l+1}}\ ,\ 
{\cal D}{\cal M}^*_U:=\qmod{[{\cal A}_l^*
\times Clif^k]^{SW_a}_U}{{\cal G}_{l+1}}\ ,\ $$
where the upper script $(\ )^*$ denotes the subspace
 with non-zero spinor component.

\begin{co} Let 
$p=(A,\Psi,\gamma)\in  [{\cal A}_l^* \times Clif^k]^{SW_a}_U$ 
such that
for some $u\in U$, $\Psi_u\ne 0$ and $(A,\Psi)$ is 
non-abelian in $u$. Then the
Zariski tangent space 
$T_{[p]}{\cal D}{\cal M}^*_U$ has infinite codimension in 
$T_{[p]}{\cal M}^*$.  In
particular, $T_{[p]}{\cal D}{\cal M}^*_X$ has infinite 
codimension in $T_{[p]}{\cal M}^*$
for  every solution $p$ with non-abelian $(A,\Psi)$-component. 
\end{co}

\pf  We have 
$$pr_{T_\gamma(Clif^k)}(T_{[p]}{\cal M}^*)=
\frac{\partial sw}{\partial\gamma}^{-1}
\left[D^1_{(A,\Psi)}(A^1(\gr_0)_l\times
 A^0(\Sigma^+)_l)\right]\ ,
$$
and the vector space \  $D^1_{(A,\Psi)}(A^1(\gr_0)_l
\times A^0(\Sigma^+)_l)$ has   finite  
codimension in
$A^0(\Sigma^-)_{l-1}
\times A^0(\ad_+\otimes \gr_0)_{l-1}$.

Therefore, also the image of $T_{[p]}{\cal M}^*$ under the projection to
$T_\gamma(Clif^k)$ has finite codimension.

But, by Theorem 3.7, the image of $T_{[p]}{\cal D}{\cal M}^*_X$ 
under the same
projection has infinite codimension. This proves the
 first assertion.

The second assertion follows from Aronszajin's unique continuation
 theorem and the
fact that the vanishing locus of an harmonic spinor cannot separate domains
[FU]. Alternatively, one can use the Unique Continuation Theorem for monopoles
[FL] to see that a mnopole with non-vanishing spinor component, and
which is abelian on a non-empty open set, must be globally abelian. 

Therefore in the condition of the proposition we can find a point $x\in X$ with
  $\Psi_x\ne 0$ such that $(A,\Psi)$ is non-abelian in $x$.
\qed 

Using this result we can prove that for a generic Clifford map $\gamma$, 
the only degenerate
solutions in the moduli space ${\cal M}^*\cap p_{Clif^k}^{-1}(\gamma)$ 
are the abelian ones.
The idea is the following: 

Let  ${\cal D}{\cal M}^\circ_X\subset {\cal D}{\cal M}^*_X$ 
be the subspace of ${\cal D}{\cal
M}^*_X$ consisting of solutions with non-abelian $(A,\Psi)$-component. We have proven
that
${\cal D}{\cal M}^\circ_X$ has infinite codimension in ${\cal M}^*$. 
  Since the projection
${\cal D}{\cal M}^\circ_X\map Clif^k$ has "index  $-\infty$",  the
  generic fibre should be
empty.  There are of course two serious problems with this argument:\\ 
1. ${\cal D}{\cal M}^\circ_X$ is not smooth.\\
2. The restriction of the projection ${\cal D}{\cal M}^\circ_X\map Clif^k$
 to the smooth 
part is not Fredholm.

The idea to proceed is to   weaken locally the equation defining 
${\cal D}{\cal
M}^\circ_X$, such that the resulting spaces of solutions become
 smooth manifolds which are
Fredholm  of negative index over $Clif^k$. This can be achieved, 
 since ${\cal D}{\cal
M}^\circ_X$ is embedded in the space ${\cal M}^*$, which, though 
possibly singular, is
Fredholm over $Clif^k$.\\

In order to carry out this idea, we will need the following two
 general lemmata.

Let $f$ be  a smooth map taking values in a Banach space,  
and  denote by $Z(f)$ its
vanishing locus. For a point $p\in Z(f)$ define the Zariski tangent
space to $Z(f)$ in $p$ by   
$$T_p(Z(f)):=\ker(d_p f )\ .
$$

\begin{lm} Let $\Sigma$ be a Banach manifold, $p\in\Sigma$, $E$ a
 Banach space, and
$s:\Sigma\map E$ a smooth map such that $s(p)=0$. Suppose\\
i)  $\ker d_p s$ has a topological complement.\\
ii) $\im d_p s$ is closed and has a topological complement.

Then there exists an open  neighbourhood  $\Sigma'$ of $p$ in $\Sigma$ 
and a submanifold $W$ of
$\Sigma$ containing $p$, such that\\
1. $\Sigma'\cap Z(s)$ is a closed subset of $W$.\\
2. $T_p(Z(s))=T_p(W)$.
\end{lm}
\pf Put $T:=\im d_p s$, and denote by $pr_T$ the projection on $T$ associated with a topological
complement of $T$.

The composition $pr_T\circ s$ is a submersion in $p$, since its derivative in $p$ is surjective
and $\ker(d_p(pr_T\circ s))=\ker(d_p  s)$ has a topological complement by assumption.
Let
$\Sigma'$ be an  open neighbourhood of $p$ such that $pr_T\circ s$ is a submersion in every
point of
$\Sigma'$.

Then
$$\Sigma'\cap Z(s)=\Sigma'\cap Z(pr_T\circ s)\cap Z(s)= Z(pr_T\circ
s|_{\Sigma'}))\cap Z(s)\ .
$$
Therefore, taking $W:=Z(pr_T\circ s|_{\Sigma'}))$, claim 1.  follows.  Clearly
$$T_p(W)=\ker(d_p(pr_T\circ s))=\ker(d_p  s)=T_p(Z(s))\ .
$$
\qed
\begin{lm}
Let $W$ be a Banach manifold, $E$ a Banach space, $p\in W$, and 
$\delta:W\map E$ a 
smooth map such that $\ker (d_p \delta)$ has infinite codimension in $T_p(W)$.

Then, for every $n\in \N$ there exists an open neighbourhood
 $W'_n$ of $p$ in $W$ and a
codimension $n$ submanifold $V_n$ of $W$ such that
$W'_n\cap Z(\delta)$ is a closed subset of $V_n$.
\end{lm}

\pf Since  $\ker (\delta_p d)$ has infinite codimension in $T_p(W)$, 
it follows that 
$\im(d_p \delta)$ has infinite dimension. Let $F_n\subset 
\im(d_p \delta)$ be a
subspace of dimension $n$, and $pr_{F_n}$ the projection 
associated with a topological
complement of $F_n$ in $E$. The composition
$pr_{F_n}\circ \delta$ is a submersion in $p$. Indeed, the 
 derivative in $p$ is
surjective and the kernel of the derivative is closed of finite 
codimension, hence it has
a topological complement. Let $W'_n$ be an open  neighbourhood 
of $p$ such that
$pr_{F_n}\circ \delta$ is a submersion in every point of $W'_n$.
Then
$$W'_n\cap Z(\delta)=W'_n\cap Z(pr_{F_n}\circ \delta)
 \cap   Z(\delta)=Z(pr_{F_n}\circ
\delta|_{W'_n})
\cap   Z(\delta)\ .
$$
Take $V_n:=Z(pr_{F_n}\circ \delta|_{W'_n})$.
\qed

\begin{lm}  Every   non-abelian point $[p]\in {\cal D}{\cal M}^*_X$ 
 has a neighbourhood
$U_{[p]}$ which is a closed analytic subspace of a submanifold 
$V_{[p]}\subset \qmod{[{\cal
A}^*_l\times Clif^k]}{{\cal G}_{l+1}}$  such that the 
projection $V_{[p]}\map
Clif^k$ is Fredholm of negative index.
\end{lm}

\pf Put $p=(\theta,\gamma)$ with $\theta\in{\cal A}_l^*$ 
and $\gamma\in Clif^k$. 
Consider
a slice  $S_\theta\subset \theta+\ker(D^0_\theta)^*\subset
  {\cal A}^*_l$ through
$\theta$ to the orbits of the ${\cal G}_{l+1}$-action, such 
that the restriction of the
canonical projection to $S_\theta$   defines a parameterization
 of the quotient
$\qmod{{\cal A}^*_l}{{\cal G}_{l+1}}$ around $[\theta]$. 

Note first, that the image $T$ of the differential 
$d_p({sw|_{S_\theta\times Clif^k}})$ is
closed and has finite codimension in the Hilbert space 
$A^0(\Sigma^-)_{l-1}\times
A^0(\ad_+\otimes\gr_0)_{l-1}$.

Indeed,   $T$ contains the image  of $\frac{\partial sw}{\partial
(A,\Psi)}|_p$ ,  which is the operator $D^1_{\theta}$ associated 
with the deformation
elliptic complex of  the solution $\theta=(A,\Psi)$,  and  the 
image of $D^1_{\theta}$ is
already closed of finite codimension.

Now put $\Sigma:=S_\theta\times Clif^k$, and note that the restriction 
$$q:\Sigma\map\qmod{[{\cal
A}^*_l\times Clif^k]}{{\cal G}_{l+1}}$$
of the canonical projection is a parametrisation of the Banach manifold  $\qmod{[{\cal
A}^*_l\times Clif^k]}{{\cal G}_{l+1}}$ around $[p]$. \\ \\
\ub{Claim:} Put $s:=sw|_\Sigma$. Then the projection 
$$T_\theta S_\theta\times T_\gamma(Clif^k)\supset
\ker(d_p s)\map T_\gamma(Clif^k)$$
is Fredholm. In particular $\ker(d_p s)$ has a topological complement
 in  the tangent space 
  $T_p(\Sigma)=T_\theta S_\theta\times T_\gamma(Clif^k)$.\\

Indeed, the  kernel of this map is $\H^1_\theta$ and its image can 
be identified
with the subspace 
$\left(\frac{\partial {sw}}{\partial \gamma}\right)^{-1}[\im
D^1_{\theta}]$, whose codimension is at most  $\dim \H^2_\theta$. 
   If $\Lambda$ is a
topological complement of $\H^1_\theta$ in $T_\theta S_\theta=
\ker(D^0_\theta)^*$ and
$F$ is a topological complement of $\left(\frac{\partial {sw}}{\partial
\gamma}\right)^{-1}[\im D^1_{\theta}]$ in $T_\gamma(Clif^k)$,
 then
$(\Lambda\times\{0\})\oplus(\{0\}\times F)$ is a topological 
complement of $\ker(d_p s)$
in $T_\theta S_\theta\times T_\gamma(Clif^k)$.
\\

Applying Lemma 3.13 to the Banach
manifold $\Sigma$ and the map $s$, we get a neighbourhood
$\Sigma'$ of $p$ and a submanifold
$W$ such that $\Sigma'\cap Z(s)$ is a  closed subset of $W$ and
$$T_p(W)=T_p(Z(s))\simeq T_{[p]}({\cal M}^*)\ .
$$

The restriction $\det|_W$ of the determinant map
 $\det:\Sigma\map
 A^0(\det(P^u))_{l}$
satisfies  the hypothesis of Lemma 3.14.

Indeed,
$$\ker d_p(\det|_W)=\ker(d_p(\det|_\Sigma))\cap T_p(W)=
\ker(d_p(\det|_\Sigma))\cap\ker
d_p(sw|_\Sigma)\simeq$$
$$\simeq T{[p]}({\cal D}{\cal M}^{*}_X)\ ,
$$
which has infinite codimension in $T_{[p]}({\cal M}^*)\simeq T_p(W)$ 
 by Corollary 3.12.

Using now Lemma 3.14 we get, for any $n\in \N$, an open 
neighbourhood $W'_n$ of $p$ in
$W$ and a codimension $n$ submanifold $V_n$ of $W$ such 
that $W'_n\cap Z(\det|_W)$ is a
closed subspace of $V_n$.

Let $\Sigma'_n\subset\Sigma'$ be an open neighbourhood of $p$ in 
$\Sigma$ such that 
$$W'_n=\Sigma'_n\cap W\ .
$$
Then we have
$$\Sigma'_n \cap q^{-1}(({\cal D}{\cal M}^{*}_X)=   Z(sw|_{\Sigma'_n})\cap
Z(\det|_{\Sigma'_n})= $$
$$
=Z(pr_T\circ sw|_{\Sigma'_n})\cap Z(sw|_{\Sigma'_n})\cap
Z(\det|_{\Sigma'_n})=W'_n\cap Z(sw|_{\Sigma'_n})\cap
Z(\det|_{\Sigma'_n})= $$
$$= [W'_n\cap Z(\det|_{W'_n})]\cap  Z(sw|_{\Sigma'_n})\ . 
$$

Therefore $\Sigma'_n \cap q^{-1}(({\cal D}{\cal M}^{*}_X)$
 is a closed subspace of $
[W'_n\cap Z(\det|_{W'_n})]$, which is closed in $V_n$.

On the other hand we know that the projection 
$$T_p(W)=\ker(d_p s )\map
T_{\gamma} Clif^k$$
 is Fredholm.  Since being Fredholm is an open property, we 
may assume (taking $\Sigma'$  
small) that the projection of $W$ on $Clif^k$ is Fredholm  of constant index.

Now   choose $n$ larger than the index of this projection, and put
$$V_{[p]}:=q(V_n)\ , U_{[p]}:=q(\Sigma'_n \cap q^{-1}
(({\cal D}{\cal M}^{*}_X))=
q(\Sigma'_n)\cap  {\cal D}{\cal M}^{*}_X\ .
$$
\qed
\begin{co}   The set
$$\{\gamma\in Clif^k|\  {\cal D}{\cal M}^*_X
\cap{\rm pr}_{Clif^k}^{-1}(\gamma)\ {\rm contains\
a\   non-abelian\ pair}\} $$
is a set of the first    category  in $Clif^k$.
 
\end{co}
\pf Indeed, let again ${\cal D}{\cal M}^\circ_X$ be the open subspace of 
${\cal D}{\cal M}^*_X$
consisting of solutions with non-abelian $(A,\Psi)$-component. By Lemma 3.15
  and the
Lindel\"of Theorem ([Ke], p. 49) we can find a countable cover $(U_i)_i$ of
${\cal D}{\cal M}^\circ_X$ such that every $U_i$ is a closed analytic subspace 
of a smooth
submanifold
$V_i\subset
\qmod{ [{\cal A}^*_l\times Clif^k]}{{\cal G}_{l+1}}$ which 
projects on the parameter space
$Clif^k$ via a Fredholm map of negative index. Since Fredholm 
maps are locally proper [Sm], it
follows that ${\rm pr}_{Clif^k}({\cal D}{\cal M}^\circ_X)$ is
 a countable union of closed sets;
each of these closed sets is contained in a set of the form 
${\rm pr}_{Clif^k}(V_i)$, which is of
the first category, by the Sard-Smale theorem.
\qed 
 \vspace{1cm}

Corollary 3.12, Lemma 3.15, Corollary 3.16 hold  for 
\ub{every} family of order 0-perturbations
of the equations which contains the perturbations  of the Clifford map 
 which  we have  
studied above.
We need the following particular case:

Define the space of parameters  ${\cal P}^k$ by 
$${\cal P}^k:={\cal C}^k(P^u\times_\pi\C^4)
\times {\cal C}^k(GL(\ad_+))\times
Clif^k\ .
$$
Recall that a section $\beta$ in the bundle
$$P^u\times_\pi\C^4=\Hom(\Sigma^+_c,\Sigma^-_c)\subset
 \Hom(\Sigma^+,\Sigma^-)$$
defines an order 0-operator $A^0(\Sigma^+)\map  A^0(\Sigma^-)$,
commuting with the gauge action.

Consider now  the equations
$$\left\{\begin{array}{lcc}
\Dr_{\gamma,a,A}^0\Psi+ \beta (\Psi)&=&0\\
\Gamma_\gamma(F_A)&=&K (\Psi\bar\Psi)_0
\end{array}\right. \eqno{(\tilde {SW_a})}
$$
for a system 
$$(A,\Psi,\beta,K,\gamma)\in \tilde {\cal A}:=
{\cal A}(\bar\delta(P^u))_l\times
A^0(\Sigma^+)_l\times {\cal P}^k\ .$$  

Let $[{\cal A}_l\times{\cal P}^k]^{\tilde {SW}_a}$ 
($[{\cal A}_l\times{\cal P}^k]^{\tilde
{SW}_a}_U$) be the space of solutions of the equations
$(\tilde {SW}_a)$ (whose spinor component is degenerate on $U$), 
and denote also   by 
$\widetilde {\cal M}^*$ ($\widetilde{{\cal D}{\cal M}}^*_U$) 
 the moduli space of solutions
(whose spinor component is degenerate on $U$) with
 non-vanishing spinor component.

\begin{pr} Let $p=(A,\Psi,\beta,K,\gamma)\in
 [{\cal A}_l\times{\cal P}^k]^{\tilde
{SW}_a}_U$  such that for some $u\in U$, $\Psi\ne 0$ and
 $(A,\Psi)$ non-abelian in $u$.
Then the  Zariski tangent  space  
$T_{[p]}\widetilde{{\cal D}{\cal M}}^*_U$ has
 infinite codimension  in
$T_{[p]}\widetilde {\cal M}^*$.
\end{pr}

\pf   
Consider  the image of  
$T_p([{\cal A}_l\times{\cal P}^k]^{\tilde {SW}_a}_U)$ under the
projection to the tangent space
$T_{(\beta,K,\gamma)}{\cal P}^k\ .$
This image has again infinite codimension. To see
 this it is enough to notice that the
intersection of this image with the subspace 
$\{0\}\times\{0\}\times T_\gamma Clif^k$
has infinite codimension in 
$\{0\}\times\{0\}\times T_\gamma Clif^k$. But 
this follows  by precisely the same
arguments as in Theorem  3.7; one just has  to replace
 the equations $(SW_a)$ by their
$(\beta,K)$-perturbations. The left hand  side in the
 crucial identity (6) will only be
modified by the 0-order term $\dot\zeta\beta(s_0^+)$.
\qed

Using this result and the same arguments as above, we get
\begin{co} The set
$$\{\pg\in{\cal P}^k|\  \widetilde{{\cal D}{\cal M}}^*_X
\cap{\rm pr}_{{\cal P}^k}^{-1}(\pg)\ {\rm
contains\ a\ non-abelian\ pair}\} $$
is a set of the first  category in ${\cal P}^k$ .
\end{co}
We can state now our generic regularity result:
\begin{thry} There is a dense second category subset ${\cal P}^k_0$ of 
${\cal P}^k$ such that for
every
$\pg\in{\cal P}^k_0$ the moduli space ${\cal M}^*_\pg:=
\widetilde{\cal M}^*\cap{\rm pr}_{{\cal
P}^k}^{-1}(\pg)$ is smooth away from the abelian locus.
\end{thry}

\pf We know by Theorem 3.1 and  Remark 3.2  that 
$\tilde{\cal M}^*\setminus \widetilde
{{\cal D}{\cal M}}^*_X$ is a smooth manifold. Applying the  
Sard-Smale theorem   to
the Fredholm map
$$\widetilde {\cal M}^*\setminus \widetilde{{\cal D}{\cal M}}^*_X\map{\cal P}^k$$
 it follows
that there exists a first category subset   
${\cal P}^k_1\subset{\cal P}^k$ such that the moduli space 
$[\widetilde {\cal M}^*\setminus
\widetilde{{\cal D} {\cal M}}^*_X]\cap{\rm pr}_{{\cal P}^k}^{-1}(\pg)$
 is smooth for every
$\pg\in{\cal P}^k\setminus{\cal P}^k_1$.  Let ${\cal P}^2_k$ be
 the first category set given
by Corollary 3.18, and take ${\cal P}^k_0:={\cal P}^k
\setminus({\cal P}_k^1\union {\cal
P}_k^2)$.
\qed
\vspace{1cm}

Finally consider the following parameterized
ASD-$Spin^c$- equations
$$\left\{\begin{array}{lcc}
\Dr^0_{\gamma,a,A}\Psi+ \beta(\Psi)&=&0\\
\Gamma_\gamma(F_A)&=&0\end{array}\right.
$$
for a system $(A,\Psi,\beta,\gamma)\in {\cal A}(\bar\delta(P^u))_l\times
A^0(\Sigma^+)_l\times{\cal C}^k(P^u\times_\pi\C^4)\times Clif^k$.

Let ${\cal M}'^*$ be the moduli space of
solutions with non-trivial spinor component, and let ${\cal P}'^k$
  be the parameter space ${\cal
P}'^k:={\cal C}^k(P^u\times_\pi\C^4)\times Clif^k$. Denote
 also by ${\cal D}{\cal M}'^*_X$ the
subspace of solutions with degenerate spinor component, and by ${\cal
M}'^*_{\rm red}$ the subspace of solution with reducible
 connection-component.

Using the methods of section 3.2, one can prove the following 
partial transversality result
\begin{pr} Suppose that the base manifold is simply connected. Then 
the moduli space ${\cal
M}'^*$ is smooth away from the union ${\cal D}{\cal M}'^*_X\union{\cal
M}'^*_{\rm red}$.
\end{pr}

\pf Indeed,    let $ p=(A,\Psi,\beta,\gamma)$ be a solution with non-degenerate
 spinor
component and non-reducible connection component, and suppose as in the 
proof of Theorem 3.1
that $(\Phi,S)$ is $L^2_{g_\gamma}$-orthogonal on
 the image of the differential in $p$ of the
map cutting out the space of solutions. Using variations 
$\dot \beta$ of $\beta$ one sees that
$\Phi$ must vanish on a non-empty open set. But using 
variations of $\Psi$, it follows that
$\Phi$ must solve a Dirac equation, hence by Aronszajin's unique 
continuation theorem, it must
vanish on $X$. Then   using variations $\dot \gamma$ of $\gamma$ 
we get as in [DK], p. 154
that $S=0$. It is enough to notice that $A$ is 
  $g_\gamma$-ASD, and that any variation of the
metric  $g_\gamma$ is induced by a variation of the Clifford map $\gamma$.
\qed 

 In the proof of Theorem
3.7  we have only used the Dirac equation and the
 ellipticity (modulo the gauge group) of the
system .  Therefore the same arguments as above give the  following important  
\begin{thry} \hfill{\break}
1. There exists a first category subset ${\cal P}'^k_2\subset{\cal
P}'^k$ such that for every $\pg\in {\cal P}'^k\setminus 
{\cal P}'^k_2$ the only    solutions
with degenerate spinor component in the moduli space
${\cal M}'^*\cap p_{{\cal P}'^k}^{-1}(\pg)$ are the abelian ones.\\
2. If the base manifold $X$ is simply connected, there
 exists a dense second category subset
${\cal P}'^k_0\subset{\cal P}'^k$ such that for every $\pg\in{\cal P}'^k_0$  the $Spin^c$-moduli
space
${\cal M}'^*\cap p_{{\cal P}'^k}^{-1}(\pg)$ is
 smooth away from ${\cal M}'^*_{\rm red}\cap
p_{{\cal P}'^k}^{-1}(\pg)$.
\end{thry}
 The results above are sufficient to go forward towards a
 complete proof of the
Witten conjecture.

Moreover, one can use the same method to prove a generic
 regularity theorem along the
abelian part of the moduli space.  

More precisely, let ${\cal M}^{\rm ab}_\pg\subset {\cal M}^*_\pg$

be the abelian part of the moduli space ${\cal M}^*_\pg$ of solutions of the monopole
equations associated with the perturbation parameter  $\pg$. 
The space ${\cal M}^{\rm
ab}_\pg$ can be identified with the disjoint union of the 
$Spin^c$-Seiberg-Witten moduli
spaces associated with the abelian reductions of $P^u$ ([OT5], [OT7], [T1]). 

Let $[p]\in  {\cal M}^{\rm ab}_\pg$ be an abelian solution. 
The elliptic deformation  complex
${\cal C}_p$ of $p$ splits as the sum
$$ {\cal C}_p={\cal C}_p^{\rm ab}\oplus {\cal N}_p 
$$
where the first summand ${\cal C}_p^{\rm ab}$ can be
 identified with the
elliptic deformation complex of $p$ regarded as solution of 
the abelian monopole equations,
and ${\cal N}_p$ is the  so called normal elliptic complex of $p$.

The union ${\cal H}^1_\pg:=\union\limits_{[p]
\in {\cal M}^{\rm ab}_\pg}\H^1({\cal N}_p)$ is a
real analytic space which fibres over ${\cal M}^{\rm ab}_\pg$,
 but in general is not locally
trivial  over ${\cal M}^{\rm ab}_\pg$,  and local triviality
 cannot be achieved in the class of
$S^1$-equivariant perturbations.

Using the   method from above one can prove
\begin{pr} For a generic parameter $\pg\in{\cal P}^k$, 
the complement of the zero section in
${\cal H}^1_\pg$ is smooth of the expected dimension in every point.
\end{pr}

\section{The Uhlenbeck Compactification}
\subsection{ Local estimates}

The essential difference between the anti-self-dual and the 
    monopole equations is that the
latter are  not conformal invariant.  Under a conformal rescaling
  of a metric 
$g\mapsto\tilde g=\rho^2 g$ on a 4-manifold $X$, the associated
 objects change as
follows 
$$\begin{array}{c} 
\tilde{g^*}=\rho^{-2}g^*\ {\rm on\ 1-forms} ;  \ \
vol_{\tilde g}=\rho^4 vol_g \ ;\ \ s_{\tilde g} =\rho^{-2}s_g+
2\rho^{-2}\Delta \rho\\
\Sigma^{\pm}_{\tilde g}=\Sigma^{\pm}_{g}\ ({\rm as\ Hermitian\ bundles}) ,  \
\tilde\gamma=\rho^{-1}\gamma\ ;\
\tilde\Gamma=\rho^{-2}\Gamma\ ;
\
\ 
 
\Dr_{\tilde g}=\rho^{-\frac{5}{2}}\Dr_g\rho^{\frac{3}{2}} \ .

\end{array}
$$

A standard procedure used in proving  regularity and 
compactness theorems  for
instantons is the following: restrict the equations on small balls in the base manifold, and
then rescale the metric. In this way, using the conformal invariance of the equations,  one
can reduce the local computations to the  unit ball endowed with a metric close to the
euclidean one.

A similar procedure will be used in the case of $PU(2)$-monopoles.
 The problem here is
that the perturbed equations depend on a much larger system of parameters
 (data). Using
\underbar{constant} rescalings of the Clifford map (and hence of the metric), 
we show
first that one can reduce the local computations to computations on the unit ball endowed
with a system close to a system of "standard data" (see Definition 4.4).

First of all notice that  if $(A,\Psi)\in {\cal A}(\bar\delta(P^u))
\times A^0(\Sigma^+)$ is a
solution of the non-perturbed $PU(2)$-monopole equations
 $SW^\sigma_a$ for the metric
$g$ with respect to the $Spin^{U(2)}(4)$-structure $\sigma$, 
and   if $\rho$ is a
\ub{constant}, then $(A,\rho^{-1}\Psi)$ is a solution of the monopole
equations $SW^{\tilde\sigma}_a$ for $\tilde g=\rho^2 g$ with 
respect to   the 
$Spin^{U(2)}(4)$-structure $\tilde\sigma$ defined by the 
correspondingly rescaled Clifford
map $\tilde\gamma=\rho^{-1}\gamma$.

The case of the perturbed equations is more delicate.  Fix a 
$Spin^{U(2)}(4)$-bundle $P^u$. 
To write down the general perturbed $PU(2)$-monopole 
equations we considered, one also
needs {\it a system of data} of the form $\pg=(\gamma,C,a,\beta,K)$, where $\gamma$
is a Clifford map (see Definition 3.3), $C$ is an $SO(4)$-connection in
$P^u\times_\pi\R^4$, $a$ is a connection in the line bundle $\det(P^u)$,    
$\beta$ is a section in $P^u\times_\pi\C^4$, and  $K$ is a section  in
$\End(\ad_+)$.  

The rescaling rule is:
\begin{re} If $(A,\Psi)\in{\cal A}(\bar\delta(P^u))\times A^0(\Sigma^+)$  
  solves the
perturbed $PU(2)$-monopole equations associated with the data 
$(\gamma,C,a,\beta,K)$.
Then $(A,\rho^{-1}\Psi)$ solves the perturbed $PU(2)$-monopole 
equations associated with
the data $$(\rho^{-1}\gamma,C,a,\rho^{-1}\beta,K).$$
\end{re}

Let $\bar B$ be the standard closed 4-ball  with interior $B$.  
Fix two copies $\H_{\pm}$ of the quaternionic skew-field $\H$ regarded
 as right complex
and quaternionic vector spaces and consider the two trivial $SU(2)$-bundles
$S_0^{\pm}:=\bar B\times\H_{\pm}$. Let also $E_0=\bar B\times\C^2$ be
the trivial Hermitian rank 2-vector bundle on $\bar B$.   

Let $P^u_0$ be the trivial $Spin^{U(2)}(4)$-bundle associated with 
$S_0^{\pm}$, $E_0$ via
the   morphism $SU(2)\times SU(2)\times  U(2)\rightarrow  
Spin^{U(2)}(4)$ (section 2.1,
Prop. 2.2).  

A Clifford map for $P^u$ is an orientation preserving linear 
isomorphism
$\gamma:\Lambda^1_{\bar B}\map \Hom_\H(S^+_0,S^-_0)=
\bar B\times\H$. To every such a
 Clifford map $\gamma$, we can associate the constant Clifford $\gamma^c$ given by the
composition 
$$\Lambda^1_{\bar B}\map\bar B\times
\Lambda^1_0\textmap{\id\times{\gamma|_{\Lambda^1_0}} } B\times\H  $$
 Note that the corresponding metric $g_{\gamma^c}$ is flat.

Denote by $ h_r:\bar B\map \bar B_r\subset\bar B$    the homothety of slope $r<1$. 
\begin{re}\hfill{\break}
  The Clifford maps $\gamma_r:=r h_r^*(\gamma|_{B_r})$
 converge in the ${\cal
C}^{\infty}$-topology to $\gamma_0$, which is a Clifford map  for the flat
metric $g_{\gamma^c}$. In particular the metrics 
$g_r:=r^{-2}h_r^*(g)$ converge to the flat
metric $g_{\gamma^c}$.
\end{re}
Indeed, one has 
$$\gamma_r(x,\lambda)=r \gamma((h_r)_*(x,\lambda))=r
\gamma(rx,r^{-1}\lambda)=\gamma(rx,\lambda)$$

The data of a $PU(2)$-connection $A\in{\cal A}(\bar\delta(P^u_0))$ 
is equivalent to the
data of a connection matrix, i.e. an element in $A^1(\bar B,su(2))$. 
Similarly, the data of a
$U(1)$-connection in $\det(P^u_0)$ is equivalent to the data of a 1-form in $A^1(\bar
B,u(1))$.

\begin{re} Let $(A,\Psi)\in{\cal A}(\bar\delta(P^u_0))\times 
A^0(\Sigma^+(P^u_0))$ be a pair
 which solves the monopole equations for  the data
$(\gamma,C,a,\beta,K)$. Then
$(h_r^*(A), r h_r^*(\Psi))$ solves the $PU(2)$-monopole
 equations for the  data  
$(\gamma_r,h_r^*(C),h_r^*(a),rh_r^*(\beta),h_r^*(K))$.  
\end{re} 

Note that, as $r\rightarrow 0$\ ,\\
 $\gamma_r\rightarrow \gamma^c$ (which is a Clifford map
 for the flat metric
$g_{\gamma^c}$), \\
$rh_r^*(\beta)\rightarrow 0$, $h_r^*(K)\rightarrow K(0)$,\\
$h_r^*(a)$ converges to the flat connection in $B\times\C=\det(P^u_0)$, \\
$h_r^*(C)$ converges to the flat connection in $B\times\H$, and \\
$\gamma^{-1}_r(h_r^*(C))$ converges to the flat connection in
$(\Lambda^1_B=B\times\R^4,g_{\gamma^c})$, which is 
precisely the Levi-Civita connection
for $g_{\gamma^c}$).

\begin{dt} A system  of data for the bundle $P^u_0$ will 
 be called a standard system, if it
has the form $(\gamma_0, C_0,0,0,K_0)$, where:\\
$\gamma_0$ is the standard identification $\Lambda^1_{B}
=B\times\R^4\map
B\times\H$, \\
 $C_0$ the flat $SO(4)$-connection in $B\times\H$,  and \\
$K_0$ is a \ub{constant} automorphism of the trivial bundle 
$su(S^+_0)=B\times
su(2)_+$. 
\end{dt}
The  metric associated with the standard identification
$\Lambda^1_{B}=B\times\R^4\rightarrow B\times\H$ is the standard 
Euclidean metric
$g_0$ on the ball.

For any  $K_0\in \End(su(2))$, let $\pg_{K_0}$ be the standard 
system of  data on $\bar B$
defined by $K_0$.\\

Let $X$ now be 4-manifold, and $P^u$ a $Spin^{U(2)}(4)$-bundle on it. 
 Let
$x_0$ be a point in $X$ and $U$ an open  neighbourhood of $x_0$. Fix an
 identification
of $P^u|_U$ with the the trivial $Spin^{U(2)}(4)$-bundle on $U$, i. e. 
with the
$Spin^{U(2)}(4)$-bundle associated with the triple $U\times\H_\pm$,
$U\times \C^2$ (see section 2.1). 

Given  a system of data $(\gamma,C,a,\beta,K)$ for $P^u$,
we consider a parameterization $B_{r_0}\stackrel{f}{\map} U\subset  X$ 
around $x_0$ such
that $f(0)=x_0$ and  $\gamma|_{\Lambda^1_{x_0}}\circ 
[f_*]_{\Lambda^1_0}$ is the
standard identification
$\Lambda^1_0=\R^4\map  \H$.  
\begin{re}
For any pair $(A,\Psi)$ solving the monopole equations for the data
$(\gamma,C,a,\beta,K)$, the pair $\left((f\circ h_r)^*(A),
r(f\circ h_r)^*(\Psi)\right)$
solve the monopole equations associated with  the system 
$$\left(f^*(\gamma)_r,(f\circ h_r)^*(C),(f\circ h_r)^*(a),
r(f\circ h_r)^*(\beta),(f\circ h_r)^*(K)\right)\ .$$
 This system converges to a system of standard data on the ball, as $r\rightarrow 0$.
\end{re}

Therefore, as long as we are interested only in local computations, we can work on the
standard ball and  assume (via the transformation defined in Remark 4.5) that our
system of data belongs to a small neighbourhood of a standard system.\\

We recall now the following important "gauge fixing" theorem 
(see Theorem 2.3.7 in [DK]). 

\begin{thry} (Gauge-fixing) There are constants
$\varepsilon_1,\ M>0$ such that  the following holds: 

Any connection $A$ on the trivial bundle $E_0$ over $\bar B$ with
$\parallel F_A\parallel_{L^2}<\varepsilon_1$ is gauge equivalent to a 
connection $\tilde A$
over $B$ with\\ 
(i) $d^*_0\tilde A=0$,  where $d^*_0$ is the normal adjoint of $d$
 with respect to the
standard flat metric $g_0$.\\
 (ii) $\lim_{r\rightarrow 1} A_r=0$ on $S^{3}$,\\  
(iii) $\parallel\tilde A\parallel_{L^2_1}\leq M\parallel F_{\tilde
A}\parallel_{L^2}$. \\
 The corresponding  gauge transformation is unique up to a constant matrix.
\end{thry}
Using this result we can prove the following
\begin{thry} (Local estimates for data close to the standard data)  There
is a positive constant $\varepsilon_2=\varepsilon_2(K_0)>0$  such that for
 any system  of
data $\pg'$ on $\bar B$ which is sufficiently ${\cal C}^2$-close to the standard
system $\pg_{K_0}$,    the following holds:

For any solution $(A,\Psi)$  of the $PU(2)$-monopole equation for
the   monopole equations associated with $\pg$ over the open ball $B$ 
 satisfying the
conditions $d^*_0 A=0$, 
$\nr(A,\Psi)\nr_{L^4}\leq\varepsilon_2$, and any interior domain 
$D\subsetint B$, one has  
estimates of the form : 
$$\nr (A,\Psi)\nr _{L^2_l(D)}\leq
C_{D,l,\pg'} \nr (A,\Psi)\nr_{L^4}\ ,\\
$$
with  positive  constants  $C_{D,l,\pg'}$,  for all $l\geq 1$.
\end{thry}
\pf  First of all we identify the ball with the upper semi-sphere of
$S:=S^4$ and we endow the sphere with a metric $g_s$ which extends
 the standard flat
metric $g_0$ on the ball, and which has non-negative sectional 
curvature\footnote{Such
a metric can be obtained as follows: consider a plane {\it convex} curve
 with
symmetry  axis $Oy$, which is horizontal in a neighbourhood of its upper
 intersection point
with $Oy$. Then  rotate this curve around the $Oy$-axis in the 
5-dimensional space
{\tenmsb R}$^4\times Oy$ . The hypersurface obtained in this 
way is also conformally
flat, by a theorem of E. Cartan (see [GHL], p. 157,  [Ch],  Th. 4.2, 
p. 162)}.

We fix  a   $Spin(4)$-structure on the sphere with spinor bundles
$S^{\pm}_s$ given by a Clifford map
$\gamma_s:\Lambda^1_S\map \Hom_\H(S_s^+ , S^-_s)$, 
 which, with respect to  fixed
trivializations $S^{\pm}|_{\bar B}=\bar B\times \H_{\pm}$, 
 extends the standard Clifford map
$\gamma_0$ on the ball. Let also $C_s$ be the Levi-Civita 
connection induced
by $\gamma_s$ in  $\Hom_\H(S_s^+ , S^-_s)$. Its restriction to the
 ball is the standard flat
connection $C_0$ in $\bar B\times \H$.  Let finally $K_s$ be an   extension 
of $K_0$ to
an endomorphism $K_s\in A^0(\End(su(S_s^+)))$.

We denote by $E_s$ the trivial $U(2)$-bundle over $S$, and by $P^u_s$ the
$Spin^{U(2}(4)$-associated with the triple $(S^{\pm}_s,E_s)$. 
$P^u_s$ comes with an
identification $P^u_s|_{\bar B}=P^u_0$, induced by the fixed trivializations
 of $S^{\pm}_s$.

The system $(\gamma_s,C_s,0,0,K_s)$ is an extension on the sphere of 
the standard
system  $\pg_s:=(\gamma_0,C_0,0,0,K_0)$. The point  is now that any 
system $\pg'$ of data
which is close to $\pg_{K_0}$ has an extension $\pg$ which is close to 
$\pg_s$.

Put $\pg=(\gamma,C,a,\beta,K)=(\qg,K)$. The system $\qg$ defines two
 first order elliptic
operators on the sphere
$$\begin{array}{ccccc}
\Dr_{\qg}&:& A^0(S^+_s\otimes E_s)& \map&
A^0(S^-_s\otimes E_s)\\
 \delta_\gamma:=d_s^*+ \Gamma_\gamma\circ d &: &A^1(su(2))&
\map& A^0(su(2))^\bot
\oplus A^0(su(S^+_s) \otimes su(E_s)) 
\end{array}$$
The symbol $d^*_s$ means the adjoint of 
$d:A^0(su(2))\map 
A^1(su(2))$ with respect to the fixed   metric $g_s$, and
$\Dr_{\qg}:=\Dr_\gamma^C+\beta+\gamma(\frac{a}{2})$.  
 $A^0(su(2))^\bot$ denotes
the $L^2_{g_s}$-orthogonal complement of the 3-dimensional 
space of constant sections.

These operators are  injective in the special case $\qg=\qg_s:=
(\gamma_s,C_s,0,0)$, by the
  Weitzenb\"ock formula for the Dirac operator and because the cohomology group
$H^1_{\rm DR}(S)$  vanishes. Since   the coefficients of both operators in local coordinates
are algebraic expressions in the components of
$\qg$, it follows by
elliptic semicontinuity that the two operators remain injective if 
$\qg$ is sufficiently
${\cal C}^0$-close to $\qg_s$. Denote by $D_\qg$ the direct
sum of these operators. We   get operator valued maps
$$ 
\qg\mapsto D_{\qg}\in $$ $$\Iso\left[A^0(S^+_s\otimes
E_s\oplus \Lambda^1(su(2)))_{k+1}, A^0(S^-_s\otimes E_s 
\oplus su(S^+_s)\otimes
su(E_s)_k\oplus A^0(su(2))_k^\bot\right] $$
which are continuous with respect the ${\cal C}^k$-topology on 
the space of data $\qg$ on the
sphere.

Therefore one has elliptic estimates
$$ 
\nr u\nr_{L^2_{k+1}}\leq const(\qg)\nr D_\qg u\nr_{L^2_{k}}\\
\eqno{(el_k)}$$
where $const(\qg)$  depends continuously on $\qg$
w. r. t. the ${\cal C}^k$-topology. In a sufficiently small 
${\cal C}^2$-neighbourhood of
$\qg_s$ one has the following estimates with $\qg$-\ub{independent} 
\ub{constants}
$$
\nr u\nr_{L^2_{k+1}}\leq const\nr D_{\qg}
u\nr_{L^2_{k}}\\
\eqno{(el)}$$

Since $D_\qg$ is a first order operator, we have an identity of the form:
$$D_\qg(\varphi v)=\varphi D_\qg (v)+  A_{\qg,\partial\varphi}(v) \eqno{(*)}
$$
where $A_{\qg,\partial\varphi}$ is an operator of order 0 
depending on $\qg$ and depending
linearly on the first order derivatives of $\varphi$.

The first step is  an input-estimate for the ${L^2_1(D)}$-norms:

Denote by $u$ the pair $(A,\Psi)$. Let $\varphi_1$ be a cut-off
 function supported in the open
ball
$B$ which is identically 1 in a neighbourhood of $\bar D$. Then $u_1:=
\varphi_1 u$
extends as section in the   bundle $\Lambda^1(su(2))\oplus S^+_s\otimes E_s$
 on the sphere. 

Taking into account that $u$ solves the monopole equations associated 
with the data $\pg'$, its
connection component is in Coulomb gauge, and that $\pg=(\qg,K)$ 
extends $\pg'$ one gets by
$(*)$
$$D_\qg(u_1)=A_{q,\partial\varphi_1}(u)+\varphi_1
\left[ \matrix{ -\gamma(A)\Psi\cr
-\Gamma_\gamma(A\wedge A)+K(\Psi\bar\Psi)_0}\right]=
$$
$$=A_{q,\partial\varphi_1}(u)+\varphi_1 B_{\gamma,K}(u) \eqno{(1)}
$$
where $B_{\gamma,K}$ is a quadratic map.

Then by $(el)$ we obtain an elliptic estimate of
the form
$$\ 
\nr u_1\nr_{L^2_1}\leq
c\nr D_{\qg} u_1\nr_{L^2} 
\leq c'(\nr u\nr_{L^4}^2+\nr
d\varphi\nr_{L^4}\nr u\nr_{L^4})\leq 
 $$
$$
\leq c''(\nr u\nr_{L^4} \nr u\nr_{L^2_1}+\nr
d\varphi\nr_{L^4}\nr u\nr_{L^4})$$
where, for the second inequality we have used on the right the    bounded 
Sobolev
embedding $L^2_1\subset L^4$. The constants $c$, $c'$
can be chosen to depend continuously on $\pg$, so that we can assume that 
they are
independent of $\pg$ on a small neighbourhood of $\pg_s$. We use now the
 standard
rearrangement procedure described in [DK], p. 60, 62.  For a  sufficiently
 small (independent of
$D$) apriori bound  $\varepsilon(K_0)$  of the norm $\nr u\nr_{L^4}$,  
  we get an estimate of
the type
$$\nr\ u_1\nr_{L^2_1}\leq const_D\nr u\nr_{L^4} \ .
$$
The constant $const_D$ in this estimate is independent of  $\pg$  in a 
sufficiently
small neighbourhood of $\pg_s$,  but it depends on $D$ via 
$\nr d\varphi_1\nr_{L^4}$.\\

In a next step we estimate the $L^2_2$-norms:

Put $u_2=\varphi_2 u$,  where $\varphi_2$ is 
identically 1 on $D$, but the support ${\rm supp}\varphi_2$
 is contained in the
interior of $\varphi_1^{-1}(1)$. Then we can also
write $u_2=\varphi_2 u_1$, and we have
$A_{\qg,\partial\varphi_2}(u)=A_{\qg,\partial\varphi_2}(u_1)$.

We estimate first the $L^2_1$-norm of the right hand side of the
formula obtained by replacing $\varphi_1$ with $\varphi_2$ in (1) . We find
$$\nr D_\qg(u_2)\nr_{L^2_1}\leq const  \nr\varphi_2
 B_{\gamma,K}(u_1)\nr_{L^2_1} 
+const_{D}\nr u_1\nr_{L^2_1}\ ,\eqno{(2)}$$
and again we can assume that the constants do not depend on 
$\qg$.  The term
$\varphi_2B_{\gamma,K}(u_1)$ can be written as 
$\tilde B_{\gamma,K}(\varphi_2
u_1\otimes u_1 )$, where $\tilde B_{\gamma,K}$ is the 
linear map defined on the tensor product
$\left(\Lambda^1(su(2))\oplus S^+_s\otimes E_s\right)^{\otimes 2}$ 
associated with the  
quadratic map $B_{\gamma,K}$.

In local coordinates we can write:
$$\partial_i [\tilde B_{\gamma,K}(\varphi_2 u_1\otimes 
u_1) ]=\partial_i(\tilde B_{\gamma,K}) 
(\varphi_2\otimes
u_1)\otimes u_1 +  \tilde B_{\gamma,K}
 \left[\partial_i(\varphi_2 u_1)\otimes u_1+
u_1\otimes(\varphi_2\partial_i u_1)\right] 
$$
$$=\partial_i(\tilde B_{\gamma,K}) (\varphi_2\otimes
u_1)\otimes u_1+\tilde B_{\gamma,K} 
\left[\partial_i(\varphi_2 u_1)\otimes
u_1+u_1\otimes\partial_i  (\varphi_2 u_1)-
\partial_i(\varphi_2) u_1\otimes u_1\right]
$$
This gives an estimate of the form
$$ 
\nr \tilde B_{\gamma,K}(\varphi_2 u_1\otimes 
u_1)  \nr_{L^2_1}  \leq const \nr  u_2\nr_{L^4_1} 
\nr u_1\nr_{L^4} +const_D \nr
u_1\nr_{L^4}^4\ ,
 $$
which together with (2) and $(el)$ gives 
$$\nr u_2\nr_{L^2_2}\leq const \nr  u_2\nr_{L^4_1}
 \nr u_1\nr_{L^4} +const_D(\nr
u_1\nr_{L^2_1}+\nr u_1\nr_{L^4})\ .
$$
By the same rearrangement argument and using the existence of a bounded 
inclusion
$L^2_2\subset L^4_1$, we get, for a sufficiently small, 
independent of $D$, apriori bound of
$\nr u\nr_{L^4}$, an estimate of the form 
$$\nr u_2\nr_{L^2_2}\leq const_D \nr u\nr_{L^4}\ .
$$

The estimates for the third step can be proved by  the same 
algorithm, using the existence of
a bounded inclusion $L^2_3\subset L^4_2$. 

Since $L^2_3$ is already a Banach algebra, the
estimates for the higher Sobolev norms follow by the usual
 bootstrapping procedure using the
estimates $(el_k)$. Note in particular that we no longer need
 to use the rearrangement
argument, so we do \underbar{not}   have  to take  smaller bounds 
for $\nr u\nr_{L^4}$ to get
estimates of the higher Sobolev norms, so that a positive number
$\varepsilon_2=\varepsilon(K_0)$ (independent of
$l$  and $D$ !) with the required  property   does exist.
\qed

Let $V_+$, $F$ Hermitian vector spaces of rang 2. One can easily
 check that there exists a
universal constants $\eg>0$, $C$, $C_1>0$,  $C_2>0$ such that
 for every
$K\in\End(su(V_+))$ with
$|K-\id|<\eg$, and every $\Psi\in V_+\otimes F$ the following 
inequalities hold 
$$C_1|\Psi|^2\leq |K(\Psi\bar\Psi)_0|\leq  C_2|\Psi|^2 \eqno{(3)}
$$
$$C|\Psi|^4\leq \left(K(\Psi\bar\Psi)_0,
(\Psi\bar\Psi)_0\right)=
\left(K(\Psi\bar\Psi)_0(\Psi),\Psi\right) 
\eqno{(4)}
$$

From now on we'll \ub{always} assume the last component $K$ of a system of data
$(\gamma,C_0,a,\beta,K)$ satisfies in every point $x$ the inequality 
$|K(x)-\id_{\ad_+}|<\eg$.

\begin{co} (Estimates in terms of the curvature) There exists a constant 
$\varepsilon>0$, such
that for any system $\pg'$ of data on the closed ball which is sufficiently 
${\cal
C}^2$-close to a system of standard data $(\gamma_0,C_0,0,0,K_0)$ 
with $|K_0-\id|<\eg$
the following holds:

 For  any interior ball
$D\subsetint B$ and any
$l\geq 1$ there exist  a  positive constants $C_{D,l,\pg'}$, $C'_{D,l,\pg'}$ 
such that every solution
$(A,\Psi)$ of the $PU(2)$-monopole equations on $\bar B$   associated with $\pg'$ 
 satisfying $\nr F_A\nr_{L^2}\leq\varepsilon$, is gauge equivalent on $B$ to a pair
 $(\tilde
A,\tilde\Psi)$ satisfying the estimates
$$\nr\tilde A\nr_{L^2_l(D)}\leq C_{D,l,\pg'}\nr
 F_A\nr_{L^2}\ ,\ \ \nr\tilde
\Psi\nr_{L^2_l(D)}\leq C'_{D,l,\pg'}\nr F_A
\nr_{L^2}^{\frac{1}{2}} \ .
$$
\end{co}
\pf  Note first that all the pairs $(A,\Psi)$ with $\nr F_A\nr_{L^2}\leq 
\varepsilon_1$ are
gauge equivalent to pair $(\tilde A,\tilde\Psi)$ whose connection component
 is in the Coulomb
gauge with respect to the trivial connection and such that 
$$\nr\tilde A \nr_{L^2_1}\leq M\nr
F_{\tilde A}\nr_{L^2}\eqno{(5)}$$

  Since now the constant $K_0$ is supposed to belong to the bounded
set $B(\id,\eg)$ the conclusion of Theorem 4.7 holds for a constant 
$\varepsilon_2$
which can be chosen {\it independently} of $K_0$.  

On the other hand, by the estimate (3) and the second monopole equation, 
one has
$$\nr \tilde\Psi\nr_{L^4_{g_{\gamma'}}}\leq \frac{1}
{C_1^{\frac{1}{2}}} \nr
\Gamma_{\gamma'}(F_{\tilde A})
\nr_{L^2_{g_{\gamma'}}}^{\frac{1}{2}}
= \frac{\sqrt 2}{C_1^{\frac{1}{2}}} \nr F_{\tilde
A} ^{+_{g_{\gamma'}}}\nr_{L^2_{g_{\gamma'}}}^{\frac{1}{2}}
\leq\frac{\sqrt 2}{C_1^{\frac{1}{2}}} \nr F_{\tilde
A}  \nr_{L^2_{g_{\gamma'}}}^{\frac{1}{2}}
\eqno{(6)}
$$
Since $\gamma'$ is supposed to belong to a small
 neighbourhood 
of $\gamma_0$ this gives
an uniform estimate of $\nr \tilde\Psi\nr_{L^4}$ in terms of $\nr
  F_{\tilde
A} \nr_{L^2}^{\frac{1}{2}}$. Using the bounded inclusion
$L^2_1\subset L^4$, and the estimates (5), (6) we see now that the
$L^4$ norm of the pair $(\tilde A,\tilde\Psi)$ can be made as small 
as we please by choosing
$\varepsilon$ small, in particular smaller than the constant 
$\varepsilon_2$.  With
this choice the conclusion of Theorem 4.7 holds, and we get 
estimates of the Sobolev
norms of the restrictions on smaller disks $D\subsetint B$ in
 terms of $\nr(\tilde
A,\tilde\Psi)\nr_{L^4}$, hence in terms of $\nr 
 F_{\tilde A} \nr_{L^2}^{\frac{1}{2}}$.  

On the other hand,   the same cutting off procedure as in the 
proof of Theorem 4.7, gives on the
sphere an identity of the form  
$$(d_s^*+\Gamma_\gamma d)(\varphi_1\tilde
A)=A'_{\qg,\partial\varphi }(\tilde A)+\varphi [-\Gamma_\gamma
(\tilde A\wedge
\tilde A)+(\tilde \Psi\bar{\tilde \Psi})_0]\ ,
$$
 which is similar to the identity (1) in the proof of the theorem. Using
 Theorem 4.7 to estimate
the quadratic term on the right, it follows    that the
$L^2_l$-norm of $\tilde A|_D$ can be estimated in terms of the 
$L^2_{l-1}$-norm of the
restriction of $\tilde A$ to a slightly larger disk
$D_l \subsetint B$ and $\nr(\tilde A,\tilde \Psi)\nr_{L^4}^2$.  
Inductively
we get an estimate of the $L^2_l$-norm of $\tilde A|_D$ in terms 
of the $L^2_1$-norm of $\tilde
A$ and of $\nr(\tilde A,\tilde \Psi)\nr_{L^4}^2$. But both terms
 can be estimated now in terms
of  $\nr F_{\tilde A}\nr_{L^2}$.
\qed

Note that the estimate in terms of $\nr F_{\tilde A}
\nr_{L^2}^{\frac{1}{2}}$ which we
obtained by applying directly Theorem 4.7, is in fact fully 
sufficient for our purposes.
However it is interesting to notice that the Sobolev norms of
 the  connection component
$\tilde A$ can be estimated as in the instanton case in terms of 
$\nr F_{\tilde
A}\nr_{L^2}$.

\begin{co} (Local compactness) There exists   a constant
$\varepsilon>0$ such that the following holds:

For any pair system of data $\pg$ which is sufficiently close to a system 
of standard data
$\pg_{K_0}$ on the ball with $|K_0-\id|<\eg$ , and any sequence
$(A_n,\Psi_n)$ of solutions of the
$PU(2)$-monopole equations    for $\pg$   with $\nr
F_{A_n}\nr_{L^2}\leq
\varepsilon$, there is a subsequence
$m_n$ of
$\N$ and gauge equivalent solutions $(\tilde A_{m_n},
\tilde \Psi_{m_n})$ 
converging in the
${\cal C}^{\infty}$-topology on the open ball $B$.
\end{co}
\qed

We can prove now  the following result, which is the analogon of  
 Proposition 4.4.9 p. 161 [DK].
\begin{co} (Global compactness)  Let $\Omega$ be a  4-manifold
   and let  $P^u$  be a $Spin^{U(2)}(4)$-bundle on $\Omega$ such 
that $\Lambda^1_\Omega\simeq
P^u\times_\pi\R^4$ as oriented 4-bundles. Let  
$\pg=(\gamma,C,a,\beta,K)$ be an arbitrary
system of data for $(\Omega,P^u)$ satisfying the condition 
$|K(x)-\id_{\ad_{+}}|<\eg$ in every
point $x\in\Omega$.

Let $(A_n,\Psi_n)$ be a sequence of solutions of the 
$PU(2)$-monopole equations associated
with $\pg$ such that every point $x\in
\Omega$ has a geodesic ball neighbourhood
$D_x$ such that for all large enough $n$,
$$\int_{D_x}| F_{A_n}|_{g_\gamma}^2 
vol_{g_\gamma}< \varepsilon^2$$
where $\varepsilon$ is the constant in Corollary 4.9. 
Then there is a subsequence
$(m_n)\subset\N$ and gauge transformations $u_{n}\in{\cal G}_0$ 
such that
$ u_{n} (A_{m_n},\Psi_{m_n})$ converges in the 
${\cal C}^{\infty}$-topology on
$\Omega$.
\end{co}
\pf First of all note that every point has a geodesic ball 
neighbourhood $D'_x\subset D_x$ such
that   for a suitable subsequence $(m^x_n)_n\subset\N$ and 
suitable gauge transformations
$u^x_{n}$ over $D'_x$ the sequence $(u^x_n(A_{m^x_n}|_{D'_x},
\Psi_{m^x_n}|_{D'_x})_n$ converges
in the ${\cal C}^{\infty}$ topology on $D'_x$. This follows from 
 Remark 4.5, Corollary  4.9
and the conformal invariance of the $L^2$-norm of 2-forms.

Using now    Corollary 4.4.8 p. 160  [DK]  we get a subsequence $(m_n)_n$
 of $\N$  and
gauge transformations $u_n$ such that $u_n(A_{m_n})$ converges in the ${\cal
C}^{\infty}$-topology on $\Omega$ to a connection $A$. But using the first
 monopole
equation we see that the convergence of the connection component together 
with the local
$L^4$-bound of the spinor component (provided by the local $L^2$-boundedness 
of the
curvature) implies the local boundedness of the spinor component in any 
$L^2_l$-norm.
 \qed
\begin{pr} (Apriori ${\cal C}^0$-boundedness of the spinor) Let $X$ be a 
compact oriented
4-manifold, $ P^u$  a $Spin^{U(2)}(4)$-bundle  on $X$ with 
$P^u\times_\pi\R^4\simeq
\Lambda^1_X$ as oriented 4-bundles, and
$\pg=(\gamma,C,a,\beta,K)$ a system of data for the pair $(X,P^u)$ 
 satisfying the condition
$|K(x)-\id_{\ad_{+}}|<\eg$ in every point $x\in X$. \\  
1. If $\beta=0$, and $C$ is induced via $\gamma$ by the Levi-Civita 
connection in
$(\Lambda^1,g_\gamma)$, then for any solution
$(A,\Psi)\in{\cal A}(\bar\delta(P^u))\times A^0(\Sigma^+(P^u))$
 of the $PU(2)$-monopole
equations associated with $\pg$, the following apriori estimate  holds:
$$\sup\limits_X|\Psi|^2_{g_\gamma}\leq \max\left(0,
C^{-1}\sup\limits_X(-\frac{s}{4}+c|F_a^+|_{g_\gamma})\right)
$$ 
Here $s$ stands for  the scalar curvature of $g_\gamma$, $c$ is a 
universal
positive constant, and 
$C$ is the  universal positive constant in (4) . \\ \\
2. In the general case one has an apriori estimate of the form 
$$\sup\limits_X|\Psi|^2_{g_\gamma}\leq \max\left(0,
C^{-1}\left[\sup\limits_X(-\frac{s}{4}+c|F_a^+|_{g_\gamma})+
\sigma(\gamma,C,\beta)\right]
\right)
\ ,
$$
where   $\sigma(C,\beta,\gamma)$ depends continuously on
the coefficients of $\gamma,C,\beta$ with respect to the ${\cal C}^2\times{\cal
C}^1\times{\cal C}^1$-topology.
\end{pr}
\pf  We prove  the second assertion. Using Remark 3.4, it follows
 that, modifying $\beta$ if
necessary, we may assume that $C$ is induced via $\gamma$ by the
 Levi-Civita connection in 
$(\Lambda^1,g_\gamma)$, so that the Dirac operator $\Dr^C_{\gamma,a,A}$ 
associated
with $C$ coincides with the standard Dirac operator
$\Dr_{\gamma,a,A}$.   

The Weitzenb\"ock formula for coupled Dirac operators gives for any 
triple
$(A,a,\Psi)\in {\cal A}(\bar\delta(P^u))\times
 {\cal A}(\det(P^u))\times
 A^0(\Sigma^+(P^u))$ 
$$\Dr_{\gamma,A,a}^2\Psi=\nabla_{A,a}^*\nabla_{A,a}\Psi+
\Gamma_\gamma[(F_A+\frac{1}{2}F_a)^{+_{g_\gamma}}]\Psi+
\frac{s}{4}\ \Psi \cdot
$$
 On the other hand
$$\Dr_{\gamma,a,A}(\Dr_{\gamma,a,A}+\beta)=
\Dr_{\gamma,a,A}^2+
\gamma\cdot\nabla_{a,A}\circ\beta$$

If $(A,\Psi)$ solves the $PU(2)$-monopole equations for the 
system of data $\pg$, it most hold
pointwise
$$\begin{array}{c}(\nabla_{A,a}^*\nabla_{A,a}\Psi,\Psi)+
 (K (\Psi\bar\Psi)_0(\Psi),\Psi)+
\frac{1}{2}(\Gamma_\gamma(F_a)(\Psi),\Psi)+\\ \\
 +\frac{s}{4}|\Psi|^2+
(\gamma\cdot \nabla_{a,A}\circ\beta(\Psi),\Psi)=0 \ .
\end{array}$$
Using the inequality  (4), we get
$$
\begin{array}{l}
\frac{1}{2}\Delta|\Psi|^2= (\Delta_{A,a}\Psi,\Psi)-
|\nabla_{A,a}\Psi|^2\leq  \\ \\ 
\ \ \ \ \  \leq -C|\Psi|^4+ (c
|F_a^+|-\frac{s}{4})|\Psi|^2
+|(\gamma\cdot\nabla_{a,A}\circ\beta(\Psi),\Psi)|-
|\nabla_{A,a}\Psi|^2\ .
\end{array}\eqno{(7)}$$
On the other hand 
$$\gamma \cdot \nabla_{a,A}\circ\beta(\Psi)=\gamma\cdot
[(\nabla_{C}\beta)(\Psi)+\beta\nabla_{A,a}\Psi]\ .$$
Therefore the term $(\gamma\cdot\nabla_{a,A}\circ\beta(\Psi),\Psi)$ 
can be estimated   as
follows
$$|(\gamma \nabla_{a,A}\circ\beta(\Psi),\Psi)|
\leq c'\left(|\nabla_{C}(\beta)| |\Psi|^2 +
|\beta||\nabla_{A,a}\Psi| |\Psi|\right)\leq $$ 
$$\leq c'\left[|\nabla_{C}(\beta)| |\Psi|^2 +
|\beta| \left(\varepsilon\ |\nabla_{A,a}\Psi|^2 +
\frac{1}{\varepsilon}\  |\Psi|^2\right)\right] \ ,
\eqno{(8)}$$
where $c'$ is a universal constant and $\varepsilon$
 is any positive number.
Choose now
$\varepsilon:=\frac{1}{2(c'\sup|\beta|+1)}$, so that
 the total coefficient of 
$|\nabla_{A,a}\Psi|^2$ in the expression obtained by 
replacing
$|(\gamma\cdot \nabla_{a,A}\circ\beta(\Psi),\Psi)|$ in    (7)
 with the right hand term of (8)
becomes negative. Then we get an inequality of the form
$$\frac{1}{2}\Delta|\Psi|^2\leq -C|\Psi|^4+ \sup\left(c|
F_a^+|+c'|\nabla_{C}(\beta)|  -\frac{s}{4}\right)|\Psi|^2
+\frac{c'\sup|\beta|}{\varepsilon}|\Psi|^2\ ,
$$
and the assertion follows easily by the maximum principle.

\qed
\begin{co} If $\Omega$ is compact, the condition "$\int\limits_{D_x}|
F_{A_n}|_{g_\gamma}^2 vol_{g_\gamma}<\varepsilon^2$ 
for all sufficiently large $n$"  in Corollary 
4.10 can be replaced by the condition
 $$"\int_{D_x}| F_{A_n}^{-_{g_\gamma}}|_{g_\gamma}^2 
vol_{g_\gamma}<\frac{\varepsilon^2}{2}
\   for\ all\ sufficiently\ large\ n\ ".$$
\end{co}
\pf By Proposition 4.11 and the inequality (3),  the pointwise norm
$|F_{A_n}^{+_{g_\gamma}}|$ of the $g_\gamma$-self-dual component  
of the curvature is  apriori
bounded by a constant (depending on $s_{g_\gamma}$ and  $\pg$) hence 
$\int_{D_x}| F_{A_n}^{+_{g_\gamma}}|_{g_\gamma}^2$ can be made 
arbitrarily small, by
replacing eventually
$D_x$ with a smaller ball.
\qed

\subsection{Regularity} 

We begin with the following simple

\begin{re} Let $X$ be a   4-manifold and $g$, $g'$ two metrics on $X$.  
Then the operator
$d_{g}^*+d^{+_{g'}}:A^1\map  A^0 \oplus A^2_{+_{g'}}$ is elliptic. If $X$ is compact then the
kernel of this operator is the harmonic space $\H^1_{g}$. The image of its extension
$L^2_{k+1}\map  L^2_{k}$ is $(A^0)_{k}^{\bot}\oplus
 (A^2_{+_{g'}})^{\bot}_{k}$, where
$(A^0)_{k}^{\bot}$ is the
$L^2_{g}$-orthogonal complement of  $\R\subset (A^0)_{k}$, and
$(A^2_{+_{g'}})^{\bot}_{k}$ is the $L^2_{g'}$-orthogonal 
complement of $\H^2_{+_g'}\subset
(A^2_{+_{g'}})_{k}$.    
\end{re}
Indeed, one checks easily that the symbol $\sigma$ of $d_{g}^*+d^{+_{g'}}$ 
is injective for
non-vanishing cotangent vectors $\xi$. Indeed, if 
$\sigma_\xi(\alpha)=0$, then
$(\xi\wedge\alpha)_{+_{g'}}=0$, hence $ \xi\wedge\alpha=0$. Therefore 
$\alpha$ has the
form $\alpha=c\ \xi$, $c\in\R$. Using now the $\Lambda^0$-component 
of the equation
$\sigma_\xi(\alpha)=0$, we get $c\ |\xi|_{g}^2=0$, i. e. $c=0$. 
But $\Lambda^1$,
$\Lambda^0\oplus \Lambda^2_{+_{g'}}$ have both rang 4, so 
$\sigma_\xi$ must be
isomorphism. 

On compact 4-manifolds one has
$\ker d^{+_{g'}}=\ker d$. Therefore 
$\ker(d_{g}^*+d^{+_{g'}})=\ker(d_{g}^*+d)=\H^1_{g}(X)$.  The
image of the 
$L^2_{k+1}\map  L^2_{k}$ extension of $d_{g}^*+d^{+_{g'}}$
 is obviously contained in
$(A^0)_{k}^{\bot}\oplus (A^2_{+,{g'}})^{\bot}_{k}$.
 Therefore it must coincide with  this space,
because 
${\rm index}(d_{g}^*+d^{+_{g'}})=
{\rm index}(d_{g}^*+d^{+_{g}})=b_1-b_+-1$.

\qed

As in the section above we fix $SU(2)$-bundles $S^{\pm}_s$ on the 4-sphere
 $S$ such that
$\Lambda^1_S\simeq\RSU(S^+_s,S^-_s)=\Hom_\H(S^+_s,S^-_s)$ as 
oriented 4-bundles. The
pairs consisting  of a metric on the sphere and a
$Spin(4)$-structure for that metric are parameterized by linear
 isomorphic Clifford maps 
$$\gamma:\Lambda^1_S\map \Hom_\H(S^+_s,S^-_s)\ .$$

We denote by $Clif(S)$ the space of Clifford maps on the sphere. 
Let again $E_s$ be the
trivial $U(2)$ bundle on $S$.  

We fix a Clifford map $\gamma_s:\Lambda^1_S\map \Hom_\H(S^+_s,S^-_s)$
  such that
$g_s:=g_{\gamma_s}$ has non-negative scalar curvature, 
strictly positive in the south pole
$\infty$.  Therefore the associated selfadjoint Dirac operator 
$\Dr_{\gamma_s}$ is
injective, by the Weitzenb\"ock formula. Denote by $C_s$ the Levi-Civita 
connection
induced by $\gamma_s$ in the $SO(4)$-bundle 
$P^u_s\times_\pi\R^4=\Hom_\H(S^+_s,S^-_s)$ and
denote  by
$\qg_s$ the system of data 
$$\qg_s:=(\gamma_s,C_s,0,0)\in Clif(S)\times{\cal
A}(P^u_s\times_\pi\R^4) \times{\cal A}(\det(P^u_s))\times
 A^0(P^u_s\times_\pi\R^4)\ ,
$$
where we used as usually the identification 
${\cal A}(\det(P^u_0))=A^1(u(1))$.

Denote by 
$$sw_\pg:{\cal A}(\bar\delta(P^ u_s))\times
 A^0(S^+_s\otimes E_s)\map
A^0(S^-_s\otimes E_s)\times A^0(su(S^+_s)\otimes su(2))$$
 the Seiberg-Witten map
associated with a system of data $\pg$ for the pair $(S,P^u_s)$.

\begin{pr} (Regularity of $L^4$-small
$L^2_1$-almost solutions with connection component in
 Coulomb gauge) 
Let $g$ be an arbitrary fixed metric on the sphere. There are positive constants  
 $\alpha$,
$\mu$, $c$ (depending on $g$ and $\gamma_s$) such that for any system of data 
$\pg=(\qg,K)$
with
$\qg$ sufficiently close to $\qg_s$ and $|K-\id_{su(S^+_s)}|<\eg$ the following
 holds.

Any  pair $u=(A,\Psi)\in L^2_1(\Lambda^1 (su(2)))\times 
L^2_1(S^+_S\otimes E_s )$ 
satisfying:\\
(i) $d^*_g(A)=0$,\\
(ii) $\nr u\nr_{L^4}< \alpha$,\\
satisfies the inequality
$$\nr u\nr_{L^2_1}\leq c \ \nr sw_{\pg}(u)\nr_{L^2}\ .$$
If,  moreover\\
(iii) $\nr sw_{\pg}(u)\nr_{L^2}<\mu$,\\
(iv)   $sw_{\pg}(u)$ is smooth,\\
then $u$ is also smooth. 

\end{pr}

\pf We use the method of continuity as in the proof of 4.4.13 [DK].  
 The essential fact used
in the proof of that theorem is that the map 
$$B\longmapsto (d^* B, F_B^+)
$$
can be written as the sum of an \ub{injective} elliptic first order operator and a
quadratic map. By   Remark 4.13, the map $(d_g^*,sw_\pg)$ has the same
 property. Note
that we do \ub{not} require  the metric $g$ to be close to $g_{\gamma_s}$.

As in the proof of Theorem 4.7,  the system $\qg=(\gamma,C,a,\beta)$ defines 
an elliptic
first order operator  on the sphere
$$D_{\qg}:=
\begin{array}{ccccc}
\Dr_{\qg}&:&A^0(S^+_s\otimes E_s)&\map&
A^0(S^-_s\otimes E_s)\\
\oplus&& \oplus&&\oplus\\
d_g^*+\Gamma_\gamma\circ d &:&A^1(su(2))&\map& A^0(su(2))^\bot
\oplus A^0(su(S^+_s)\otimes su(E_s)) 
\end{array}$$

Here   $\Dr_{\qg}$ stands for the Dirac operator $\Dr_\gamma^C+\beta+
\gamma(\frac{a}{2})$,
and $A^0(su(2))^\bot$ for the $L^2_{g}$-orthogonal complement of the
 3-dimensional
space of constant $su(2)$-valued functions.

By Remark 4.13 and elliptic semicontinuity, it follows   that $D_{\qg}$
is injective if $\qg$ is sufficiently
${\cal C}^0$-close to $\qg_s$. Moreover, the $L^2_{k+1}\map L^2_{k}$
 extension of
$D_{\qg}$ is an isomorphisms depending continuously on $\qg$ with respect to
 the ${\cal
C}^k$-topology.   

We extend the operator $d_g^*$ on pairs by putting $d_g^*(B,\Phi):=d_g^*(B)$. 
With this
convention  note that the map $d_g^*+sw_\pg$ can be written as
$$(d_g^*+sw_\pg)(B,\Phi)=D_\qg(B,\Phi)+\left[ \matrix{\gamma(B)\Phi\cr
\Gamma_\gamma(B\wedge
B)-K(\Phi\bar\Phi)_0}\right]=D_\qg(B,\Phi)+B_{\gamma,K}(B,\Phi)\ ,
$$
where $B_{\gamma,K}$ is the quadratic map defined by the square bracket.\\

\ub{Claim 1:} If $\alpha$ is sufficiently small, there exists a constant 
$c=c(g,\gamma_s)$
such that for any $L^2_1$-pair $v$ with $d_g^*v=0$, $\nr v\nr_{L^4}<\alpha$,  
one has the
estimate
$$\nr v\nr_{L^2_1}< c \nr sw(v)\nr_{L^2}\ . \eqno{(1)}
$$

Indeed, the Coulomb gauge condition $d_g^*(v)=0$ implies
$$D_\qg(v)=-B_{\gamma,K}(v)+sw(v)\ . \eqno{(2)}
$$
This gives an estimate of the form
$$ \nr v\nr_{L^2_1}\leq C_\qg\nr D_\qg(v)\nr_{L^2}\leq C_\qg C_{\gamma,K}\nr
v\nr_{L^4}^2+\nr sw(v)\nr_{L^2}\leq $$
$$\leq C C_\qg C_{\gamma,K}\nr v\nr_{L^4}\nr
v\nr_{L^2_1}+\nr sw(v)\nr_{L^2}\ ,
$$
Since   $\qg$ is assumed to be close to $\qg_s$
and $K$ belongs to a bounded family, it follows  that the constants $C_\qg$, 
$C_{\gamma,K}$ can
be chosen independently of $\pg$. The claim follows by the  same rearrangement
 argument used
in the proof of Theorem 4.7, taking
$\alpha\leq \frac{1}{2C C_\qg C_{\gamma,K}}$. This proves the claim and the 
first part of the
theorem.\\

\ub{Claim 2:} If $\alpha$ is sufficiently small, then for any two $L^2_1$-sections 
$v_1$, $v_2$
with 
$d_g^*(v_i)=0$, $\nr v_1\nr_{L^4}<\alpha$, $\nr v_2\nr_{L^4}<\alpha$ and 
 $sw(v_1)=sw(v_2)$ it
follows
$v_1=v_2$.\\

Indeed, let $b_{\gamma,K}$ be the $\R$-bilinear map associated with
 $B_{\gamma,K}$. One
has
$$D_\qg(v_1-v_2)=b_{\gamma,K}((v_2-v_1),v_1)+
b_{\gamma,K}(v_1,(v_2-v_1))\ ,\
$$
hence, by the injectivity  of $D_\qg$, we get an estimate of the form
$$\nr v_1-v_2\nr_{L^4}\leq C \nr v_1- v_2\nr_{L^2_1}\leq  C C_\qg
\nr b_{\gamma,K}((v_2-v_1),v_1)+
b_{\gamma,K}(v_1,(v_2-v_1))\nr_{L^2}
$$
$$\leq C C_\pg (\nr v_1\nr_{L^4}+\nr v_2\nr_{L^4})\nr v_2-
v_1\nr_{L^4}
$$
where $C_\pg$ is a constant  depending continuously of $\pg$  with 
respect to the ${\cal
C}^0$-topology. Therefore, we may   suppose as above that $C_\pg=C_1$ 
is independent of
$\pg$. Take
$\alpha\leq\frac{1}{4C C_1}$.\\

\ub{Claim 3:} \ If $\alpha$ is sufficiently small, then for any \ \ub{smooth} \  pair \ $v$ \ with \
$d_g^*(v)=0$, $\nr v\nr_{L^4}<\alpha$ one has estimates of the form
$$\nr v \nr_{L^2_{k+1}}\leq C_{\pg,k}  \nr sw(v) \nr_{L^2_k} +
 P_{\pg,k}(\nr sw(v)
\nr_{L^2_{k-1}})\ ,
$$
where $C_{\pg,k}$ is a positive constant and $P_{\pg,k}$ is a polynomial with
positive coefficients and without constant term.\\

To see this use again  the rearrangement argument above to estimate 
the $L^2_2$ and the
$L^2_3$ norms of $v$ (compare with the proof of Theorem 3.7). For the 
higher Sobolev norms
apply the usual bootstrapping procedure to the elliptic equation $(2)$. \\

\ub{Claim 4:}   If $\alpha$ is sufficiently small, there exists a positive 
number $\mu$ such that
for every smooth section \ $f\in A^0(S^-_s\otimes E_s \oplus su(S^+_s)
\otimes su(2))$ \ with \
$\nr f\nr_{L^2}<\mu$,   the equation 
$$ sw(v)=f\ ,\ d_g^*(v)=0$$
has a smooth solution $v$ satisfying $\nr v\nr_{L^4}<\alpha$.\\

Indeed,  choose first $\alpha$   such that the conclusions of Claims 1-3 hold.  
We use the
continuity method to find a smooth solution of the equations 
$sw(v)=f$, $d_g^*(v)=0$ . Let
$(SW^t)$ be the equation
$$(d_g^*+sw_\pg)(v)=t\ f \ .\eqno{(SW^t)}
$$
We have to find a smooth solution of $(SW^1)$ whose $L^4$-norm
 is bounded by $\alpha$. Let
$N$ be the set  
$$N:=\{t\in [0,1]| \ (SW^t)\ {\rm has\ a\ smooth\ solution}\ 
v\ {\rm with}\ \nr
v\nr_{L^4}<\alpha\}$$ 

 The set $N$ contains $0$. We assert that, taking a smaller
 bound $\alpha$ if necessary, $N$
becomes an open set. We use the implicit function theorem.
Let $v_0$ be a solution of $(SW^{t_0})$ satisfying $d_g^*(v_0)=0$, $\nr
v_0\nr_{L^4}<\alpha$. We have
$$\frac{\partial}{\partial v}(d_g^*+sw_\pg)(\dot v)=D_\qg(\dot v)+ 
b_{\gamma,K}(\dot v,v)+b_{\gamma,K}(v,\dot v)
$$
This shows that, for $v=0$, the operator $ \frac{\partial}{\partial v}|_{_0}
(d_g^*+sw_\pg)$
defines an isomorphism:
$$
\begin{array}{ccc}
L^2_1(S^+_s\otimes E_s)&\map&
L^2(S^-_s\otimes E_s)\\
 \oplus&&\oplus\\
L^2_1(\Lambda^1(su(2)))&\map& L^2(su(2))^\bot
\oplus L^2(su(S^+_s)\otimes su(E_s)) 
\end{array}$$

If $\nr v\nr _{L^4}$ is sufficiently  small, then the $L^2_1 \map L^2$ 
extension of 
$\frac{\partial}{\partial v}(d_g^*+sw_\pg) $ is still an isomorphism.
 By the Fredholm
alternative it follows that the $L^2_3 \map L^2_2$ extension is  an isomorphism, too.
Therefore, there exists $\varepsilon>0$ and an $L^2_3$ solution $v_t$ of 
$(SW^t)$ for any
$t\in (t_0-\varepsilon,t_0+\varepsilon)$  such that $v_{t_0}=v_0$. 
Using the usual bootstrapping
procedure, it follows that $v_t$ must be smooth. 

We claim that $N$   is closed, if the bound $\mu$ of $\nr f\nr_{L^2}$ is
sufficiently  small. Indeed, if $t_n\rightarrow t_0$, and if $v_n$ is a smooth
 solution of
$(SW^{t_n})$ with $\nr v_{n}\nr_{L^4}<\alpha$, then  Claim 3. shows that
there is a subsequence $(v_{n_m})_{m\in \N}$ converging in the ${\cal
C}^{\infty}$-topology to a smooth section $v_0$, which must solve the 
equation 
$(SW^{t_0})$. Of course, it is not clear that the strict inequality 
$\nr v_{n_m}\nr_{L^4}<\alpha$
is preserved at the limit.  On the other hand,  using the estimate (1) proved
 in Claim 1. and the
boundedness of the inclusion $L^2_1\subset L^4$, we see that,  choosing
$\mu$ sufficiently small, we can assure that
$$\nr v_n\nr_{L^4}\leq \frac{\alpha}{2}\ .
$$
Therefore   $v_0$ satisfies the  stronger inequality 
$\nr v_0\nr_{L^4}\leq\frac{\alpha}{2}$.
Now the second assertion in the theorem follows immediately: If $\nr u\nr_{L^4}<
\alpha$,
$d_g^*(u)=0$, $\nr sw_{\pg}(u)\nr_{L^2}<\mu$, and  $sw_{\pg}(u)$ is
 smooth, we can find a
smooth solution $v$ of the equations $d_g^* v=0$, $sw_{\pg}(v)=sw_{\pg}(u)$
 with $\nr
v\nr_{L^4}<\alpha$. But, by Claim 2., this solution must coincide with $u$.
\qed
\begin{co} With the notations and assumptions of the theorem, the following holds: 
There
exists a positive constant $\alpha_1$ (depending on $(g,\gamma_s)$) such that any
$L^2_1$-pair
$u=(A,\Psi)$ with $d_g^*(A)=0$, $\nr u\nr_{L^2_1}\leq \alpha_1$ and 
$sw_\pg(u)$  smooth, is
also smooth.
\end{co}

\subsection{Removable singularities}

We notice first that Corollary 4.8  (Estimates in terms of the curvature) can be 
easily
generalized to an  {\it arbitrary} system  of data $\pg'=(\qg',K')$   for the pair
$(\bar B, P^u_0)$, not necessarily   close to a standard system.
The only difference is that the constant $\varepsilon$ in the conclusion of the
 theorem will
depend
 on $\pg'$.   To see this it is enough to notice that the operator $D_\qg$ 
constructed in the proof of Theorem 4.7   is always elliptic by Remark 4.13  
(even if the metric
$g_\gamma$ is not close to the   metric $g_s$).  We can use in fact the standard 
 constant
curvature metric on the sphere for the Coulomb condition, as in [DK]. $D_\qg$ 
will be in
general non-injective,   but the injectivity of this operator is not essential in the
 proof
of 4.7: the corresponding elliptic estimates $(el)$, $(el)_k$ will contain on the
 right the
additional term
$\nr u\nr_{L^2}$, which can be estimated in terms of  $\nr u\nr_{L^4}$ using 
the volume of the
sphere endowed with the metric $g_\gamma$.

An alternative argument uses a division of the unit ball in small balls, the scale
 invariance
of the equations (Remark 4.5), the original Theorem 4.7, and the patching arguments
explained on p. 162 [DK] in the instanton case.

 Using this generalization of Corollary 4.8,  we get the following analogon of  
Proposition
4.4.10 [DK]:

\begin{lm} Let $\Omega$ be a strongly simply connected    4-manifold
endowed with a $Spin^{U(2)}(4)$-bundle $P^u$ with $P^u\times_\pi\R^4\simeq
\Lambda^1_\Omega$,  $\bar\delta(P^u)\simeq \Omega\times PU(2)$. Fix a 
trivialization of the
$PU(2)$-bundle $\bar\delta(P^u)$. 

Let $\pg=(\gamma,C,a,\beta,K)$ be a system of data for the bundle $P^u$ such that
pointwise $|K-\id|<\eg$.

There exists a positive constant
$\varepsilon_\pg$, and for every precompact interior domain 
$\Omega'\subsetint\Omega$ there exists a positive constant $M_{\pg,\Omega'}$ such
that any solution $(A,\Psi)$
 of the $PU(2)$-monopole equations for $\pg$    with $\nr
F_A\nr_{L^2_{g_\gamma}}<\varepsilon_\pg$ is gauge equivalent  over $\Omega'$ 
to a pair
$(A',\Psi')$ satisfying 
$$\nr A'\nr_{L^4_{g_\gamma}(\Omega')}<
M_{\pg,\Omega'}\nr F_A\nr_{L^2_{g_\gamma}}\ .$$
\end{lm}
\qed
\begin{re} Given a fixed   system of data $\pg_0$,   we can find constants 
$\varepsilon_0$,
$M_{0,\Omega'}$ (independent of $\pg$) such that the conclusion of the
 theorem holds with
these constants,  for every
$\pg$ sufficiently close to $\pg_0$. Moreover, the statement is true if we  
use the fixed
metric $g_{\gamma_0}$ to compute the Sobolev norms.
\end{re}

We will need these results in the following particular case:

  Let ${\cal N}$,
${\cal N}'$ be the annuli
$${\cal N}:=\{x\in B|\ \frac{1}{2}<|x|<1\}\ ,\ \ {\cal N}':=\{x\in B|\
\frac{4}{6}<|x|<\frac{5}{6}\}  \ . $$

Denote by ${\cal N}_r$, ${\cal N}'_r$ the images of ${\cal N}$, ${\cal N}'$ 
under the
homothety $h_r$. We recall that we denoted by $P^u_0$ the trivial
$Spin^{U(2)}(4)$-bundle on $\bar B$, which is associated with the triple of
$SU(2)$-bundles $S^{\pm}_0:=\bar B\times\H_{\pm}$, $E_0:=\bar
 B\times\C^2$.
\begin{lm} Let $\pg=(\gamma,C,a,\beta,K)$ be a system of data for the trivial
$Spin^{U(2)}(4)$-bundle $P^u_0$ on the ball $\bar B$, such that pointwise
$|K-\id|<\eg$, and   such that
$\gamma|_{\Lambda^1_0}:\R^4\map\H=(P^u_0\times_\pi\R^4)_0$ is
 the standard
identification..  Then there exists
constants $\varepsilon(K_0)>0$,   $M(K_0)$  such that for any sufficiently small
$r>0$ the following holds:\\
 
Any solution $(A,\Psi)$
 of the $PU(2)$-monopole equations for \ $\pg|_{{\cal N}_r}$  \  with \ 
$\nr F_A\nr_{L^2({\cal N}_r)}<\varepsilon(K_0)$ is gauge equivalent  over
 ${\cal N}'_r$ to a pair
$(A',\Psi')$ satisfying 
$$\nr A'\nr_{{L^4}({\cal N}_r')}<
M(K_0)\nr F_A\nr_{{L^2}({\cal N}_r)}\ .$$
\end{lm}

The  constants $\varepsilon(K_0)>0$,   $M(K_0)$ are
independent of $r$, and the Sobolev norms are computed with respect to  the standard
euclidean metric.

\pf We use the same argument as in Remark 4.5. Let 
$$h_r:({\cal N}',{\cal N})\map ({\cal N}'_r,{\cal N}_r)$$
 the homothety of slope $r$.

The pair $(h_r^*(A),rh_r^*(\Psi))$ solves the monopole equations 
associated with the
system of data $(\gamma_r:=rh_r^*(\gamma),h_r^*(C),h_r^*(a),
rh_r^*(\beta),h_r^*(K))$, which
converges to the standard system $\pg_{K_0}$ restricted to ${\cal N}$, as
$r\rightarrow 0$.  

The result follows now  from 4.16, 4.17 and the conformal invariance of the
$L^4$-norm  on 1-forms  and of the $L^2$-norm on 2-forms.
\qed

We shall use the following notations 
$$\Omega_r:=B\setminus \bar B(r)\ ,\ \ B^{\bullet}=B\setminus\{0\}\ , \ \
B^{\bullet}(R)=B(R)\setminus\{0\}\  \ , S^{\bullet}=S\setminus\{0\}\ .
$$

\begin{lm} Let $\pg=(\gamma,C,a,\beta,K)=(\qg,K)$ be a system of data for the
trivial bundle $P^u_0$ on the ball $B$, and let $(A,\Psi)$ be a pair on
$B^{\bullet}$ solving the monopole equations for $\pg|_{B^{\bullet}}$ such that
$$\int\limits_{B^\bullet}|F_A|^2<\infty
$$
Then for any sufficiently small $r>0$, there exist an $SU(2)$-bundle $E_r$ over $B$,
a pair $(A_r,\Psi_r)\in{\cal A}(E_r)\times A^0(S^+_0\times E_r)$ and an
$SU(2)$-isomorphism
$$\rho_r:E_r|_{\Omega(r)}\map E|_{\Omega(r)}
$$
such that:\\
i) $\rho_r^*(A,\Psi)=(A_r,\Psi_r)$,\\
ii) $\nr sw_\pg(A_r,\Psi_r)\nr_{L^2(B)}\rightarrow 0$ as $r\rightarrow 0$.
 \end{lm}

\pf  Let $\varphi$ be a cut-off map $\varphi: B\map  [0,1]$ which is identically 1 on
$B\setminus B(\frac{5}{6}-\varepsilon)$ 
and identically 0 on $B(\frac{4}{6}+\varepsilon)$.

Put $\varphi_r:=\varphi\circ h_r^{-1}$. Note first that, by the conformal
 invariance of
the
$L^4$-norm on 1-forms, the norm
$\nr d\varphi_r\nr_{L^4}$ (computed with the euclidean metric) does not depend on
$r$.

Consider now the restriction of the pair $(A,\Psi)$ to ${\cal N}_r$.  Since the
total integral of $|F_A|^2$ on the ball is finite, it follows that for any
 sufficiently small $r>0$  we have 
$$\nr F_A\nr_{L^2({\cal N}_r)}<\varepsilon(K_0)\ , $$
so that Lemma 4.18 applies. The conclusion of this Lemma 
 can be reformulated as follows: There exists an $SU(2)$-trivialization
${\cal N}'_r\times\C^2\textmap{\tau_r} E_0|_{{\cal N}'_r}$ such that  the connection
matrix of  $\tau_r^*(A)$ (which we also denote by $\tau_r^*(A)$) satisfies the
estimate
$$\nr \tau_r^*(A)\nr_{L^4({\cal N}'_r)}\leq M(K_0)\nr F_A\nr_{L^2({\cal N}_r)}
\eqno{(1)}
$$

We define the $SU(2)$-bundle $E_r$ by gluing (over the annulus 
 ${\cal N}'_r$) the
trivial bundles $B(0,\frac{5r}{6})\times\C^2$, 
$E_0|_{\Omega(\frac{4r}{6})}$  via the
isomorphism $\tau_r$.

Let $P^u_r$ be the $Spin^{U(2)}(4)$-bundle associated with the triple
$(S^{\pm}_0,E_r)$. The system $\pg$ can be also  regarded as a 
system of data for the bundle
$P^u_r$. 

Now denote by $u$ the initial pair $u:=(A,\Psi)$, and by $u_r$ the pair 
 $$u_r\in {\cal A}(E_r)\times A^0(S^+_0\otimes E_r)\ ,$$  
which coincides with $u$ on $\Omega((\frac{5}{6}-\varepsilon)r)$ and 
 with the
cut-off $\varphi_r \tau_r^*(u)$ of $\tau_r^*(u)$  on $B(0,\frac{5r}{6})$.

The section $sw_\pg(u_r)$ vanishes identically on 
${\Omega(\frac{5r}{6})}$, where
$u_r$ coincides with $u$. Therefore, in order to prove $ii)$ we 
only have to estimate the $L^2$
norm of
 $sw_\pg(\varphi_r \tau_r^*(u))$ on $B(0,\frac{5r}{6})$, where $E_r$ 
coincides with
the trivial bundle  $B(0,\frac{5r}{6})\times\C^2$.  

On $B(0,\frac{5r}{6})$ the Seiberg-Witten map $sw_\pg$ can be
written as a sum between a first order differential operator and a quadratic map:
$$sw_{\pg}(B,\Phi)=\left[\matrix{\Dr_\qg \Phi\cr
 \Gamma_\gamma(d B)  }\right]+
\left[\matrix{\gamma(B)(\Phi)\cr \Gamma_\gamma(B\wedge
B)-K(\Phi\bar\Phi)_0}\right]=T_\qg(B,\Phi)+B_{\gamma,K}(B,\Phi)
$$

Since $T_\qg$ is a first order operator, we have an identity of the form
$$T_\qg(f v)=A_\qg(df)(v)+fT_\qg(v)\ ,
$$
where $A_\qg(df)$ is a 0-order operator whose coefficients depend linearly on the
first order derivatives of the real function $f$. 

Therefore 
$$sw_\pg(\varphi_r \tau_r^*(u))=A_\qg(d\varphi_r)(\tau_r^*(u))+\varphi_r
T_\qg(\tau_r^*(u))+\varphi_r^2 B_{\gamma,K}(\tau_r^*(u))=
$$
$$=\varphi_r
sw_\pg(\tau_r^*(u))+A_\qg(d\varphi_r)(\tau_r^*(u))+
(\varphi_r^2-\varphi_r)B_{\gamma,K}(\tau_r^*(u))=
$$
$$A_\qg(d\varphi_r)(\tau_r^*(u))+
(\varphi_r^2-\varphi_r)B_{\gamma,K}(\tau_r^*(u))\ .
$$

Therefore, taking into account that $d\varphi_r$ and $(\varphi_r^2-\varphi_r)$  
vanish outside ${\cal N}'_r$, we get 
$$\nr sw_\pg(u_r)\nr_{L^2(B)}=
\nr sw_\pg(\tau_r^*(u))\nr_{L^2(B(\frac{5r}{6}))}\leq
$$
$$\leq C_\qg\nr d\varphi_r\nr_{L^4} 
\nr\tau_r^*(u)\nr_{L^4({\cal N}'_r)}+
 C'_\pg \nr\tau_r^*(u)\nr_{L^4({\cal N}'_r)}^2
$$
Since $\nr d\varphi_r\nr_{L^4}$ does not depend on $r$ we have 
only to prove that 
$\nr\tau_r^*(u)\nr_{L^4({\cal N}'_r)}$ converges to 0 as $r\rightarrow 0$. 
 But the
estimate $(1)$ shows that the $L^4$-norm of the connection component of
$\tau_r^*(u)$ converges to 0 as  $r\rightarrow 0$.

On the other hand, by the inequality (3) Section 4.1 and the second 
monopole equation,
one has pointwise in ${\cal N}'_r$. 
$$|\tau_r^*(\Psi)|^4=|\Psi|^4\leq
 \left[C_1^{-1}|\Gamma_\gamma(F_A)|\right]^2 \ .
$$
This gives an estimate of $\nr\tau_r^*(\Psi)\nr_{L^4({\cal N}'_r)}$
 in terms of 
$\nr F_A\nr_{L^2({\cal N}'_r)}^{\frac{1}{2}}$, which obviously 
converges to 0 as
$r\rightarrow 0$.
\qed

We recall from [DK] the following important

\begin{thry}  (Gauge fixing on the sphere \footnote{ Note that in [DK] it is stated 
 a slightly
weaker form of this theorem (Proposition 2.3.13 p. 63): The hypothesis 
requires that $A$
can be joined to the flat connection by a path of  connections with $L^2$-small
curvature  . However, the second proof of this result, which is given in section 
2.3.10, 
{\it does not} use this additional assumption. I am grateful  to Peter Kronheimer for
pointing me out this important detail.  On the other hand, note that this second proof
 works  only for the standard constant curvature metric, and can be generalized to 
conformally
flat metrics with non-negative sectional curvature. Since our regularity theorem 
works for
solutions whose connection component is in Coulomb gauge with respect to 
 {\it any}
metric, not necessary close to the metric defined by the $Spin^{U(2)}$-structure, 
we don't need
this generalization}) Let
$g_c$ be the standard constant curvature metric on the sphere $S^4$.  
Then there are constants
$\varepsilon_c$,
$M_c$ such that any connection
$A$ in  the trivial $SU(2)$-bundle $E_s$ with 
$\nr F_A\nr_{L^2}<\varepsilon_c$ is gauge 
equivalent to a connection $\tilde A$ satisfying
$$d_{g_c}^*(\tilde A)=0\ ,\ \ \nr \tilde A\nr_{L^2_1}<M_c\nr F_A\nr_{L^2} \ .
$$
\end{thry}
\qed 

We can prove now 

\begin{thry} (Removable singularities)  Let $\pg=(\qg,K)$ a system of
 data for
the trivial $Spin^{U(2)}(4)$-bundle $P^u_0$ on $\bar B$ and let $u=(A,\Psi)$ 
be a pair on the
punctured ball solving the monopole equations for  $\pg|_{B^{\bullet}}$ such that
$$\int_{B^{\bullet}}|F_{A}|^2<\infty\ .$$ 
There exists  an $SU(2)$-bundle  $F$ on the ball, and an $SU(2)$-isomorphism
$\rho:F|_{B^{\bullet}}\map E_0|_{B^{\bullet}}$ such that $\rho^*(A,\Psi)$ 
extends to a
global smooth solution   of the monopole equations associated with $\pg$ and the
$Spin^{U(2)}(4)$-bundle defined by $(S^{\pm}_0,F)$.
\end{thry}

 \pf   We use similar arguments as in the proof of the "Removable singularities"
 theorem for
the instanton equation (Theorem 4.4.12 [DK]).  The only difference is that the 
$L^2_1$-bound of
the approximate solutions we construct, does not follow directly from 
Theorem 4.20 (Gauge
fixing on the sphere).

Identify $\bar B$ with the upper hemisphere of the 4-sphere $S$, and extend
 the system $\pg$ to a
system for the $Spin^{U(2)}(4)$-bundle $P^u_s$. The extended system will
 be denoted by the same
symbol $\pg$, and we can assume that $\pg$ has the form $(\qg,K)$ with
 $\qg$ close to the
system $\qg_s$ constructed in the proof of Theorem 4.7, so that Theorem 
4.14 and Corollary
4.15 applies. We shall use these results in the particular case $g=g_c$; 
with respect to
this metric  connections with $L^2$-small curvature can be brought in 
the Coulomb gauge, by
4.20. \\

Step 1.  For a sufficiently small positive  number $R<1$ we  use   
Lemma 4.18 to get a
trivialization  of $E_0|_{{\cal N}'_R}$, such that  the $L^4$-norm 
of the corresponding
connection matrix is controlled by $\nr F_A\nr_{L^2({\cal N}_R)}$. 
By the same gluing
procedure we get a bundle $E^R$ on the punctured sphere $S^{\bullet}$,
 trivialized on $S\setminus
\bar B(\frac{4R}{6})$. We cut off the pair $u$ this time towards the outer
 boundary of the ball,
and we get a pair $u^R=(A^R,\Psi^R)$. It holds
$$\lim\limits_{R\rightarrow 0} \nr sw_\pg(u^{R})\nr_{L^2}=
\lim\limits_{R\rightarrow 0}
\nr F_{A^{R}}\nr_{L^2}=\lim\limits_{R\rightarrow 0} 
\nr \Psi^{R}\nr_{L^4}=0 \ , \eqno{(2)}
$$
The first two relations follow as in the proof of Lemma 4.19, since both 
maps $sw_\pg(\cdot)$,
$F_{\cdot}$ can be written as the sum of a first order operator and a 
quadratic map, hence the
perturbations produced by  of the two cut-off operations can be estimated 
in terms of the
$L^2$-norm of the curvature restricted to the corresponding annuli. 

To get the third formula, it is enough to notice that the
pointwise norm of the spinor is invariant under bundle isomorphisms,  and that
 the $L^4$-norm
of $\Psi|_{B^\bullet(R)}$ can be estimated in terms of  $\nr
F_{A|_{B^\bullet(R)}}\nr_{L^2}^\frac{1}{2}$.

Suppose now that  $r<R<1$ and use the same procedure (to modify the bundle
 and cut off
the solution), but this time in both directions. 

We get   
$SU(2)$-bundles, $E_r^R$ on the sphere, which come  with trivializations 
over
$B(\frac{5r}{6})$, $S\setminus \bar B(\frac{4R}{6})$, and with an
 isomorphism
$$E_r^R|_{B(\frac{5R}{6})\setminus \bar B(\frac{4r}{6})}
\textmap{\simeq\ \rho_{r,R}}
E_s|_{B(\frac{5R}{6})\setminus \bar B(\frac{4r}{6})}\ ,$$
as well as cut-off pairs
$$u_r^R=(A_r^R,\Psi_r^R)\in {\cal A}(E_r^R)
\times A^0(S^+_s\otimes E_r^R)\ .$$
With this construction, it holds
$$\lim\limits_{r\rightarrow 0} \nr sw_\pg(u_r^R)\nr_{L^2}=\nr
sw_\pg(u^R)\nr_{L^2}\ ,\ \ \lim\limits_{r\rightarrow 0}
\nr F_{A_r^R}\nr_{L^2}=\nr F_{A^R}\nr_{L^2}\ , $$ 
$$\lim\limits_{r\rightarrow 0} \nr
\Psi_r^R\nr_{L^4}=\nr \Psi ^R\nr_{L^4} \ . \eqno{(3)}
$$

Note that the double gluing-procedure we used could apriori give rise to
 a \ub{non}-\ub{trivial}
$SU(2)$-bundle $E_{r,R}$ on the sphere. But since the curvature
 $F_{A^R_r}$ can be made
as small as we please, it follows that  all the bundles
$E_{r,R}$ become trivial, if $R$ is small.

Step 2. Using (2), (3)  and Theorem 4.20 it follows that, once $R$ is small, 
 there exists
an
$SU(2)$-isomorphism $\theta_r^R:E_s\map E_r^R$ such that
$B^R_r:=\theta_{r,R}^*(A^R_r)$ satisfies
$$d_{g_c}^*(B^R_r)=0\ , \ \ \nr B^R_r\nr_{L^2_1}\leq M_c \nr
 F_{A^R_r}\nr_{L^2} \eqno{(4)}
$$
Put $\Phi^R_r:=(\theta^R_r)^*(\Psi_r^R)$, $v^R_r:=(B^R_r,\Phi^R_r)$.

Step 3. Using   (2),  (3), (4) and the boundedness of the embedding 
$L^2_1\subset L^4$, it
follows that, if
$R$ is small enough, the $L^4$-norm of the pair $v^R_r$
 can be made
smaller as the constant $\alpha$ in the Regularity  Theorem 4.14, so that
 we get an estimate of
the form
$$\nr  v^R_r\nr_{L^2_1}\leq c \nr sw_\pg( v^R_r) \nr_{L^2}=
\nr  sw_\pg(u^R_r)\nr_{L^2}\
.\eqno{(5)}
$$
The relations (2), (3) imply now that , choosing $R$ small, we can 
 assure that 
$$\nr  v^R_r\nr_{L^2_1}\leq \alpha_1\ , \eqno{(6)} 
$$
where $\alpha_1$ is the constant in Corollary 4.15. From this point the 
proof goes further like in
the instanton case: We  choose $R$ sufficiently small such that all the mentioned 
properties
are fulfilled, and we let $r$ tend to 0. Using the $L^2_1$-boundedness obtained in 
  (6) it
follows that we can find a sequence
$r_i\rightarrow 0$ such that $v_i=(B_i,\Phi_i):=v^R_{r_i}$ converges 
weakly in $L^2_1$ to an
$L^2_1$-pair $v=(B,\Phi)$.

Step 4. We want to prove that $v$ is smooth. The weak limit $v$ must
 also satisfies $\nr 
v
\nr_{L^2_1}\leq \alpha_1$ by the weak-semicontinuity  of the norm in 
reflexive Banach spaces.
Therefore, by Corollary 4.15, we only have to prove that the $L^2$-section 
$sw(v)$ is smooth.  

But on any small ball $D$ , $\bar D\subset S^{\bullet}$, the pairs 
$v_i=(B_i,\Phi_i)$ remain in
the same gauge equivalence class. Recall now from [DK] that   the Sobolev 
norms of any
connection $H$ in Coulomb gauge can be estimated in terms of the gauge 
invariant expressions
$$\nr F_H\nr_{L^\infty}\ , \ \nr \nabla_H^{(i)} F_H\nr_{L^2} \ ,
$$
as soon as its $L^4$-norm is sufficiently small. Using the estimate (4) and 
the scale
invariance of the $L^4$-norm on 1-forms,   this condition will be also 
fulfilled (for all
 small balls  $D$), if $R$ is sufficiently  small.

On the other hand one can easily bound the Sobolev norms of a spinor $\Xi$ 
in terms of the gauge
invariant expressions $\nr \nabla^{(i)}_H\Xi\nr_{L^2}$ and the Sobolev
 norms of the connection
$H$.

Therefore, taking a subsequence if necessary, we can assume that $v_i$ 
converges in the Fr\'echet
${\cal C}^{\infty}$-topology on $S^{\bullet}$, so that $sw(v)$ is smooth 
on the punctured sphere.

But, by Lemma 4.19,  $\lim\limits_{i}
 \nr sw_\pg(v_i|_{B(\frac{4R}{6})})\nr _{L^2}\rightarrow
0$, so $sw(v)$, which is the limit of $sw(v_i)$ in the distribution sense, 
vanishes in a
neighbourhood of $0$.

On the other hand, for any ball $D$, $\bar D\subset B^\bullet(\frac{4R}{6})$,
 the isomorphism
$\theta^R_{r_i}$ intertwines the connection matrices $A$, $B_i$, and $B_i$ 
converges in the
${\cal C}^{\infty}$ topology on such a ball. Therefore a subsequence 
$\theta^R_{r_{i_n}}$  
converges in the ${\cal C}^{\infty}$ topology  on $B^\bullet(\frac{5R}{6})$ to a 
smooth bundle
isomorphism $\theta$, such that
$$\theta^*(A|_{B^\bullet(\frac{5R}{6})})=B|_{B^\bullet(\frac{5R}{6})}\ .$$
 Taking the limit of
$[\theta^R_{r_{i_n}}]^*(\Psi|_{{B(\frac{5R}{6})\setminus \bar
B(\frac{4r_{i_n}}{6})})})=\Phi^R_{r_{i_n}}|_{{B(\frac{5R}{6})
\setminus \bar
B(\frac{4r_{i_n}}{6})})}$ for $n\rightarrow\infty$, we also get 
$$\theta^*(\Psi|_{B^\bullet(\frac{5R}{6})})=
\Phi|_{B^\bullet(\frac{5R}{6})})
$$
\qed

\subsection{Compactified moduli spaces}

Let $X$ be a closed oriented   4-manifold.  For a 
$Spin^{U(2)}(4)$-bundle $P^u$  with
$P^u\times_\pi\R^4\simeq \Lambda^1$ and a system of data 
$\pg=(\gamma,C,a,\beta,K)$ for
$P^u$  denote by ${\cal M}_\pg(P^u)$ the moduli space of pairs
$(A,\Psi)\in{\cal A}(\bar\delta(P^u))\times A^0(\Sigma^+(P^G))$  
solving the $PU(2)$-monopole
equations associated with $\pg$.

By Proposition 2.1, the data of a 
$Spin^{U(2)}(4)$-bundle  $P'^u$ on $X$ with $\det(P'^u)\simeq \det(P^u)$, 
$P'^u\times_\pi\R^4\simeq P^u\times_\pi\R^4$ is equivalent via 
the map $\bar\delta$ to the
data of $PU(2)$-bundle
$\bar P'$ whose Pontrjagin class satisfies
$$p_1(\bar P')\equiv (w_2(X)+ \bar c_1(\det(P^u)))^2\  {\rm mod}\ 4 \ .
$$

For every   number $l\in \N$ we fix:\\
1. A $Spin^{U(2)}(4)$-bundle $P_l^u$ with
$$l=\frac{1}{4}\left(p_1(\bar\delta(P_l^u))- 
(p_1(\bar\delta(P^u))\right)$$
2. Identifications 
$$P_l^u\times_\pi\R^4\textmap{\simeq}P^u\times_\pi\R^4\ ,\ \
\det(P_l^u)\textmap{\simeq}\det(P^u). \eqno{(1)}
$$
These bundle isomorphisms allow  us to identify the   spaces  of 
perturbations-data
associated  with the bundles $P^u$, $P^u_l$. 
\begin{dt}   An \ ideal
$PU(2)$-monopole  \ of  \ type  \ $(P^u,\pg)$  \ is a pair 
$([A',\Psi'],\{x_1,\dots,x_l\})$
consisting of an element $\{x_1,\dots,x_l\}$    in   a symmetric power
 $S^l(X)$ of
$X$ and a monopole $[A',\Psi']\in{\cal M}_\pg(P_l^u)$.
\end{dt}

We denote by  $I{\cal M}_\pg(P^u)$ the space of ideal monopoles 
of type $(P^u,\pg)$.\\

 Let $\delta_{x}$ be  the Dirac measure associated with a point $x\in X$. If
$\pg=(\gamma,A,a,\beta,K)$, we   always use the metric $g_\gamma$ to
 compute the norms
and to define (anti-)self-duality for 2-forms.
\begin{lm}  The map $F:I{\cal M}_\pg(P^u)\map [{\cal C}^0(X,\R)]^*$, 
defined  by 
$$F([A',\Psi'],\{x_1,\dots,x_l\})=|F_{A'}|^2+8\pi^2
\sum\limits_{i=1}^l\delta_{x_i} \ ,
$$
is bounded with respect to the strong topology in the dual space
 $ [{\cal C}^0(X,\R)]^*$.
\end{lm}
\pf Let $\varphi\in{\cal C}^0(X,\R)$ with $\sup\limits_X|\varphi|\leq 1$. 
Then 
$$\begin{array}{ll}|\langle F([A',\Psi'],\{x_1,\dots,x_l\}),
\varphi\rangle|\leq& \left[\nr
F_{A'}^-\nr_{L^2}^2-
\nr F_{A'}^+\nr_{L^2}^2\right] +2\nr F_{A'}^+\nr_{L^2}^2 +8\pi^2 l
 \\ \\
&=-2\pi^2 p_1(\bar\delta(P^u))+2 C\nr \Psi'\nr_{L^4}^4\ ,
\end{array}$$
where $C$ is a universal positive constant. The assertion follows from 
the apriori ${\cal
C}^0$-boundedness of the spinor component of a solution (Proposition 4.11).
\qed

Let $\mg'=([A',\Psi'],s')$ be an ideal monopole of type $(P^u,\pg)$ 
with $s'\in S^{l'}(X)$ and
$[A',\Psi']\in{\cal M}_\pg(P^u_{l'})$. For a positive number $\varepsilon$
 we define
$U(\mg',\varepsilon)$ to be the set of ideal monopoles $\mg''=([A'',\Psi''],s'')$
 of type
$(P^u,\pg)$ with $s''\subset s'$,  and which have the  following property:
 
There exists an
isomorphism of $Spin^{U(2)}(4)$-bundles
$$\varphi:P_{l''}^u|_{X\setminus s'}\map P_{l'}^u|_{X\setminus s'}$$
which is  compatible with the
identifications (1) such that
$$d_1(\varphi^*(A',\Psi'),(A'',\Psi''))<\varepsilon\ ,$$
 where $d_1$ is a metric defining the
Fr\'echet ${\cal C}^{\infty}$-topology in  the product 
$${\cal A}(\bar\delta(P_{l''}^u|_{X\setminus
s'}))\times A^0(\Sigma^+(P_{l''}^u|_{X\setminus s'})))\ .$$

Let $M>0$ be a bound  for the map $F$ defined above. The weak topology 
in the ball of radius
$M$ in 
$[{\cal C}^0(X,\R)]^*$ is metrisable (see [La], Theorem 9.4.2). Let $d_2$ be a metric
defining this topology.

We endow    $I{\cal M}_\pg(P^u)$  with a metric topology by taking as
 basis of
open neighbourhoods for  an ideal monopole$\mg'$  of type  $(P^u,\pg)$
  the sets of the form
$U(\mg',\varepsilon)\cap F^{-1}(B_{d_2}(F(\mg'),\varepsilon))$, 
$\varepsilon>0$. 

\begin{thry}
With respect to the metric topology defined above the moduli space  
${\cal M}_\pg(P^u) \subset 
I{\cal M}_\pg(P^u) $  is an open subspace with compact closure
$\overline{{\cal M}_\pg(P^u)}$.
\end{thry}
\pf The first assertion is obvious. For the second, we use the same argument as
 in the
instanton case, but we make use in an essential way of the
 ${\cal C}^0$-boundedness of the
spinor:

Let $\mg_n$ a sequence of ideal monopoles. It is easy to see that we 
can reduce the general
case to the case  where $\mg_n=[A_n,\Psi_n]\in {\cal M}_\pg(P^u)$. 
By Lemma 4.23, the
sequence of measures $\mu_n:=F(\mg_n)$ is bounded, so  after replacing 
$\mg_n$ by a
subsequence, if necessary,  it converges weakly to a (positive) measure
 $\mu$ of total
volume
$\mu(1)\leq M$. The set 
$$S_\varepsilon:=\{x\in X|\exists n\in\N\ \forall m\geq n\ 
 (\mu_m(D)\geq
\varepsilon^2\ {\rm for\ every}\ {\rm geodesic\ ball}\ D\ni x)\} $$
contains at most $\frac{M}{\varepsilon^2}$ points, so it is finite for
 every positive number
$\varepsilon$. Choosing the constant $\varepsilon$ provided by  the  
 "Global
compactness" theorem  (Corollary 4.10), it follows by a standard 
diagonal procedure that
there exists a subsequence $(\mg_{n_m})_m$ and gauge transformations
$f_m$ on $X\setminus S_\varepsilon$, such that  $f_m^*(\mg_{n_m})$ 
converges to a 
solution
$(A_0,\Psi_0)$ of the  monopole equations $SW_\pg$ restricted  
 to $X\setminus
S_\varepsilon$. By the "Removable Singularities" theorem, we can 
extend this solution to a
global solution $(\tilde A_0,\tilde \Psi_0)$ of the monopole equations 
associated with $\pg$
and a new $Spin^{U(2)}(4)$-bundle $P'^u$, which comes with 
  identifications
$$P'^u\times_\pi\R^4\textmap{\simeq}P^u\times_\pi\R^4\ ,\ \
\det(P'^u)\textmap{\simeq}\det(P^u).  
$$
 We have 
$$ |F_{\tilde A_0}|^2 =|F_{A_0}|^2=\mu-8\pi^2
\sum_{x\in S_\varepsilon}\lambda_x\delta_x  
$$
with positive numbers $\lambda_x$. It remains to prove that
 the $\lambda_x$ are integers.
 Since $F(\mg_{n_m})\rightarrow\mu$, we have for small enough $r>0$
$$\lambda_x=\lim\limits_{m\rightarrow\infty}\frac{1}{8\pi^2}
\int\limits_{B(x,r)}|F_{A_{n_m}}|^2-|F_{\tilde A_0}|^2=
$$
$$\lim\limits_{m\rightarrow\infty}\frac{1}{8\pi^2}
\int\limits_{B(x,r)}-\tr (F_{A_{n_m}}^2)+\tr (F_{\tilde A_0})^2+
2 \left(| (\Psi_{n_m}\bar\Psi_{n_m})_0|^2-|
 (\tilde \Psi_0\bar{\tilde \Psi_0})_0|^2\right) \ .
$$

As in the instanton case we get  
$$\lim\limits_{m\rightarrow\infty}\frac{1}{8\pi^2}
\left[\int\limits_{B(x,r)}-\tr (F_{A_{n_m}}^2)+
\tr (F_{\tilde A_0})^2\right]_{{\rm mod}\
\Z}=\lim\limits_{m\rightarrow\infty}
(\tau_{S(x,r)}(A_{n_m})-\tau_{S(x,r)}({\tilde A_0}))$$
$$=0\ {\rm in}\
\qmod{\R}{\Z}
$$
by the convergence 
$f_m^*(A_{n_m}|_{X\setminus S_\varepsilon})\rightarrow {\tilde
A_0}|_{X\setminus S_\varepsilon}$. Here
$\tau_S(B)$ denotes the Chern-Simons invariant of 
the connection $B$ on a 3-manifold $S$
([DK]). 

On the other hand, by  the apriori
${\cal C}^0$-bound of the spinor component on the space of 
monopoles, the term
$\int\limits_{B(x,r)}2
\left(| (\Psi_{n_m}\bar \Psi_{n_m})_0|^2-
| (\tilde \Psi_0\bar{\tilde \Psi_0})_0|^2\right)$ can be
made as small as we please by choosing $r$ sufficiently small.  
This shows that the
$\lambda_x$ are integers, and that  
$$\sum_{x\in S_\varepsilon}\lambda_x=
\frac{1}{4}\left(p_1(\bar\delta(P'^u))-
(p_1(\bar\delta(P^u))\right)\ ,
$$
which completes  the proof.
\qed

\newpage

\centerline{\large{\bf References}}
\vspace{5 mm}
\parindent 0 cm

[Ad] Adams, R.: {\it Sobolev spaces}, Academic Press (1975)

[AHS] Atiyah M., Hitchin N. J., Singer I. M.: {\it Self-duality in
four-dimensional Riemannian geometry}, Proc. R. Lond. A. 362, 425-461 (1978)

[B] Besse, A.: {\it Einstein manifolds}, Springer Verlag (1987)

[BB] Booss, B.; Bleecker, D. D.: {\it Topology and Analysis. 
The Atiyah-Singer index formula
and Gauge-Theoretic Physics}, Springer Verlag, (1985)

[C] Chen, B.: {\it Geometry of submanifolds}, Pure and Applied Mathematics, New
York (1973)
 
[D] Donaldson, S.: {\it Anti-self-dual Yang-Mills connections over 
complex algebraic surfaces and stable vector bundles}, Proc. London Math.
Soc. 3, 1-26 (1985)
 
[DK] Donaldson, S.; Kronheimer, P.B.: {\it The Geometry of 
four-manifolds}, Oxford Science Publications   (1990)

[F] Feehan, P.  {\it   Generic metrics, irreducible rank-one
$PU(2)$-monopoles, and transversality}, preprint,  dg-ga/9809001,
(1998), to appear in Comm. Anal. Geom.

[FL] Feehan, P.; Lennes, Th. {\it $PU(2)$-monopoles. I : Regularity,
Uhlenbeck compactness, and transversality}, Journal of Differential
Geometry (1998)

[FU] Freed D. S. ;  Uhlenbeck, K.:
{\it Instantons and Four-Manifolds.}
Springer-Verlag (1984) 

[GHL] Gallot, S.; Hulin, D.; Lafontaine, J.: {\it Riemannian Geometry}, 
Springer Verlag
(1987)

[HH] Hirzebruch, F.; Hopf,  H.: {\it Felder von Fl\"achenelementen in
4-dimensiona\-len 4-Mannigfaltigkeiten}, Math. Ann. 136 (1958)

[H] Hitchin, N.: {\it  Harmonic spinors}, Adv. in Math. 14, 1-55 (1974)

[Ke] Kelley, J.: {\it General Topology}, Springer Verlag (1955)
 
[K] Kobayashi, S.: {\it Differential geometry of complex vector bundles}, 
Princeton University Press  1987 

[KM] Kronheimer, P.; Mrowka, T.: {\it The genus of embedded surfaces in 
the projective plane}, Math. Res. Letters 1, 797-808 (1994)

[LL] Li, T.; Liu, A.: {\it General wall crossing formula}, Math. Res. Lett. 2,    
  797-810  (1995).
 
[LM] Labastida, J. M. F.;  Marino, M.: {\it Non-abelian monopoles on
 four manifolds}, Preprint,
Departamento de Fisica de Particulas,   Santiago de Compostela, April
 (1995) 

[La]  Larsen, R.: {\it Functional analysis, an introduction}, Marcel Dekker, Inc., 
New York, 1973

[LMi] Lawson, H. B. Jr.; Michelson, M. L.: {\it Spin Geometry}, Princeton 
University Press, New
Jersey, 1989 

[LT] L\"ubke, M.; Teleman, A.: {\it The Kobayashi-Hitchin 
correspondence},  
 World Scientific Publishing Co.  1995

[M]  Miyajima, K.: {\it Kuranishi families of
vector bundles and algebraic description of
the moduli space of Einstein-Hermitian
connections},   Publ. R.I.M.S. Kyoto Univ.  25,
  301-320 (1989) 

[OST]  Okonek, Ch.; Schmitt, A.; Teleman, A.: {\it Master spaces for stable pairs},
 Preprint,
alg-geom/9607015
  
[OT1]   Okonek, Ch.; Teleman, A.: {\it The Coupled Seiberg-Witten 
Equations, Vortices, and Moduli Spaces of Stable Pairs},    Int. J. Math.
Vol. 6, No. 6, 893-910 (1995)

[OT2] Okonek, Ch.; Teleman, A.: {\it Les invariants de Seiberg-Witten 
et la conjecture de Van De  Ven}, Comptes Rendus Acad. Sci. Paris, t.
321,  S\'erie I, 457-461 (1995) 

[OT3] Okonek, Ch.; Teleman, A.: {\it Seiberg-Witten invariants and 
rationality of complex surfaces}, Math. Z., to appear

[OT4] Okonek, Ch.; Teleman, A.: {\it Quaternionic monopoles},  Comptes 
Rendus Acad. Sci. Paris, t. 321, S\'erie I, 601-606 (1995) 

[OT5]  Ch, Okonek.;  Teleman, A.: {\it Quaternionic monopoles},
Commun.\ Math.\ Phys., Vol 180, Nr. 2, 363-388  (1996) 

[OT6] Ch, Okonek.;  Teleman, A.: {\it Seiberg-Witten invariants for manifolds with $b_+=1$, 
and the universal wall crossing formula},
Int. J. Math. 7 (6), 811-832 (1996) 

[OT7] Ch, Okonek.;  Teleman, A.: {\it Recent developments in Seiberg-Witten Theory 
and Complex Geometry},
 preprint, Z\"urich (1996).
 
[PT1]  Pidstrigach, V.; Tyurin, A.: {\it Invariants of the smooth 
structure of an algebraic surface arising from the Dirac operator},
Russian Acad. Izv. Math., Vol. 40, No. 2, 267-351 (1993)

[PT2]  Pidstrigach, V.; Tyurin, A.: {\it Localisation of the Donaldson 
invariants along the Seiberg-Witten classes}, dg-ga/9507004 (1995)

[Sm] Smale, S.: {\it  An infinite dimensionale version of Sard's theorem}, American Journal of
Math. 87,  861-866 (1965)

[T1] Teleman, A. :{\it Non-abelian Seiberg-Witten theory}, Habilitationsschrift, Universit\"at
Z\"urich, 1996

[T2] Teleman, A. :{\it Non-abelian Seiberg-Witten theory and   stable oriented  pairs},
Preprint, Universit\"at Z\"urich, (alg-geom/9609020) (1996)

[T3] Teleman, A. :{\it Fredholm $L^p$-theory for coupled Dirac
operators}, Preprint (1998), to appear in Comptes Rendus de l'Ac. de
Paris

[T4]  A.  Teleman, {\it Almost virtual $PU(2)$-monopoles},
in preparation.

[W] Witten, E.: {\it Monopoles and four-manifolds}, Math.  Res.
Letters 1,  769-796  (1994)
\vspace{1cm}\\
Author's address : \\
 
Mathematisches Institut, Universit\"at Z\"urich,  Winterthurerstr. 190, \\ 
CH-8057 Z\"urich,  {  e-mail}:  teleman@math.unizh.ch

\end{document}